%% file: computability.tex
\documentclass[11pt]{amsart}
\usepackage{eucal}
\usepackage[all]{xy}
\usepackage{amsfonts,amssymb}
\usepackage{epsfig}
\usepackage[usenames]{color}
\usepackage{float} 
\xyoption{arc}

\newtheorem{theorem}{Theorem}[section]
\newtheorem{lemma}[theorem]{Lemma}

\newtheorem{corollary}[theorem]{Corollary}
\newtheorem{defn}[theorem]{Definition}

\newcommand{\C}{\mathbb C}
\newcommand{\R}{\mathbb R}
\newcommand{\K}{\mathbb K}

\newcommand{\D}{\mathbb D}

\newcommand{\Z}{\mathbb Z}
\newcommand{\N}{\mathbb N}

\newcommand{\nhd}{\mbox{nhd}}

\newcommand{\bracket}[1]{ \{#1\} }

\title{Computability and the growth rate of symplectic homology}
\author{Mark McLean}

\begin{document}

\begin{abstract}
For each $n$ greater than $7$ we explicitly construct a sequence of Stein manifolds
diffeomorphic to complex affine space of dimension $n$
so that there is no algorithm to tell us in general whether a given such Stein manifold
is symplectomorphic to the first one or not.
We prove a similar undecidability result for contact structures on the $2n-1$ dimensional sphere.
We can generalize these results by replacing complex affine space with any smooth affine variety of dimension $n$
and the $2n-1$ dimensional sphere with any smooth affine variety intersected with a sufficiently large sphere.
We prove these theorems by using an invariant called the growth rate of symplectic homology to reduce these problems
to an undecidability result for groups.
\end{abstract}

\maketitle

\bibliographystyle{halpha}


\tableofcontents

\section{Introduction}

A Stein manifold is a closed properly embedded complex submanifold
of $\C^N$. We can put a symplectic structure on such a manifold by restricting the standard
one on $\C^N$ to this submanifold.
Smooth affine varieties over $\C$ are important examples of Stein manifolds.
Stein manifolds have been studied symplectically in 
 \cite{Eliashberg:steintopology},\cite{Eliashberg:symplecticgeometryofplushfns}
and \cite{EliahbergGromov:convexsymplecticmanifolds}.
A good way of describing these manifolds in a symplectic way is in
terms of something called Weinstein handle attaching.
We start with the standard unit ball in $\C^n$
and attach handles of dimension $\leq n$ along certain spheres
and extend the symplectic structures in a particular way over these
handles (see
\cite[Section 2.2]{Cieliebak:handleattach} or \cite{Weinstein:contactsurgery}
or
Sections \ref{section:liouvilldomaindefinition} and \ref{section:handlelefschetz}).
This creates a manifold with boundary called a Stein domain.
If $M$ is a Stein domain then we can form a new manifold $\widehat{M}$
called the {\it completion} of $M$ by attaching a cylindrical end $[1,\infty) \times \partial M$
with a certain symplectic structure on it
(see section \ref{section:liouvilldomaindefinition}).
The completion of a Stein domain is symplectomorphic to a Stein manifold
and such Stein manifolds are called {\it finite type} Stein manifolds.
%
The boundary of a Stein domain has a natural contact structure.
Such contact manifolds are called Stein fillable contact manifolds.
This means we can describe many contact manifolds in terms of Weinstein handle attaching.

One question is the following:
If we have two different Stein manifolds that are described in some explicit way
then is there an algorithm to tell us whether they are symplectomorphic or not?
Similarly we can ask the same question of Stein fillable contact manifolds.
The word {\it explicit} is very important here otherwise these questions are
not very interesting.
For instance we can start with two symplectically (even homotopically) different
Stein manifolds $S_1$ and $S_2$
and then for every group presentation $P$ we define
$S_P$ to be $S_1$ if $G_P$ was trivial and $S_2$ if $G_P$ was non-trivial.
Here $G_P$ is the group described by the presentation $P$.
Then there would be no algorithm telling us if a given $S_P$
is symplectomorphic to $S_1$ (otherwise we could solve the word problem).
In this case the Stein manifolds $S_P$ are not described in an explicit way.
For our purposes explicit will mean building our Stein manifolds
(starting from a single point) using the following operations:
\begin{enumerate}
\item Taking the product of our Stein manifold with some smooth affine variety
defined as the zero set of some explicit polynomials.
\item Attaching Weinstein handles along spheres where we know the exact location
of these spheres.
\end{enumerate}

It turns out that the answer to both these questions is no for the following reason:
It is possible (in complex dimension $3$ or higher) to explicitly construct a Stein manifold $S_P$
for every group presentation $P$ with the property that $S_P$ is symplectomorphic
to $S_{\left<|\right>}$ if and only if $P$ represents a trivial group
($\left<|\right>$ is the empty presentation).
The point is that we start with a Weinstein $0$-handle which is just the unit ball in $\C^n$
and then attach a Weinstein $1$-handle for each generator in $P$
and a $2$-handle for each relation so that the fundamental group
of $S_P$ is $G_P$.
Again if we had an algorithm which inputs $P$ and tells us whether
$S_P$ is symplectomorphic to $S_{\left<|\right>}$ or not then it would
also tell us whether $G_P$ is trivial or not which is impossible.
A similar argument works for Stein fillable contact manifolds.

We can strengthen the above questions
by asking if there is an algorithm telling us when
diffeomorphic Stein manifolds are symplectomorphic to each other or not
and similarly for contact manifolds.
The results in \cite[Section 6]{Seidel:biasedview} show that the answer is no
for both Stein manifolds and for Stein fillable contact manifolds.
Here Stein manifolds $A_P$ are constructed explicitly using Weinstein handle attaching for each group presentation $P$
that are all diffeomorphic to each other and so that
$A_P$ is symplectomorphic to $A_{\left<|\right>}$ if and only if $G_P$ is trivial.
Also for each group presentation $P$, a contact manifold $C_P$ is constructed
with similar properties.
The manifolds $A_P$ are all diffeomorphic to the unit disk cotangent bundle of a sphere
of dimension $6$ or higher with Weinstein $2$-handles attached.
Here the unit disk cotangent bundle is the manifold of cotangent vectors on the sphere
whose length is less than or equal to one with respect to some chosen Riemannian metric.

In this paper we prove the same results as in \cite[Section 6]{Seidel:biasedview}
but all of our Stein manifolds are diffeomorphic to Euclidean space
and all of our contact manifolds are diffeomorphic to a sphere.
\begin{theorem} \label{theorem:mainsymplectic}
Let $n \geq 8$.
For each group presentation $P$ there is a finite type Stein manifold
$S_P$ explicitly constructed using smooth affine varieties and Weinstein handle attaching
which is diffeomorphic to $\R^{2n}$ such that
$S_P$ is symplectomorphic to $S_{\left<|\right>}$ if and only if
$G_P$ is trivial.
In particular there is no algorithm to tell us when $S_P$ is symplectomorphic to
$S_{\left<|\right>}$ or not.
None of these symplectic manifolds are symplectomorphic to $\C^n$.
\end{theorem}

\begin{theorem} \label{theorem:maincontact}
Let $n \geq 8$.
For each group presentation $P$ there is Stein fillable contact structure
$\xi_P$  on $S^{2n-1}$ explicitly constructed using smooth affine varieties and Weinstein handle attaching so that
$\xi_P$ is contactomorphic to $\xi_{\left<|\right>}$ if and only if
$G_P$ is trivial.
In particular there is no algorithm to tell us when $\xi_P$ is contactomorphic to
$\xi_{\left<|\right>}$ or not.
None of these contact structures are contactomorphic to $S^{2n-1}$ with the standard contact structure.
The associated contact hyperplane fields are all homotopic as hyperplane subbundles
to the standard contact structure.
\end{theorem}

Here the standard contact structure is the one on the boundary of the Stein domain
given by the unit ball in $\C^n$.
The contact structure of this unit sphere $S^{2n-1}$
is equal to $TS^{2n-1} \cap JTS^{2n-1}$
($J$ is the standard complex structure on $\C^n$).
We will prove these results for a larger class of manifolds.
We can construct a Stein domain $\overline{A}$ by intersecting a smooth affine variety
$A$ with a very large ball (see section \ref{section:liouvilldomaindefinition}).
We can attach Weinstein handles of dimension strictly less than $n$
(called subcritical handles) to create a new Stein domain
$\widetilde{A}$ called an {\it algebraic Stein domain with subcritical handles attached}.
The completion of such a Stein domain is symplectomorphic to
a Stein manifold which we will call
a {\it smooth affine variety with subcritical handles attached}.
The boundary of such a Stein domain is naturally a contact manifold
which we will call a contact manifold {\it fillable by
an algebraic Stein domain with subcritical handles attached}.

\begin{theorem} \label{theorem:mainsymplecticgeneral}
Theorem \ref{theorem:mainsymplectic} is true if we replace
$\R^{2n}$ with some fixed (independent of $P$)
smooth affine variety with subcritical handles attached
with trivial first Chern class, but $\C^n$ in the statement remains the same.
\end{theorem}
Note that this Theorem is a generalization of Theorem \ref{theorem:mainsymplectic}
because $\R^{2n}$ is a smooth affine variety with subcritical handles attached.

\begin{theorem} \label{theorem:maincontactgeneral}
Let $Q'$ be a $2n-1$ dimensional contact manifold fillable by an algebraic Stein domain
with subcritical handles attached which also has trivial first Chern class.
For each group presentation $P$ there is a contact structure $\xi_P$ on $Q'$
constructed explicitly using smooth affine varieties and Weinstein handle attaching
such that $\xi_P$ is contactomorphic to
$\xi_{\left<|\right>}$ if and only if $G_P$ is trivial.
All these contact structures are homotopic as hyperplane subbundles of $TQ'$
to the original contact structure on $Q'$.
None of the contact structures are contactomorphic to the standard contact $2n-1$ dimensional sphere.
\end{theorem}
The above theorem is a generalization of Theorem \ref{theorem:maincontact}
where we have $Q'$ equal to the unit sphere inside $\C^n$.
We get the following direct corollary of Theorem \ref{theorem:maincontactgeneral}:
\begin{corollary} \label{corollary:atleasttwocontactstructures}
Let $n \geq 8$.
Given a $2n-1$ dimensional contact manifold $Q'$
fillable by an algebraic Stein domain with subcritical handles attached,
there are at least two contact structures on $Q'$
that are not contactomorphic to each other but are homotopic as vector
subbundles of $TQ'$ to the original contact structure on $Q'$.
If $Q'$ is the standard $2n-1$ dimensional contact sphere then there are at least $3$ contact structures.
\end{corollary}
\proof of Corollary \ref{corollary:atleasttwocontactstructures}.
Let $P$ represent a non-trivial group (i.e. $P = \left< x | \right>$, the free group
generated by $x$). We have that $\xi_P$ is not contactomorphic to
$\xi_{\left< | \right>}$ by Theorem \ref{theorem:maincontactgeneral}.
This means that there are two contact structures: $\xi_{\left<|\right>}$ and
$\xi_P$.
If $Q'$ is the standard $2n-1$ dimensional sphere we have that
$\xi_P$ and $\xi_{\left< | \right>}$ are not contactomorphic to the standard
$S^{2n-1}$ dimensional contact sphere, so there are $3$ contact structures:
$\xi_{\left<|\right>}$, $\xi_P$ and the standard contact structure on the sphere.
\qed

Here is a summary of some previous results concerning exotic Stein fillable
contact structures on the sphere.
The main theorem in \cite{Ustilovsky:infinitecontact}
tells us that each sphere $S^{4n+1}$ for $n \geq 1$
has infinitely many contact structures up to contactomorphism
whose contact plane distributions are all homotopic to the standard contact one.
But this theorem does not tell us anything about dimensions
$4n+3$.
We also have by \cite[Corollary 5.4]{Seidel:biasedview}
that $S^{2n-1}$ for $n \geq 4$ has at least two non-contactomorphic Stein fillable
contact structures.
This result actually holds for $n=3$ as well by combining
this work with an example from \cite[Section 3.1]{McLean:thesis}.
It turns out from \cite{Eliashberg:twentyyearssinsmartinet} that there is only
one Stein fillable contact structure on $S^3$
so the lowest possible dimension for exotic Stein fillable
contact structures on the sphere is $5$.
Theorems 
\ref{theorem:mainsymplecticgeneral} and \ref{theorem:maincontactgeneral}
will be proven in Section \ref{section:proofofourcomputabilityresult}.

We will now briefly describe how the proof of Theorems 
\ref{theorem:mainsymplecticgeneral} and \ref{theorem:maincontactgeneral} work.
For simplicity we will sketch the proof of Theorems
\ref{theorem:mainsymplectic} and \ref{theorem:maincontact} instead of Theorems
\ref{theorem:mainsymplecticgeneral} and \ref{theorem:maincontactgeneral}.
The main technical tool used is growth rates.
These are invariants of Stein manifolds up to symplectomorphism
(along with a small amount of additional data) taking
values in $\{-\infty\} \cup [0,\infty]$.
Because every Stein domain gives us a unique Stein manifold by completing it,
it is also an invariant of Stein domains.
Growth rates satisfy the following properties:
\begin{enumerate}
\item the growth rate of a product of Stein manifolds
is the sum of the growth rate of its factors.
\item Cotangent bundles naturally have the structure of a Stein manifold
and if $T^*Q$ has fundamental  group given by the product
of at least $3$ non-trivial groups then the growth rate
is in fact infinite.
\item If the boundary of a Stein domain is
fillable by
an algebraic Stein domain with subcritical handles attached
then its growth rate is finite.
\item If we attach subcritical Weinstein handles to a Stein domain then
its growth rate does not change.
\end{enumerate}
These properties are stated in more detail in Section \ref{section:growthrateproperties}.
Let $n \geq 8$.
Novikov in the appendix of \cite{VolodinKuznecovFomenco:algorithmicthreesphere}
constructed (for each group presentation $P$
satisfying some additional mild conditions)
a homology sphere $M_P$ of dimension $n-2$ whose fundamental group is $P$.
We consider the Stein domain $D^* M_P$
(the set of cotangent vectors of length $\leq 1$ with respect to some fixed metric).
We can attach $2$ and $3$ dimensional Weinstein handles until
the Stein domain $D^* M_P$ has trivial fundamental group and such that
all of its homology groups are trivial except in degree $0$
and $n-2$ where it is equal to $\Z$.
We then take the cross product of $D^* M_P$ with a contractible
algebraic Stein domain $T$ and attach a Weinstein $n-1$ handle
so that the homology in degree $n-2$ is killed.
We let $N^4_P$ be this Stein domain and we define $N_P$
to be equal to $N^4_{P'}$ where $P'$ is the free product of $3$ copies of $P$.
By the $h$-cobordism theorem we have that $N_P$ is diffeomorphic to the closed unit ball.
The algebraic Stein domain $T$ must have growth rate greater than or equal to $0$.
An example of $T$ is in \cite[Theorem 3.1]{McLean:thesis}. We could have used
the example from \cite{SS:rama}, although some extra work would have to be done then.
It turns out that using the properties of growth rates as stated
above we have that the growth rate of $N_P$ is finite
if and only if $P$ represents a trivial group.
Also using these properties we have that if $P$ represents a non-trivial
group then the boundary of $N_P$ is not contactomorphic to the boundary
of $N_{\left<|\right>}$
($\left<|\right>$ is the trivial presentation).
It can be shown that $N_P$ is isotopic to $N_{\left<|\right>}$
if $P$ represents a trivial group.
Putting all of this together we have that the completion of $N_P$ which we write as $\widehat{N}_P$
is symplectomorphic to the Stein manifold $\widehat{N}_{\left<|\right>}$ if and only if $P$
gives us a trivial group.
Similarly the boundary of $N_P$ is contactomorphic to the boundary
of $N_{\left<|\right>}$ if and only if $P$ represents a trivial group.
Because there is no algorithm telling us whether $P$ gives us a trivial
group we get that there is no algorithm telling us if
$\widehat{N}_P$ is symplectomorphic to $\widehat{N}_{\left<|\right>}$ or not
and similarly if $\partial N_P$ is contactomorphic to $\partial N_{\left<|\right>}$ or not.

We need to show that $\widehat{N}_P$ is not symplectomorphic to $\C^n$
and $\partial N_P$ is not contactomorphic to the standard contact structure on $S^{2n-1}$.
Because $\C^n$ is constructed entirely using subcritical handles we have that
its growth rate is $-\infty$ but it turns out that the growth rate of $N_P$
is greater than or equal to $0$.
Hence $\widehat{N}_P$ is not symplectomorphic to $\C^n$.
If the boundary of $N_P$ was contactomorphic to $S^{2n-1}$ with its standard contact structure
then basically by \cite[Corollary 6.5]{Seidel:biasedview} we have that its growth rate
vanishes which is a contradiction.
This proves Theorems \ref{theorem:mainsymplectic} and \ref{theorem:maincontact}.

The results about growth rates might be of independent interest.
For instance if we have growth rate greater than $0$ then the contact
boundary must have at least one Reeb orbit.
In fact if the growth rate is greater than $1$ then the boundary of our Stein
domain has infinitely many Reeb orbits even if the contact form has degenerate orbits.
This is the subject of future work 
\cite{McLean:localfloer}.
A similar result using contact homology will be proven in \cite{HryniewiczMacarini:localcontacthomology}.
There are other useful facts about growth rates that are consequences of this paper.
For instance if we have some symplectomorphism $\phi$ of a Liouville domain $F$
which is the identity on the boundary then we can assign to it a Floer homology group
$HF^*(\phi)$.
If $M$ is a Liouville domain whose boundary has an open book with page $F$
and monodromy $\phi$ then the growth rate of $M$ has an upper bound
given by looking at how fast the rank of $HF^*(\phi^k)$ grows as $k$ gets large.
This is true by combining Theorem \ref{theorem:lefschetzgrowthratebound},
Lemma \ref{lemma:partiallefschetzfibrationexistence},
\cite[Theorem 1.2]{McLean:spectralsequence} and
\cite[Theorem 1.3]{McLean:spectralsequence} (we can only use coefficients
in a field of characteristic $0$).
It should be possible to put an open book on the boundary of
any algebraic Stein domain $M$
with monodromy $\phi$ such that the rank of $HF^*(\phi^k)$
is bounded above by some polynomial of degree at most $\text{dim}_\C (M)+1$.
We do not prove this directly here but it can be proven by
showing an inequality of the form
\[\text{sup}\{|i| \text{ } | HF^i(\phi^k) \neq 0\} \leq Ck\] for some constant
$C$ combined with Theorem \ref{theorem:affinelefschetzupperbound},
\cite[Theorem 1.2]{McLean:spectralsequence} and
\cite[Theorem 1.3]{McLean:spectralsequence}.
This open book is constructed using algebraic Lefschetz fibrations
(defined in Section \ref{section:fillingsofalgebraiclefschetzfibrations}).

The paper is organized as follows:
In section \ref{section:themainargument} we state the main definitions,
construct our Stein domains $N_P$, state the properties of growth
rates precisely (without proving them), and prove the main theorems
stated in this introduction.
We spend the remaining sections defining growth rates precisely,
and proving that growth rates satisfy the properties that we stated
earlier.

{\bf Acknowledgements:}
I would like to thank Paul Seidel and Ivan Smith for useful comments concerning this paper.
The author was partially supported by
NSF grant DMS-1005365.

\section{The main argument} \label{section:themainargument}
\subsection{Liouville domains and handle attaching}
\label{section:liouvilldomaindefinition}

A {\it Liouville domain} is a compact manifold $N$
with boundary and a $1$-form $\theta_N$ satisfying:
\begin{enumerate}
\item $\omega_N := d\theta_N$ is a symplectic form.
\item The $\omega_N$-dual of $\theta_N$ is transverse
to $\partial N$ and pointing outwards.
\end{enumerate}
We will write $X_{\theta_N}$ for the $\omega_N$-dual of $\theta_N$
(i.e. so that $\iota(X_{\theta_N})\omega_N = \theta_N$).
We say that two Liouville domains $N_1$ and $N_2$ are {\it Liouville deformation equivalent}
if there is a diffeomorphism $\Phi : N_1 \rightarrow N_2$
and a smooth family of Liouville domain structures $(N_1,\theta^t_{N_1})$
so that $\theta^0_{N_1}= \theta_{N_1}$ and $\theta^1_{N_1} = \Phi^* \theta_{N_2}$.
Sometimes we have manifolds with corners with $1$-forms
$\theta_N$ satisfying the same properties as above.
We view these as Liouville domains by smoothing the corners slightly.
By flowing $\partial N$ backwards along $X_{\theta_N}$ we have a collar
neighbourhood of $\partial N$ diffeomorphic to $(0,1] \times \partial N$
with $\theta_N = r_N \alpha_N$.
Here $r_N$ parametrizes $(0,1]$ and $\alpha_N$ is the contact form on the boundary
given by $\theta_N |_{\partial N}$.
The completion $\widehat{N}$ is obtained
by gluing $[1,\infty) \times \partial N$ to this collar neighbourhood
and extending $\theta_N$ by $r_N \alpha_N$.
By abuse of notation we will write $\theta_N$ for this $1$-form on $\widehat{N}$.
Two Liouville domains are said to be deformation equivalent
if there is a smooth family of Liouville domains joining them.
If we have two Liouville domains that are deformation equivalent
then their completions are exact symplectomorphic.
An {\it exact symplectomorphism} is a symplectomorphism $\Phi$
between two symplectic manifolds $(M_1,d\theta_1)$,$(M_2,d\theta_2)$
such that $\Phi^* \theta_2 = \theta_1 + df$ for some smooth function $f : M_1 \rightarrow \R$.

We will now describe handle attaching.
Weinstein handles were originally described in \cite{Weinstein:contactsurgery}.
An isotropic sphere inside $\partial N$ is a sphere whose
tangent space lies inside the kernel of the contact form $\alpha_N$.
Such a sphere is called a {\it framed isotropic sphere}
if it has some additional framing data
which we will describe in Section \ref{subsection:attachinghandlelefschetz}.
Given such a sphere, we can attach a handle along it and
extend the Liouville domain structure over this handle in a particular way.
Such a handle is called a {\it Weinstein handle}.
The dimension of this handle has to be less than or equal
to half the dimension of $N$ because any isotropic sphere
must have dimension less than half the dimension of $N$.
If the dimension of the Weinstein handle is less than
half the dimension of $N$ then we call such a handle
a {\it subcritical handle}.
Instead of using isotropic spheres to attach Weinstein handles,
we will use handle attaching triples (HAT's).
This is a triple $(f,\beta,\gamma)$ where $f : S^{k-1} \rightarrow \partial N$
is a smooth embedding where $k$ is less than or equal to half the
dimension of $N$.
Also $\beta$ is a normal framing for $f$ inside $\partial N$
(i.e. a bundle isomorphism $\beta : S^{k-1} \times \R^{2n-k} \rightarrow \nu_f$
where $\nu_f$ is the normal bundle to $f$).
Here $\gamma : S^{k-1} \times \C^n \rightarrow f^* TN$
is a symplectic bundle isomorphism where we
give $\C^n$ the standard symplectic structure.
There is an injective bundle homomorphism
$df : TS^{k-1} \hookrightarrow f^*TN$.
Let $\underline{\R}$ be the trivial $\R$ bundle over $S^{k-1}$.
We also have a bundle morphism
$Df : TS^{k-1} \oplus \underline{\R} \hookrightarrow f^*TN$
given by $df + L$ where $L$ sends the positive unit
vector in $\R$ to an inward pointing vector.
The bundle $TS^{k-1} \oplus \underline{\R}$ has a natural trivialization
$\tau$ where we view $S^{k-1}$ as the unit sphere in $\R^k$ and $\underline{\R}$
as the inward pointing vector field.
We say that $(f,\beta,\gamma)$ is a {\it handle attaching triple} or HAT if the map
$\gamma$ is isotopic
to $(Df \circ \tau) \oplus \beta$ through real bundle trivializations of $f^* TN$.
An isotopy of HAT's is a smooth family
of HAT's $(f_t,\beta_t,\gamma_t)$.
For any HAT $(f,\beta,\gamma)$, there is
a framed isotropic sphere $\iota : S^{k-1} \hookrightarrow \partial N$
isotopic to the map $f$ and such that $\iota$ is $C^0$ close to $f$
(this is due to an h-principle, see \cite[Sections
2.1, 2.2, 2.3]{Eliashberg:steintopology}).
So if we wish to attach a Weinstein handle along a sphere
and we wish that sphere to be in a particular homotopy
class of spheres then all we need to do is find a HAT
$(f,\beta,\gamma)$ such that $f$ is in this homotopy class.

An important class of Liouville domains are called Stein domains.
These are constructed as follows:
Suppose we have some complex manifold $A$ with complex structure
$J$. Let $\phi : A \rightarrow \R$ be an exhausting
(i.e. proper and bounded from below) function such that
$-dd^c \phi(X,JX) > 0$ for all non-zero vectors $X$.
Here $d^c\phi(X) := d\phi(JX)$.
Such a function is called an {\it exhausting plurisubharmonic function}.
If a complex manifold admits such a function then it is called
a {\it Stein manifold}.
Let $c$ be a regular value of $\phi$.
Then the compact manifold $\phi^{-1}(-\infty,c]$
with $1$-form $-d^c\phi$ has the structure of a Liouville domain.
We call such a Liouville domain a {\it Stein domain}.
All such Liouville domains can be constructed using Weinstein handle
attaching.
This implies that they are homotopic to a cell complex of dimension
$\leq n$ where $n$ is the complex dimension.
Also we have that any Liouville domain constructed using Weinstein handles
is deformation equivalent to a Stein domain (see \cite{Eliashberg:steintopology}).
Any smooth affine variety $A$ has the structure
of a Stein manifold.
If $\iota : A \hookrightarrow  \C^N$ is a natural embedding
coming from its defining polynomials then
$\iota^* \sum_{i=1}^N |z_i|^2$ is our plurisubharmonic function.
This means that the natural symplectic structure on $A$
is the one induced by the standard one in $\C^N$.
This symplectic structure is unique up to biholomorphism
(see \cite{EliahbergGromov:convexsymplecticmanifolds}).
From \cite[Lemma 2.1]{McLean:affinegrowth} we have that
$A$ is symplectomorphic to the completion of some Stein domain
$\overline{A}$ obtained by intersecting $A$ with a very large closed ball in $\C^N$.
Such a Stein domain is called an {\it algebraic Stein domain}.
We have that if $A_1$ and $A_2$ are isomorphic smooth affine varieties
then $\overline{A}_1$ and $\overline{A}_2$ are Liouville deformation equivalent.

\subsection{Brief description of growth rates and its properties}
\label{section:growthrateproperties}

Let $M$ be a Liouville domain.
We choose an almost complex structure $J$ on $M$ compatible with the symplectic form
and a trivialization $\tau$ of the top exterior power of the
$J$ complex bundle $TM$. We also choose a class $b \in H^2(M,\Z / 2\Z)$.
From this data we can define the growth rate
$\Gamma(M,\tau,b)$ which is an invariant of $\widehat{M}$
up to symplectomorphisms preserving the class $b$ and our trivialization $\tau$
(up to homotopy).
We will suppress the notation $b$ and $\tau$ and just write
$\Gamma(M)$ when the context is clear.

Growth rates satisfy a few important properties which we will now state.
If we have some Riemannian manifold $Q$ then we can define
its unit disk bundle $D^* Q$. 
This is the manifold of covectors whose length is less than or equal to $1$.
It is a Liouville domain with Liouville form $\sum_i q_i dp_i$
where the $q_i$ are position coordinates and $p_i$ are the momentum coordinates.
In fact it is a Stein domain (see \cite{CieliebakEliashberg:symplecticgeomofsteinmflds}).
This has a trivialization of the canonical bundle induced by the natural Lagrangian fibration
structure of $D^* Q$.
We choose our class $b$ to be the second Stiefel Whitney class of $Q$
which we pull back to $D^*Q$ under the projection to $Q$.
For any finitely generated group $G$ we can define the following growth rate
$\Gamma^{\text{cong}}(G)$:
Choose generators $g_1,\cdots,g_k$ of $G$.
Let $f(x)$ the number of conjugacy classes $[g]$ of elements
$g$ such that $g$ can be expressed as a product of at most $x$
generators $g_1,\cdots,g_k$.
We define 
$\Gamma^{\text{cong}}(G)$
to be $\varlimsup_{x} \frac{\log{f(x)}}{\log{x}}$.
\begin{theorem} \label{theorem:cotangentgrowthratebound}
\cite[Lemma 4.15]{McLean:affinegrowth}
$\Gamma(D^*Q) \geq \Gamma^{\text{cong}}(\pi_1(Q))$.
\end{theorem}
We now need a theorem relating growth rates to products.
Suppose we have two Liouville domains $N$ and $N'$.
Then we have the product $\widehat{N} \times \widehat{N}'$.
This has $N \times N'$ as a submanifold with corners.
We can smooth these corners slightly to create a new Liouville
subdomain $N''$ whose completion $\widehat{N}''$ is symplectomorphic to
$\widehat{N} \times \widehat{N}'$.
Hence we can define $\Gamma(\widehat{N} \times \widehat{N'},(\tau,b))$
for some choice of trivialization $\tau$ of the canonical bundle
and $b \in H^2(\widehat{N} \times \widehat{N}',\Z / 2\Z)$.
Choose trivializations $\tau,\tau'$ of the canonical bundles
of $\widehat{N},\widehat{N}'$ and also classes
$b \in H^2(\widehat{N},\Z / 2\Z)$, $b' \in H^2(\widehat{N},\Z / 2\Z)$.

\begin{theorem} \label{theorem:growthrateproduct}
We have
\[\Gamma(\widehat{N} \times \widehat{N'},(\tau \oplus \tau',b \otimes b'))
 = \Gamma(\widehat{N},(\tau,b)) + \Gamma(\widehat{N}',(\tau',b'))\]
\end{theorem}
We will prove this in Section \ref{subsection:products}.
This is basically a growth rate version of the main result in \cite{Oancea:kunneth}.
We need to know what the growth rate is for Liouville domains
whose boundary is contactomorphic to the boundary of a smooth affine variety with subcritical handles attached.
\begin{theorem} \label{theorem:subcriticalhandleattachaffinevariety}
Suppose that we attach a series of subcritical handles to an algebraic Stein domain
$\overline{A}$ to create a Liouville domain $N'$ and that
$H^2(N',\Z/2\Z) \rightarrow H^2(\partial N',\Z/2\Z)$ is surjective.
Let $N''$ be any Liouville domain whose boundary is contactomorphic
to $\partial N'$.
If $M$ is a Liouville domain
such that $\widehat{M}$ is
symplectomorphic to $\widehat{N''}$
then $\Gamma(M) \leq \text{dim}_\C A$.
\end{theorem}
We will prove this in Section \ref{section:affinesubcriticalhandles}.
The surjectivity assumption is not needed if our coefficient field $\K$ is of characteristic $2$.
This theorem is true for any choice of $(\tau,b)$ on $M$.

\begin{theorem} \label{theorem:subcriticalhandlelefschetzgrowth}
Let $M$ be a Liouville domain whose boundary supports
an open book whose pages are homotopic to $CW$-complexes of dimension less
than half the dimension of $M$ and let $M'$ be a Liouville domain
with a subcritical handle attached.
Then $\Gamma(M) = \Gamma(M')$.
\end{theorem}
We will prove this in Section \ref{section:attachingsubcriticalhandles}.
We will not define an open book here
(they are defined in Section \ref{subsection:openbookpartiallefschetz}).
From \cite[Theorem 10]{Giroux:openbooks} we have that every contact manifold admits an open book
supporting the contact structure whose pages are homotopic to an $n-1$-dimensional $CW$ complex.
Note that if a Liouville domain is constructed entirely
using subcritical handles then its growth rate is $-\infty$.
This is because it is the empty Liouville domain with subcritical handles attached
and we define the growth rate of the empty Liouville domain to be $-\infty$.

\subsection{Construction of our Liouville domains}
\label{section:constructionofourliouvilledomains}

Our Liouville domains will be constructed using cotangent bundles,
smooth affine varieties and subcritical Weinstein handle attaching.
We need some preliminary lemmas first.
%

\begin{lemma} \label{lemma:normalbundle2dimsphere}
Suppose that $V_1$ and $V_2$ are two trivial vector bundles on $S^2$
and $V_2$ is a subbundle of $V_1$ of codimension $> 2$.
Then $V_1 / V_2$ is also trivial and if we choose any trivializations
\[\tau_1 : \R^{n_1} \rightarrow V_1,\]
\[\tau_2 : \R^{n_2} \rightarrow V_1,\]
then there is a trivialization
\[\tau_3 : \R^{n_1-n_2} \rightarrow V_1/V_2\]
so that $\tau_2 \oplus \tau_3$ is isotopic through trivializations to $\tau_1$.
\end{lemma}
\proof of Lemma \ref{lemma:normalbundle2dimsphere}.
We will first show that $V_1 / V_2$ is trivial.
Because $S^2$ is the union of two disks along their boundary $S^1$ we have that
$V_1 / V_2$ is determined by an element $q$ of $\pi_1( O(n_1 - n_2) ) \cong \Z / 2 \Z$
as $n_1 - n_2 > 2$.
Because $V_1$ is isomorphic to $V_2 \oplus (V_1 / V_2)$ and $V_1$, $V_2$  are trivial we have that
the image of $q$ under the natural map:
\[ \pi_1(O(n_1 - n_2) ) \hookrightarrow\]\[ \pi_1(O(n_1 - n_2) ) \times \pi_1(O(n_2))
\cong \pi_1( O(n_1 - n_2) \times O(n_2) ) \rightarrow \pi_1(O(n_1) ) \]
is zero.
This natural map is an injection which means that $q$ must be trivial and hence $V_1 / V_2$ is trivial.

We now need to find a trivialization $\tau_3$ for $V_1 / V_2$.
Choose any trivialization $\tau_3$ of $V_1 / V_2$.
We have a bundle isomorphism $\iota$ from $V_2 \oplus (V_1 / V_2)$ to $V_1$.
We have that $\tau_1^{-1} \circ \iota \circ (\tau_2 \oplus \tau_3)$
is a section of the trivial bundle $\text{Aut}(\R^{n_1})$ which we view as some map
$\kappa$ from $S^2$ to $O(n_1)$. After possibly conjugating $\tau_3$ by a reflection
we can assume that $S^2$ maps to the connected component of $O(n_1)$ containing the identity element.
By looking at the natural fibration
\[O(k-1) \hookrightarrow O(k) \twoheadrightarrow S^{k-1} \]
coming from the action of $O(k)$
on $S^{k-1}$ we see that $\pi_2(O(k)) = 0$ for $k >2$.
Hence the map $\kappa$ is isotopic to the constant map which implies that
$(\tau_2 \oplus \tau_3)$ is isotopic through trivializations to $\tau_1$.
\qed

\begin{lemma} \label{lemma:trivialcherclass2dimsphere}
Let $M$ be any Liouville domain whose first Chern class is trivial
and such that $M$ has dimension greater than $4$.
Let $f : S^2 \rightarrow M$ be any map of the two sphere into $\partial M$.
Then there exists $\beta$ and $\gamma$ so that $(f,\beta,\gamma)$ is a HAT.
\end{lemma}
\proof of Lemma \ref{lemma:trivialcherclass2dimsphere}.
Because $c_1(TM) = 0$, we have that
$f^* TM$ can be trivialized as a complex vector bundle by a trivialization $\gamma$.
Let $\underline{X}$ be the oriented real line bundle spanned
by the inward pointing vector field along $\partial M$.
Let $\nu_f$ be the normal bundle to $f$ inside $\partial M$.
This is a real vector bundle of dimension greater than $2$.
Because $\underline{X}$ is a trivial vector bundle, we have a trivialization
$\tau$ of $TS^2 \oplus f^* \underline{X}$ where we view $S^2$ as the unit sphere in $\R^3$
and $f^* \underline{X}$ as the inward pointing vector field along this sphere.
Here $TS^2 \oplus f^* \underline{X}$ is a subbundle of $f^* TM$ and
$\nu_f$ is the normal bundle to this subbundle.
By Lemma \ref{lemma:normalbundle2dimsphere} we have a trivialization $\beta$
of $\nu_f$ such that $\tau \oplus \beta$ is isotopic to $\gamma$ through real bundle trivializations.
This means that $(f,\beta,\gamma)$ is a HAT.
\qed

Let $M$ be a Liouville domain.
A {\it trivially framed sphere} is a sphere
$\iota : S^k \hookrightarrow M$ along with a chosen symplectic bundle isomorphism
$\gamma : S^k \times \C^n \rightarrow TM$
and a trivialization of the normal bundle of $S^k$
given by $\beta$.
We write $\underline{\R}$ for the trivial $\R$ line bundle over $S^k$.
Recall from Section \ref{section:liouvilldomaindefinition} that
we have a canonical trivialization $\tau$ of $TS^k \oplus \underline{\R}$.
We require that $\tau \oplus \beta$ is isotopic
to the trivialization $\gamma \oplus -\underline{\R}$ through real bundle trivializations.
\begin{lemma} \label{lemma:movingoursphere}
Let $(f',\beta',\gamma' )$ be a trivially framed sphere
and let $f : S^k \rightarrow M$ be any smooth map which is isotopic through such smooth
maps to $f'$.
Then there exists a trivially framed sphere $(f,\beta,\gamma)$ which is isotopic through
trivially framed spheres to $(f',\beta',\gamma')$.
\end{lemma}
\proof of Lemma \ref{lemma:movingoursphere}.
Our isotopy can be represented by a map from $[0,1] \times S^k$
to $M$.
Because $[0,1] \times S^k$ deformation retracts onto
$\{0\} \times S^k$ we can extend our trivializations $\beta'$
and $\gamma'$ over the whole of this map to $\tilde{\beta}$
and $\tilde{\gamma}$ respectively.
Hence we define $\beta$ and $\gamma$
to be $\tilde{\beta}$ and $\tilde{\gamma}$
restricted to $\{1\} \times S^k$.
\qed

\begin{lemma} \label{lemma:symplecticistropictocontactisotropic}
Let $M$ be a Liouville domain. We write $\alpha_M$ for the contact form $\theta_M|_{\partial M}$
on the boundary.
Let $W$ be a codimension $1$ submanifold of $\partial M$ such that $d\alpha_M|_W$ is a symplectic form.
Suppose that $f : S^k \hookrightarrow W$ is a trivially framed sphere in $W$
then there is a HAT $(f,\beta,\gamma)$ inside $M$.
\end{lemma}
\proof of Lemma \ref{lemma:symplecticistropictocontactisotropic}.
Because $d\alpha_M|_W$ is non-degenerate we have that
the Reeb vector field $R$ is transverse to $W$.
Also every hyperplane in $\partial M$ transverse to $R$ has the property
that $d\alpha_M$ is non-degenerate.
Hence the vector subbundle $TW$ of $T\partial M|_W$ is isotopic to the contact plane distribution $\xi$
through hyperplane bundles where $d\alpha_M$ is non-degenerate.
The sphere $f$ is trivially framed in $W$ so there is an associated
trivialization $\beta'$ of its normal bundle and
a trivialization $\gamma'$ of the symplectic bundle $TW$ such that
$TW \oplus -\underline{\R}$ is isotopic through real bundles to
$\tau \oplus \beta'$.
Here $\tau$ is the natural trivialization of $TS^k \oplus \underline{\R}$ where we view
this bundle as $T\R^{k+1}$ restricted to the unit sphere.
All of this means that we have a trivialization of the symplectic bundle
$f^* \xi$ given by $\gamma''$ such that $\gamma'' \oplus -\underline{\R}$ is isotopic through real bundles to
$\tau \oplus \beta'$.
Let $\underline{R}$ be the real line bundle spanned by $R$.
We define $V_M$ to be the symplectic bundle spanned by $X_{\theta_M}$
and $R$.
This has a canonical symplectic trivialization $\nu$ induced by $X_{\theta_M}$
and $R$.
We define $\beta$ to be $\beta' \oplus \underline{R}$
and $\gamma$ to be $\gamma'' \oplus \nu$.
Our HAT is $(f,\beta,\gamma)$.
The reason why the trivialization $\tau \oplus \beta$
is isotopic to $\gamma$ is because
$\gamma$ splits up as $\gamma'' \oplus \underline{X_{\theta_N}} \oplus \underline{R}$
where $\underline{X_{\theta_M}}$ is the real bundle spanned by $X_{\theta_M}$.
So we identify $\underline{\R}$ with $-\underline{X_{\theta_M}}$, then use the isotopy from
$\gamma'' \oplus -\underline{\R}$
to $\tau \oplus \beta'$ to give us an isotopy
from $\gamma'' \oplus -\underline{\R} \oplus \underline{R}$
to $\tau \oplus \beta' \oplus \underline{R}$.
This is the isotopy we want because
$\gamma$ is exactly the same as the trivialization $\gamma'' \oplus -\underline{\R} \oplus \underline{R}$
and $\beta$ is the trivialization $\beta' \oplus \underline{R}$.
\qed

\begin{lemma} \label{lemma:steindomainhomotopy}
Let $M$ be a Stein domain of dimension $2n > 4$.
The map $\pi_1(\partial M) \rightarrow \pi_1(M)$ is an isomorphism
and so is $H_i(\partial M) \rightarrow H_i(M)$ for $i < n-1$.
\end{lemma}
\proof of Lemma \ref{lemma:steindomainhomotopy}.
The Stein domain $M$ admits a plurisubharmonic Morse function $\rho$
where $\partial M$ is a regular level set.
The index of all its critical points is $\leq n$.
This means that $M$ is homotopic to $\partial M$
with cells of dimension $\geq n$ attached.
Attaching a cell of dimension $\geq n > 2$ does not change $\pi_1$
or $H_i$ for $i < n-1$. This gives us our result.
\qed

{\bf Construction:}

From now on all of our manifolds are assumed to be oriented unless stated otherwise.
If we have some finite group presentation
\[P := \left< g_1,\cdots,g_k | r_1,\cdots,r_l \right>\]
were $g_i$ are generators and $r_i$ are relations
then we write $G_P$ for the associated group.
We write $\left<|\right>$ for the empty presentation.
We will let $n \geq 8$.
Novikov in the appendix of \cite{VolodinKuznecovFomenco:algorithmicthreesphere}
(see also \cite[Chapter 2]{NabutovskyWeinberger:algorithmichomeomorphism}
and \cite{Nabutovsky:einstein}) constructed for each group presentation
$P$ with $H_1(P)= H_2(P) = 0$ an $n-1$ dimensional
compact manifold with boundary $M_P$ such that:
\begin{enumerate}
\item $M_P$ is acyclic.
\item $\pi_1(M_P) = G_P$.
\item The inclusion map $\partial M_P \hookrightarrow M_P$
induces a fundamental group isomorphism.
\item $\partial M_P$ is a homology sphere.
\end{enumerate}
These manifolds are constructed explicitly using handle attaching.
We can also explicitly find loops in $\partial M_P$
corresponding to the generators $g_1,\cdots,g_k$
of $P$.
By the h-cobordism theorem we have that
if $G_P$ is trivial then $M_P$ is diffeomorphic to a closed ball.
We choose some metric on $\partial M_P$.

\begin{lemma} \label{lemma:adding1and2handles}
We can add Weinstein $2$ and $3$ handles to $D^* \partial M_P$
giving us a Liouville domain which we will call $N^2_P$
such that  $H_i(N^2_P) = 0$ for all $i \neq 0$ or $n-2$
and $\pi_1(N^2_P) = 0$. We also have $H_0(N^2_P) = H_{n-2}(N^2_P) = \Z$.
\end{lemma}
\proof of Lemma \ref{lemma:adding1and2handles}.
We first attach $2$-handles to kill $\pi_1$.
Choose a loop $\tilde{g_i}$ in the unit cotangent bundle $S^* \partial M_P$
which is isotopic to the loop representing $g_i$ in $\partial M_P \subset D^* \partial M_P$.
Because all oriented bundles on one dimensional spheres are trivial we
have HAT's corresponding to the loops $\tilde{g_i}$.
This means we can attach Weinstein $2$-handles along each loop $\tilde{g_i}$.
We can also choose the framing of these handles so that the Chern class
of the resulting Liouville domain is trivial.
This is because each loop $\tilde{g_i}$
has a canonical trivialization of $T\partial M$ (as it is uniquely determined by the trivialization of the canonical
bundle of $D^* \partial M$) and so
we want our HAT trivialization $\gamma$ to coincide with this trivialization.
We define $N^1_P$ to be the resulting Liouville domain.

Because we have attached $2$-handles to all these generators we have that
$N^1_P$ is simply connected so its boundary is also simply connected by Lemma \ref{lemma:steindomainhomotopy}.
Hence by Hurewicz we have that the natural map $\pi_2(\partial N^1_P) \rightarrow H_2(\partial N^1_P)$
is an isomorphism.
We have that $\pi_2(\partial N^1_P) = H_2(\partial N^1_P)$ is a free abelian group generated
by $k$ generators.
Choose $k$ embeddings of the $2$-sphere $f_i : S^2 \hookrightarrow \partial M_P$
representing each of these chosen generators.
By Lemma \ref{lemma:trivialcherclass2dimsphere} we give these maps $f_i$ the structure of a HAT
and then attach Weinstein handles along these HATs creating a new Liouville domain
$N^2_P$.
This Liouville domain has trivial homology groups in all degrees
except $0$ and $n-2$ and in these degrees it is isomorphic to $\Z$.
\qed

By \cite[Theorem 3.1]{McLean:thesis} we can find a contractible
smooth affine variety $T$ of complex dimension $2$ which
has non-zero symplectic homology.
This smooth affine variety is called the tom-Dieck Petrie surface.
Let $\overline{T}$ be the associated algebraic Stein domain whose completion is
$T$.
Let $N^3_P$ be the Liouville domain obtained by smoothing the corners
of $\overline{T} \times N^2_P$ slightly.
\begin{lemma} \label{lemma:addingnhandle}
We can attach a Weinstein $n-1$ handle to $N^3_P$
creating a new Stein domain $N^4_P$ which is diffeomorphic to the unit ball.
\end{lemma}
\proof of Lemma \ref{lemma:addingnhandle}.
Consider the manifold $A_P := [0,1] \times \partial M_P$.
This can be viewed as a collar neighbourhood of $\partial M_P$
where we identify $\partial M_P$ with $\{1\} \times \partial M_P$.
If we give $A_P$ the product metric then we have that
$D^* A_P$ is naturally a submanifold of $D^* [0,1] \times D^* \partial M_P$
which in turn is naturally a submanifold of
$D^* [0,1] \times N^2_P$.
Hence we can create a new manifold $X_P$ with corners
which is the union of $D^* M_P$ and $D^* [0,1] \times N^2_P$
along the common submanifold $D^* A_P$.
We have that $X_P$ is contractible because it is acyclic with trivial fundamental group.
The point is that $N^2_P$ is homotopic to $\partial M_P$ with $2$ and $3$
cells attached killing the fundamental group and so when we attach these cells to the boundary
of $M_P$ we get something which is acyclic and homotopic to $X_P$.
Also the fundamental group and all the homology groups of $N^2_P$
are trivial except in degrees $0$ and $n-2$ where they are $\Z$.
Hence by Hurewicz, we can find an $n-2$ dimensional sphere
representing a generator of $H_{n-2}(N^2_P)$.
We view this as a sphere $f : S^{n-2} \rightarrow D^* [0,1] \times N^2_P$
inside $D^* [0,1] \times N^2_P$.
Inside $X_P$ we also have a trivially framed $n-2$ sphere constructed as follows:
choose a very small chart around some point in $M_P$
and let $f_S : S^{n-2} \rightarrow M_P \subset X_P$ be a small sphere
around this point.
The coordinates $q_1,\cdots,q_{n-1}$ (respecting the orientation of $M_P$) around this point give us a framing $\gamma_S$
for the symplectic bundle $T(D^*M_P)$ over this sphere using the coordinates $q_1,p_1,q_2,p_2,\cdots,q_{n-1},p_{n-1}$ where
$p_i$ are the respective momentum coordinates.
We also have a trivialization $T$ for the outward pointing vector field along this sphere
inside $M_P$ and a trivialization $T'$ of the normal bundle of $M_P$ inside $X_P$
given by $p_1,\cdots,p_{n-1}$ hence we get some trivialization $\beta_S$
of the normal bundle of this sphere given by the sum of these trivializations $T$ and $T'$.
After possibly composing $\beta_S$ with a reflection we have that
$(f_S,\beta_S,\gamma_S)$ is our standard trivially framed sphere.
Because $X_P$ is contractible we then get that
that our $n-2$ sphere $f$ is isotopic inside $X_P$ to our
trivially framed isotropic $n-2$ sphere $f_S$.
Hence by Lemma \ref{lemma:movingoursphere}
we have a trivially framed sphere $(f,\beta_f,\gamma_f)$
inside $D^* [0,1] \times N^2_P$.
Because $D^* [0,1]$ is a small contractible subset of $\C$
we have that our sphere $f$ is contained inside $\D \times N^2_P$
where $\D$ is a small disk inside $\C$.

Choose a small codimension $1$ submanifold $D$ of $\partial \overline{T}$
such that $d\alpha_T$ is a symplectic form on $D$ where $\alpha_T$ is the contact
form on $\partial \overline{T}$.
We can assume that $D$ is symplectomorphic to a small symplectic disk inside $\C$.
The boundary of $N^3_P$ is a smoothing of
\[\partial N^2_P \times \overline{T} \cup N^2_P \times \partial \overline{T}.\]
Hence the boundary of $N^3_P$ has a codimension $1$ submanifold
with symplectic form $d\alpha_{N^3_P}$ symplectomorphic to the product
$N^2_P \times D$ (we might need to shrink $N^2_P$ very slightly but this does not matter).
Here $\alpha_{N^3_P}$ is the natural contact form on the boundary of $N^3_P$.
Because $\overline{T}$ is contractible we have that $N^3_P$
is homotopic to $N^2_P$ and hence is simply connected with trivial homology
groups in all degrees except $0$ and $n-2$.
By Lemma \ref{lemma:symplecticistropictocontactisotropic}
there is a HAT $(f',\beta,\gamma)$ representing the generator
of $H_{n-2}(N^3_P)$.
Hence we can attach a Weinstein $n-1$ handle along this HAT
to create a new Stein domain $N^4_P$ which is simply connected and acyclic.
Hence it must be contractible.
This Liouville domain is diffeomorphic to a $2n$ dimensional ball
(See \cite[Page 174]{DimcaChoudary:complexhypersurfaces}, \cite{Ramanujam:affineplane} and
\cite[Proposition 3.2]{Zaidenberg:1998exot} or \cite[Corollary 2.30]{McLean:symhomlef}).
\qed

If we have two group presentations
$P_1:= \left< g^1_1,\cdots,g^1_{k_1} | r^1_1,\cdots,r^1_{l_1} \right>$,
$P_2:= \left< g^2_1,\cdots,g^2_{k_2} | r^2_1,\cdots,r^2_{l_2} \right>$
then we can form their free product $P_1 \ast P_2$ as follows:
\[\left< g^1_1,\cdots,g^1_{k_1},g^2_1,\cdots,g^2_{k_2} | r^1_1,\cdots,r^1_{l_1},
 r^2_1,\cdots,r^2_{l_2} \right>.\]
Our Stein domains $N_P$ are defined to be equal to $N^4_{P \ast P \ast P}$.

\subsection{Proof of our computability results using growth rates}
\label{section:proofofourcomputabilityresult}

The aim of this section is to prove Theorems
\ref{theorem:mainsymplecticgeneral} and \ref{theorem:maincontactgeneral}.
We will prove several lemmas first.
Throughout this section we will be mentioning growth rates of Liouville domains.
We will use coefficients in $\Z / 2\Z$.
This means that the growth rate of a Liouville domain $M$ is independent
of the choice of our class $b \in H^2(M,\Z / 2\Z)$.
But it does depend on the choice of trivialization of the canonical bundle.
We will be calculating the growth rate of our Liouville domains $N_P$.
These are contractible so there is only one choice of trivialization.
We will also be calculating the growth rate of cotangent bundles
$D^* M_P$ and in this case our trivialization is the natural one
induced by the Lagrangian fibration structure
(i.e. we choose an almost complex structure making the fibers totally
real and then use some volume form on these fibers).

We will use exactly the same notation as in Section \ref{section:constructionofourliouvilledomains}.

\begin{lemma} \label{lemma:trivialfundamentalgroup}
$N_P$ is Liouville deformation equivalent to $N_{\left< | \right>}$
if $G_P$ is a trivial group.
\end{lemma}
\proof of lemma \ref{lemma:trivialfundamentalgroup}.
Throughout this proof, $G_P$ is trivial.
From the last section we showed that $N_P$ is constructed in 3 main stages.
First we attach $2$ and $3$ handles to $D^* \partial M_P$ to create $N^2_P$.
Then we cross with a contractible Stein domain $\bar{T}$ and attach an $n-1$ handle to create $N^4_P$.
Then $N_P = N_{P \ast P \ast P}$.
Similarly $N_{\left<|\right>}$ is created in 3 main stages giving
us two Liouville domains $N^2_{\left<|\right>}$
and $N^4_{\left<|\right>}$.
We will prove this Lemma in 2 steps. In the first step
we will show that $N^2_P$ is Liouville deformation equivalent to
$N^2_{\left<|\right>}$. In fact they are both Liouville deformation equivalent to $D^* S^{n-2}$.
In the second step we will show that $N^4_P$ is Liouville deformation
equivalent to $N^4_{\left<|\right>}$ which implies our result.

{\it Step 1}: 
By the $h$-cobordism theorem we have that
our manifold $\partial M_P$ is diffeomorphic to the $n-2$ sphere.
The Liouville domain $N^2_P$ is Liouville deformation equivalent
to $D^* \partial M_P$ for the following reason:
We are attaching Weinstein $2$-handles
along $k$ disjoint contractible loops creating $N^1_P$.
Our manifold $N^1_P$ is then homotopic to an $n-2$ sphere wedged
with $k$ copies of the $2$-sphere.
We then choose some basis for the free abelian group
$\pi_2(N^1_P)$ and attach $3$-handles along spheres corresponding to this basis
to create $N^2_P$.
By handle sliding we can ensure that these $3$-handles are cancelling
handles for our $2$-handles.
Note that we can handle slide through Weinstein handles
basically because this is equivalent to handle sliding through HAT's
(by using a $1$-parameter version of ideas from \cite[Sections 2.1, 2.2, 2.3]{Eliashberg:steintopology}).
By \cite[Lemma 3.6 b]{Eliashberg:symplecticgeometryofplushfns}
we then get that this Liouville domain is in fact deformation equivalent
to $D^* S^{n-2} = D^* \partial M_P$.
Similar reasoning ensures that $N^2_{\left<|\right>}$ is Liouville deformation equivalent
to $D^* \partial M_P$ which implies that $N^2_P$ is Liouville deformation equivalent to $N^2_{\left<|\right>}$.
This Liouville deformation fixes $\partial M_P \cong \partial M_{\left<|\right>}$.

{\it Step 2}:
To create $N^4_P$, we cross $D^* \partial M_P$ by $\overline{T}$, smooth out the corners
and then we attach an $n-1$ handle.
The Liouville domain obtained by smoothing the corners of $D^* \partial M_P \times \overline{T}$
is called $N^3_P$. We have a similar Liouville domain $N^3_{\left<|\right>}$.
We have that $N^3_P$ is Liouville deformation equivalent to $N^3_{\left<|\right>}$
because  $N^2_P$ is Liouville deformation equivalent to $N^2_{\left<|\right>}$.

In the proof of Lemma \ref{lemma:addingnhandle} we created a manifold with corners $X_P$
which was the union of $D^* M_P$ and $D^* [0,1] \times N^2_P$.
We then created a trivially framed sphere $S_1$ inside $D^* [0,1] \times N^2_P \subset X_P$
which was isotopic to some standard sphere denoted by $(f_S,\beta_S,\gamma_S)$ inside $M_P \subset D^* M_P$.
Because $M_P$ is diffeomorphic to a ball, we can actually assume that
$(f_S,\beta_S,\gamma_S)$ is equal to the sphere $\partial M_P$ with standard framings
induced by coordinates parameterizing our ball $M_P$ and the interval $[0,1]$.
Also we can assume that $S_1$ is isotopic inside $D^* [0,1] \times N^2_P$ to this sphere,
so in fact we may as well assume that $S_1$ is equal to $(f_S,\beta_S,\gamma_S)$.
We also have a similar trivially framed sphere $S_2$ inside
$D^* [0,1] \times N^2_{\left<|\right>}$ which we can assume is the sphere
$\partial M_{\left<|\right>} \subset D^* [0,1] \times N^2_{\left<|\right>}$.
Recall that $D^*[0,1] \times N^2_P$ is naturally a submanifold
of the boundary of $\partial N^3_P$ where the symplectic form is $d\alpha_{N^3_P}$
restricted to this submanifold.
Similarly $D^* [0,1] \times N^2_{\left<|\right>}$ is a submanifold of $\partial N^3_{\left<|\right>}$.
We have that $N^3_P$ and $N^3_{\left<|\right>}$ are Liouville deformation equivalent to each other
and we can ensure that this deformation restricted to the submanifold
\[D^* [0,1] \times N^2_P \subset \partial N^3_P\]
is equal to a product deformation from
$D^* [0,1] \times N^2_P$ to $D^* [0,1] \times N^2_{\left<|\right>}$
(i.e. the symplectic structure on $D^* [0,1]$ remains fixed
and we have the Liouville deformation from $N^2_P$ to $N^2_{\left<|\right>}$
on the other factor).
If $B_t$ is this product deformation and $Q_t$ is the Liouville deformatioin from $N^3_P$
to $N^3_{\left<|\right>}$ then we have a smooth family
of trivially framed spheres $L_t$ in $B_t$ joining $S_1$ and $S_2$ and hence
by Lemma \ref{lemma:symplecticistropictocontactisotropic} this gives us a smooth family of HATs
$A_t$ on $Q_t$.
If we attach a Weinstein handle along $A_0$ in $Q_0 = N^3_P$ then we get $N^4_P$
and similarly we get $N^4_{\left<|\right>}$ by attaching a Weinstein handle along $A_1$.
Hence we get that $N^4_P$ is Liouville deformation equivalent to $N^4_{\left<|\right>}$.
Hence $N_P = N^4_{P \ast P \ast P}$ is Liouville deformation equivalent to
$N_{\left<|\right>} = N^4_{\left<|\right> \ast \left<|\right> \ast \left<|\right>}$.
\qed

\begin{lemma} \label{lemma:cotangentspherealgebraicstein}
$D^* S^{n-1}$ is Liouville deformation equivalent to an algebraic Stein domain.
\end{lemma}
\proof of Lemma \ref{lemma:cotangentspherealgebraicstein}.
Consider the smooth affine variety $V$ given by
\[\left\{\sum_{i=0}^{n-2} z_i^2 = 1 \subset \C^n\right\}.\]
We let $(x_j + i y_j)_{j=0}^n$ be coordinates for $\C^n$.
The equation for our complex hypersurface then becomes:
\[ \sum_{j=1}^n (x_j^2 - y_j^2) = 1, \sum_{j=1}^n x_jy_j = 0.\]
We write $x$ for the vector $(x_j)$ and $y$ for the vector $(y_j)$.
We will show that the Liouville domain obtained by intersecting this complex hypersurface
with a large ball is Liouville deformation equivalent to $D^* S^{n-1}$.
Consider the symplectic manifold $\R^n \times \R^n$
with coordinates $u = (u_j)_{j=1}^n$,$v = (v_j)_{j=1}^n$
and symplectic form $\sum_{j=1}^n du_j \wedge dv_j$.
We will view $T^* S^{n-1}$
as a symplectic submanifold of $\R^n \times \R^n$
by the equations
\[|v| = 1, u.v = 0\]
where $.$ is the standard inner product with respect to these coordinates
and $|.|$ is the standard Euclidean norm.
There is a symplectomorphism $\Phi$
from $V$ (with the standard symplectic form)
to $T^* S^{n-1} \subset \R^n \times \R^n$ given by $v = x / |x|$ and $u = -y|x|$
(see the proof of \cite[Lemma 1.10]{Seidel:longexactsequence} or
\cite[Exercise 6.20(i)]{McduffSalamon:sympbook}).

We have $D^* S^{n-1}$ is Liouville deformation equivalent to the
Liouville domain obtained by intersecting
$T^* S^{n-1}$ inside $\R^n \times \R^n = T^* \R^n$ with
the set $|u| \leq C$ for any $C \geq 0$ and by using the Liouville form
$\sum_{j=1}^n u_j dv_j$.
Another way of thinking about this set is as the subset of
$T^* \R^n$ given by covectors of length $\leq C$ on $S^{n-1} \subset \R^n$
which vanish on the normal vectors to $S^{n-1}$
and where we restrict the standard Liouville form to this subset.

The fibers of $T^* S^{n-1}$ inside the complex hypersurface turn out to be where this variety $V$
intersects the region where $x / |x|$ is constant
where $x$ is the vector $(x_j)$ and $|.|$ is the standard Euclidean norm.
If we intersect these fibers with a sphere
$|x|^2 + |y|^2 = C$ where $C \geq 1$ then they are diffeomorphic to
an $n-1$ dimensional linear hypersurface in $\R^n$ (spanned by the $y$ coordinates)
intersected with a ball of radius $\sqrt{\frac{1}{2}(C-1)}$.
Each of these spheres intersects our variety $V$ transversely
if $C > 1$.
Hence each fiber intersected with $|x|^2 + |y|^2 \leq C$
is diffeomorphic to a ball.
The set $V \cap \{|x|^2 + |y|^2 = C\}$ ($C \geq 1$) is sent under our symplectomorphism
$\Phi$ to the set
\[|u|^2 = \frac{1}{4}(C+1)(C-1).\]
Hence the set $V \cap \{|x|^2 + |y|^2 \leq C\}$ ($C \geq 1$)
is sent to the set
\[|u|^2 \leq \frac{1}{4}(C+1)(C-1)\]
which is in fact Liouville domain deformation equivalent to $D^* S^{n-1}$.
We have that for $C$ large enough that $V \cap \{|x|^2 + |y|^2 \leq C\}$
is a Liouville domain with Liouville form $\theta_V := \sum_{j=1}^n \frac{1}{2}r_j^2 d\vartheta_j$
where $(r_j,\vartheta_j)$ are polar coordinates for the $j$'th $\C$ factor
of $\C^n$ (see \cite[Lemma 2.1]{McLean:affinegrowth}).
This is our algebraic Stein domain.
If we look at $\theta_V$ restricted to the Lagrangian $y = 0$,
we have that it is an exact $1$-form because the Lagrangian
$y=0$ inside $\C^n$ is contractible and hence $\theta_V$
restricted to $\{y=0\} \cap V$ must also be exact.
This implies that our symplectomorphism $\Phi$ is exact.
Putting all of this together we get that
$\Phi$ is an exact symplectomorphism sending the Liouville domain
$V \cap \{|x|^2 + |y|^2\} \leq C\}$ to
the Liouville domain \[|u|^2 = \frac{1}{4}(C+1)(C-1).\]
This implies that these Liouville domains are Liouville deformation equivalent because
we can just use a linear homotopy between $\Phi^* \sum_{j=1}^n u_j dv_j$
and $\sum_{j=1}^n \frac{1}{2} r_j^2 d\vartheta_j$.
Hence because \[|u|^2 \leq \frac{1}{4}(C+1)(C-1)\]
is Liouville deformation equivalent to $D^* S^{n-1}$
we also get that the algebraic Stein domain $V \cap \{|x|^2 + |y|^2\} \leq C\}$
is Liouville deformation equivalent to $D^* S^{n-1}$.
Hence $D^* S^{n-1}$ is Liouville deformation equivalent to an algebraic Stein domain.
\qed

\begin{lemma} \label{lemma:productofalgebraicsteindomains}
Let $\overline{A},\overline{B}$ be two algebraic Stein domains.
Then the Liouville domain obtained by smoothing the corners
of $\overline{A} \times \overline{B}$ is Liouville deformation equivalent to an algebraic Stein domain.
\end{lemma}
\proof of Lemma \ref{lemma:productofalgebraicsteindomains}.
The Liouville domains $\overline{A}$ and $\overline{B}$
are obtained by intersecting smooth affine varieties $A \subset \C^n$ and $B \subset \C^m$
with large balls.
We will assume (see \cite[Lemma 2.1]{McLean:affinegrowth})
that there is a $C \geq 0$ such that for every ball of radius
$\geq C$ inside $\C^n$ centered at $0$ intersects $A$ transversally
and similarly with $B \subset \C^m$.
We will also assume that
all balls of radius $\geq$
$C$ inside $\C^n \times \C^m$ also must intersect $A \times B$
inside $\C^n \times \C^m$ transversally.
We can also assume that the $\omega_A$ dual of
$\sum_{j=1}^n \frac{1}{2} r_j^2 d\vartheta_j|_A$ is transverse to
the boundary of these balls of radius $\geq C$ and pointing outwards
and similarly for $B$ and $A \times  B \subset \C^n \times \C^m$.

Now consider $\C^n \times \C^m$.
We define $(r_j,\vartheta_j)$ to be polar coordinates
for the $j$'th $\C$ factor of this product.
So $(r_j,\vartheta_j)$ are polar coordinates in the $\C^m$ factor if $j > n$.
We define $\theta_{A \times B}$ to be equal to $\sum_{j=1}^{n+m} \frac{1}{2} r_j^2 d\vartheta_j$
restricted to $A \times B$.
This is equal to $\theta_A + \theta_B$ on the product $A \times B$.
Let $X_{\theta_{A \times B}}$ be the $d\theta_{A \times B}$-dual
of $\theta_{A \times B}$.
Hence $X_{\theta_{A \times B}} = X_{\theta_A} + X_{\theta_B}$.
We have that $f = \sum_{j=1}^{n+m} r_j^2$ (viewed as a function on $A \times B$)
satisfies $df(X_{\theta_{A \times B}}) > 0$ for $f \geq C^2$.
Let $f_A := \sum_{j=1}^n r_j^2$ and $f_B := \sum_{j=n+1}^{n+m} r_j^2$.
We have $X_{\theta_{A \times B}}$ is transverse to the boundary of
$\{f_A \leq c^2\} \cap \{f_B \leq c^2\}$ and pointing outwards
for all $c \geq C$.
Hence $\{f_A \leq c^2\} \cap \{f_B \leq c^2\}$ is a Liouville domain
if we smooth its corners.
Let $V$ be this smoothed Liouville domain.
We have that $V$ is a codimension $0$ exact symplectic submanifold
of the Liouville domain $\{f \leq C'\}$ for some large $C'$.
Also the Liouville form on $V$ and $\{f \leq C'\}$ is $\theta_{A \times B}$
and the associated Liouville vector field $X_{\theta_{A \times B}}$
satisfies $f(X_{\theta_{A \times B}}) > 0$ on the closure of
$\{f \leq C'\} \setminus V$.
This means we can deform $V$ through Liouville domains to $\{f \leq C'\}$
because we can smoothly deform the boundary of $V$ while keeping it transverse to
$X_{\theta_{A \times B}}$ until it becomes $\{f = C'\}$
(this can be done if we flow it along $g X_{\theta_{A \times B}}$
where $g>0$ is some function such that any point in $\partial V$
gets flowed along $g X_{\theta_{A \times B}}$ for time $1$ to $\{f = C'\}$).
Hence $V$ is Liouville deformation equivalent to an algebraic Stein domain
$\{f \leq C'\}$.
\qed

\begin{lemma} \label{lemma:trivialgroupaffinevarietywithhandlesattached}
Let $G_P$ be a trivial group.
Then $N_P$ is Liouville deformation equivalent to an algebraic Stein domain
with subcritical Weinstein handles attached.
\end{lemma}
\proof of Lemma \ref{lemma:trivialgroupaffinevarietywithhandlesattached}.
First of all we know that $N_P$ is Liouville deformation equivalent
to $N_{\left<|\right>}$ by Lemma \ref{lemma:trivialfundamentalgroup}
so we will now assume that $P = \left<|\right>$.
If we have some isotopy of Liouville domains $M_t$ ($t \in [0,1]$)
and some sequence of subcritical handles attached to $M_0$ (creating $\widetilde{M}_0$)
then basically by Gray's stability theorem we can attach
a smooth family of Weinstein handles to $M_t$ starting with the original
subcritical handles on $M_0$.
This means that we can attach subcritical handles to $M_1$
creating a Liouville domain $\widetilde{M}_1$
which is Liouville deformation equivalent to $\widetilde{M}_0$.
So if $M_1$ is an algebraic Stein domain then
$\widetilde{M}_0$ is isotopic to an algebraic Stein domain with subcritical
handles attached.

We have that $N_P$ is equal to a smoothing of
$D^* S^{n-2} \times \overline{T}$ with an $n-1$ dimensional Weinstein handle attached
and hence a subcritical handle attached
(by looking at the proof of Lemma \ref{lemma:trivialfundamentalgroup}).
This smoothing is Liouville deformation equivalent to an algebraic Stein domain
by Lemmas \ref{lemma:cotangentspherealgebraicstein} and
\ref{lemma:productofalgebraicsteindomains}.
Hence by the previous discussion we have that $N_P$
is Liouville deformation equivalent to an algebraic Stein domain
with subcritical Weinstein handles attached.
\qed

\begin{lemma} \label{lemma:differentcompletions}
\[\Gamma(\widehat{N}_P) \geq 0\]
for all $P$ and
\[\Gamma(\widehat{N}_P) < \infty\]
if and only if $G_P$ is the trivial group.
\end{lemma}
\proof of Theorem \ref{lemma:differentcompletions}.
We will first show that 
$\Gamma(\widehat{N}_P) \geq 0.$
Again we use the notation from Section \ref{section:constructionofourliouvilledomains}.
We have that every contact manifold admits an open book by \cite[Theorem 10]{Giroux:openbooks}.
Hence by Lemma \ref{theorem:subcriticalhandlelefschetzgrowth}
we have that $\Gamma(D^* \partial M_P) = \Gamma(N^2_P)$.
This is because $N^2_P$ is equal to $D^* \partial M_P$ with subcritical
Weinstein handles attached.
We have $\Gamma(D^* \partial M_P) \geq 0$ by Theorem \ref{theorem:cotangentgrowthratebound}.
By \cite[Theorem 3.1]{McLean:thesis} we know that
$SH_*(\overline{T}) \neq 0$ and hence $\Gamma(\overline{T}) \geq 0$.
This means by Theorem \ref{theorem:growthrateproduct} we have that
$\Gamma(N^3_P) \geq 0$ because it is a smoothing of the product
$N^2_P \times \overline{T}$ and so \[\Gamma(N^3_P) = \Gamma(D^* M_P) + \Gamma(\overline{T}).\]
Finally $N^4_P$ is equal to $N^3_P$ with a subcritical handle attached
and also $N^3_P$ is a Stein domain which means that 
$\Gamma(N^4_P) = \Gamma(N^3_P) \geq 0$ by
\cite[Theorem 10]{Giroux:openbooks}
and Lemma \ref{theorem:subcriticalhandlelefschetzgrowth}.
Hence $\Gamma(N_P) \geq 0$ because $N_P = N^4_{P \ast P \ast P}$.

Suppose now that $G_P$ is non-trivial.
In the previous paragraph we showed that $\Gamma(N_P)$ is equal to
$\Gamma(D^* M_{P \ast P \ast P}) + \Gamma(\overline{T})$.
Because the fundamental group of $M_{P \ast P \ast P}$
is the free product of three non-trivial groups we get
$\Gamma(D^* M_{P \ast P \ast P}) = \infty$
by Theorem \ref{theorem:cotangentgrowthratebound}
and \cite[Lemma 4.16]{McLean:affinegrowth}.
Hence $\Gamma(N_P) = \infty$.

Suppose now that $G_P$ is trivial.
Then $N_P$ is Liouville deformation equivalent to an algebraic Stein domain
with subcritical handles attached by Lemma \ref{lemma:trivialgroupaffinevarietywithhandlesattached}.
This means that $\Gamma(\widehat{N_P}) < \infty$ by
Theorem \ref{theorem:subcriticalhandleattachaffinevariety}
because the boundary of $N_P$ is contactomorphic to the boundary
of an algebraic Stein domain with subcritical handles attached.
\qed

\begin{theorem} \label{theorem:differentcontactboundaries}
The boundary of $N_P$ is not contactomorphic to the boundary of $N_{\left< | \right>}$
if $G_P$ is not trivial.
\end{theorem}

\proof of Theorem \ref{theorem:differentcontactboundaries}.
Let $P$ be the presentation of any group $G_P$ with $H_1(G_P) = H_2(G_P) = 0$.
Suppose that $\partial N_P$ is contactomorphic to $N_{\left<|\right>}$.
By Lemma \ref{lemma:trivialgroupaffinevarietywithhandlesattached},
$N_{\left<|\right>}$ is Liouville deformation equivalent to an algebraic Stein domain
with subcritical handles attached.
This means that $\Gamma(\widehat{N_P}) < \infty$ by
Theorem \ref{theorem:subcriticalhandleattachaffinevariety}.
From Lemma \ref{lemma:differentcompletions}
we have that $G_P$ is trivial if and only if $\Gamma(\widehat{N}_P) < \infty$.
This implies that $G_P$ must be trivial in our case.
Another way of saying this is as follows:
If $G_P$ is non-trivial then the boundary of $N_P$ is not contactomorphic
to the boundary of $N_{\left<|\right>}$.
\qed

Here is a statement of Theorem \ref{theorem:mainsymplecticgeneral}.
{\it Let $n \geq 8$ and $Q$ a smooth affine variety of dimension $n$ with trivial first Chern class.
For each group presentation $P$ there is a finite type Stein manifold
$\widehat{Q}_P$ explicitly constructed using smooth affine varieties and Weinstein handle attaching
which is diffeomorphic to $Q$ such that
$\widehat{Q}_P$ is symplectomorphic to $\widehat{Q}_{\left<|\right>}$ if and only if
$\widehat{Q}_P$ is trivial.
In particular there is no algorithm to tell us when $\widehat{Q}_P$ is symplectomorphic to
$\widehat{Q}_{\left<|\right>}$ or not.
None of these symplectic manifolds are symplectomorphic to $\C^n$.}

\proof of Theorem \ref{theorem:mainsymplecticgeneral}.
Because $Q$ is an algebraic Stein manifold with
subcritical handles attached, we let $\overline{Q}$ be the associated
Stein domain.
By the results in \cite[Appendix]{VolodinKuznecovFomenco:algorithmicthreesphere}
(see also \cite{Nabutovsky:einstein}), for every group presentation $P$ we can explicitly construct
(in an algorithmic way) a group presentation
$P'$ with $H_1(P') = H_2(P') = 0$ and such that $G_P$ is trivial if and only if
$G_{P'}$ is trivial.
We define the {\it end connected sum} of $A \# B$ of two Liouville domains
to be the disjoint union of $A$ and $B$ with a Weinstein $1$-handle
joining $A$ and $B$.
We define $Q_P$ to be the completion of the end connect sum
of $\overline{Q}$ and $N_{P'}$.

Because the growth rate of a disjoint union of Liouville domains
is the maximum of the growth rates of each Liouville domain,
we have that $\Gamma(Q_P)$ is equal to
$\text{max}(\Gamma(\overline{Q}),\Gamma(N_{P'}))$ by Theorem \ref{theorem:subcriticalhandlelefschetzgrowth}.
We have by Theorem  \ref{theorem:subcriticalhandleattachaffinevariety}
that $\Gamma(\overline{Q}) < \infty$.
Hence by Lemma \ref{lemma:differentcompletions}
we have that $\Gamma(Q_P)$ is finite if and only
if $G_{P}$ is trivial.
Combining this with Lemma \ref{lemma:trivialgroupaffinevarietywithhandlesattached}
we have that $\widehat{Q}_P$ is symplectomorphic to
$\widehat{Q}_{\left<|\right>}$ if and only if $G_P$ is trivial.
Here we used the fact that if $A$ and $B$ are Liouville domains with connected boundary
such that $B$ is Liouville deformation equivalent to $B'$
then the end connect sum $A \# B$ is Liouville deformation equivalent to $A \# B'$.
Hence there cannot be an algorithm telling us if $\widehat{Q}_P$
is symplectomorphic to $\widehat{Q}_{\left<|\right>}$ or not.

We also need to show that $\widehat{Q}_P$ is not symplectomorphic to $\C^n$.
Because $\C^n$ is constructed entirely using subcritical handles
(i.e. the zero dimensional subcritical handle) we have that
$\Gamma(\C^{n}) = -\infty$.
But for any group presentation $P$ we have that $\Gamma(\widehat{Q}_P) \geq 0$
which means that $\widehat{Q}_P$ cannot be symplectomorphic to $\C^{n}$
which is $\R^{2n}$ with the standard symplectic structure.
\qed

The following Lemma is really due to Ivan smith:
\begin{lemma} \label{lemma:trivialcontactplanedistibution}
Let $M$ be a Liouville domain which is diffeomorphic to a ball $B$
then there is a diffeomorphism $\Phi$ from $\partial M$
to the sphere such that the push forward via $\Phi$
of the contact distribution is isotopic through hyperplane subvector bundles
of $T\partial B$ to the standard
contact distribution on the boundary of the ball.
\end{lemma}
\proof of Lemma \ref{lemma:trivialcontactplanedistibution}.
Choose an almost complex structure $J$ on $M$
compatible with the symplectic form making the contact plane distribution
on $\partial M$ holomorphic.
Let $R$ be a vector field on $M$.
We assume that $J$ and $R$ satisfy the following properties:
\begin{enumerate}
\item $R$ has one singularity on the interior of $M$
\item $R$ is gradient like
(this can be done because $M$ is diffeomorphic to a ball).
\item The plane distribution $\psi$ spanned by $R$ and $JR$ is symplectically
orthogonal to the contact plane distribution $\xi$ on $\partial M$.
This means that $\psi^\perp$ (the symplectic orthogonal plane distribution)
is equal to $\xi$ on $\partial M$.
\item
On a small neighbourhood around the zero point of $R$,
$R$ is a Liouville vector field which is transverse to some
small codimension $1$ sphere $S$ and pointing outwards.
The contact distribution on $S$ is the standard one.
\item
$\psi^\perp$ restricted to $S$ is the contact structure on $S$.
\end{enumerate}
There is a smooth family of spheres joining this contact sphere $S$
with $\partial M$ transverse to $R$.
By looking at how $\psi^\perp$ behaves as we move these spheres
we get that the contact distribution is isotopic to
the standard one under the diffeomorphism $\Phi$
induced by our smooth family of spheres.
\qed

Here is a statement of Theorem \ref{theorem:maincontactgeneral}.
{\it Let $Q'$ be a $2n-1$ dimensional contact manifold fillable by an algebraic Stein domain
with subcritical handles attached which also has trivial first Chern class.
For each group presentation $P$ there is a contact structure $\xi_P$ on $Q'$
constructed explicitly using smooth affine varieties and Weinstein handle attaching
such that $\xi_P$ is contactomorphic to
$\xi_{\left<|\right>}$ if and only if $G_P$ is trivial.
All these contact structures are homotopic as hyperplane subbundles of $TQ'$
to the original contact structure on $Q'$.
None of the contact structures are contactomorphic to the standard contact $2n-1$ dimensional sphere.
}

\proof of Theorem \ref{theorem:maincontactgeneral}.
By the results in \cite[Appendix]{VolodinKuznecovFomenco:algorithmicthreesphere}
(see also \cite{Nabutovsky:einstein}), for every group presentation $P$ we can explicitly construct
(in an algorithmic way) a group presentation
$P'$ with $H_1(P') = H_2(P') = 0$ and such that $G_P$ is trivial if and only if
$G_{P'}$ is trivial.
We have that $Q'$ is fillable by some
algebraic Stein domain $\widetilde{A}$ with subcritical handles attached.
For each group presentation $P$ we define $\widetilde{A}_P$  to be equal to
the end connect sum (defined in the proof of Theorem \ref{theorem:mainsymplecticgeneral})
of $\widetilde{A}$ and $N_{P'}$.
By Lemma \ref{lemma:trivialgroupaffinevarietywithhandlesattached}
and the fact that the disjoint union of two algebraic Stein domains is
an algebraic Stein domain (because the disjoint union of smooth affine varieties
is a smooth affine variety), we have that
$\widetilde{A}_{\left<|\right>}$ is Liouville deformation equivalent to
an algebraic Stein domain with subcritical handles attached.
Hence by Theorem \ref{theorem:subcriticalhandleattachaffinevariety} we have
that the boundary of $\widetilde{A}_P$ is contactomorphic to
the boundary of $\widetilde{A}_{\left<|\right>}$ if and only if $\Gamma(\widetilde{A}_P) < \infty$.
Basically by the discussion in the proof of Theorem
\ref{theorem:mainsymplecticgeneral}
we have that $\Gamma(\widetilde{A}_P) < \infty$ if and only if $G_{P}$ is trivial.
Hence the boundary of $\widetilde{A}_P$ is contactomorphic to
the boundary of $\widetilde{A}_{\left<|\right>}$ if and only if $G_{P}$ is trivial.
This means that there is no algorithm telling us whether
$\partial \widetilde{A}_P$ is contactomorphic to
the boundary of $\partial \widetilde{A}_{\left<|\right>}$
or not.
Basically by Lemma \ref{lemma:trivialcontactplanedistibution}
we have an (explicit) diffeomorphism $\Phi_P$
from $\partial \widetilde{A}_P$ to $Q'$
such that the pushforward via $\Phi_P$ of the contact distribution
is isotopic through hyperplanes inside $TQ'$ to the original contact distribution
on $Q'$.
Here we used that fact that connect summing a contact manifold $A$ with a contact manifold
diffeomorphic to the sphere whose contact structure is isotopic through hyperplanes to the standard one
gives us a new contact structure on $A$ isotopic through hyperplanes to the old one.
By pushing forward our contact structure via $\Phi_P$
we get a contact structure $\xi_P$ on $Q'$
isotopic through hyperplanes to the original contact structure on $Q'$
and such that there is no algorithm telling us if $\xi_P$ is contactomorphic to $\xi_{\left<|\right>}$ or not.

We now need to show that $\partial \widetilde{A}_P$ is not contactomorphic to $S^{2n-1}$ with the standard
contact structure.
We need an additional fact about growth rates which was not mentioned in Section
\ref{section:growthrateproperties}.
This fact is that if $M$ is a Liouville domain and $\partial M$
is contactomorphic to $S^{2n-1}$ with the standard contact structure
then $\Gamma(M) = -\infty$.
This follows directly from \cite[Corollary 6.5]{Seidel:biasedview}.
This corollary tells us that an invariant called symplectic homology vanishes
but this immediately implies from the definition of growth rates that
$\Gamma(M) = -\infty$.
So for a contradiction suppose that $\partial \widetilde{A}_P$ is contactomorphic to $S^{2n-1}$ then the above statement
says that $\Gamma(\widetilde{A}_P) = -\infty$.
But $\Gamma(\widetilde{A}_P) \geq 0$ because $\Gamma(N_{P'}) \geq 0$
by Lemma \ref{lemma:differentcompletions} which is a contradiction.
Hence $\partial \widetilde{A}_P$ is not contactomorphic to $S^{2n-1}$ for all group presentations $P$.
\qed

\section{Growth rate definition}
\subsection{Symplectic homology} \label{section:symplectichomology}

Let $N$ be a Liouville domain with $c_1 = 0$.
We make some additional choices $\eta := (\tau,b)$ for $N$.
The element $\tau$ is a choice of trivialization of the canonical
bundle of $N$ up to homotopy and $b$ is an element of $H^2(N,\Z / 2\Z)$.
Let $H : \widehat{N} \rightarrow \R$ be a Hamiltonian.
Let $J$ be an almost complex structure on $\widehat{N}$
which is compatible with the symplectic structure.
Let $(S,j)$ be a complex surface possibly with boundary and with a $1$-form
$\gamma$ satisfying $d\gamma \geq 0$.
A map $u : S \rightarrow \widehat{N}$
satisfies the {\it perturbed Cauchy-Riemann equations} with
respect to $(H,J)$ if
$(du - X_H \otimes \gamma)^{0,1} = 0$.
Here
$du - X_H \otimes \gamma$ is a $1$-form on $S$ with values in
the complex vector bundle $\text{Hom}(TS,u^* T\widehat{N})$ where the complex
structure at a point $s \in S$ is induced from $j$ and $J$.
The equation $(du - X_H \otimes \gamma)^{0,1} = 0$ is written explicitly
as
\begin{equation} \label{equation:cauchyriemann}
du - X_H \otimes \gamma + J \circ (du - X_H \otimes \gamma) \circ j = 0.
\end{equation}
Here is a particular example:
Let $S = \R \times S^1 = \C / \Z$. We let $\gamma = dt$
where $t$ parameterizes $S^1 = \R / \Z$.
Then the perturbed Cauchy-Riemann equations become
\[\partial_s u + J \partial_t u = J X_H\]
which is just the Floer equation.
A $1$-periodic orbit of $H$ is a smooth map $o$ from $S^1 = \R / \Z$
to $\widehat{N}$ such that $\frac{do(t)}{dt} = X_H$.
To each orbit we can associate a real number called its action
(which we will define soon).
The pair $(H,J)$ on $\widehat{N}$
is said to satisfy a {\it maximum principle} with respect to an open set
$U^H$ if:
\begin{itemize}
\item
 $U^H$ contains all the $1$-periodic orbits of $H$ of action
greater than some negative constant.
\item
There is a compact set $K' \subset \widehat{N}$
containing $U^H$ such that for any compact complex surface $(S,j)$ with $1$-form $\gamma$
($d\gamma \geq 0$) and map $u : S \rightarrow \widehat{N}$,
satisfying:
\begin{enumerate}
\item $u$ satisfies the perturbed Cauchy-Riemann equations.
\item $u(\partial S) \subset U^H$
\end{enumerate}
we have $u(S) \subset K'$.
\end{itemize}

For any pair $(H,J)$ satisfying the maximum principle with respect to some $U^H$
and pair $\eta := (\tau,b)$ we will define a group $SH_*^{\#}(H,J,\eta)$.
When the context is clear we will suppress the $\eta$ term
and write $SH_*^{\#}(H,J)$.

A $1$-periodic orbit $o : S^1 \rightarrow \widehat{N}$
is non-degenerate if
$D\Phi^1_{X_{H_t}} : T_x\widehat{N} \rightarrow T_x\widehat{N}$
has no eigenvalue equal to $1$ ($\Phi^1_{X_{H_t}}$ is the time $1$ flow of our vector field $X_{H_t}$).
We can perturb $H$ to $H_t$ so that all the orbits
are non-degenerate and so that $H_t = H$ outside a closed subset of $U^H$
(see \cite[Lemma 2.2]{McLean:affinegrowth}).
We can subtract a small constant from $H_t$ so that $H_t < H$.
Let $J_t$ be an almost complex structure so that $J_t = J$
outside a closed subset of $U^H$.
%
We can also perturb $J_t$ so that it is regular and $J_t = J$ outside
some closed subset of $U^H$.
Being regular is a technical condition which will enable us to
define symplectic homology.
We will call the pair $(H_t,J_t)$ an {\it approximating pair} for $(H,J)$
if it satisfies the properties stated above.
We will now define $SH_*^{[c,d]}(H_t,J_t)$ first.

Because we have a trivialization $\tau$ of the canonical bundle
of $\widehat{N}$,
this gives us a canonical
trivialization of the symplectic bundle $T\widehat{N}$
restricted to an orbit $o$ of $H_t$.
Using this trivialization, we can define an index of $o$ called the
{\bf Robbin-Salamon} index (This is equal to the Conley-Zehnder
index taken with negative sign).
We will write $i(o)$ for the index of this orbit $o$.
For a 1-periodic orbit $o$ of $H_t$ we define the {\bf action} $A_{H_t}(o)$:
\[A_{H_t}(o) := -\int_0^1 H_t(o(t))dt -\int_o \theta_N.\]

Choose a coefficient field $\K$.
Let
\[CF_k^d(H_t,J_t,\eta) := \bigoplus_{o} \K \langle o \rangle\]
where we sum over $1$-periodic orbits $o$ of $H_t$ satisfying
$A_{H_t}(o) \leq d$ whose Robbin-Salamon index is $k$.
We write \[CF_k^{(c,d]}(H_t,J_t,\eta) := CF_k^d(H_t,J_t,\eta) / CF_k^c(H_t,J_t,\eta).\]
We need to define a differential for the chain
complex $CF_k^d(H_t,J_t,\eta)$ such that the natural inclusion
maps $CF_k^c(H_t,J_t,\eta) \hookrightarrow CF_k^d(H_t,J_t,\eta)$ for $c<d$ are chain maps.
This makes $CF_k^{(c,d]}(H_t,J_t,\eta)$ into a chain complex as well.
This differential is only well defined for generic $J_t$.

We will now describe the differential
\[\partial : CF_k^d(H_t,J_t,\eta) \rightarrow CF_{k-1}^d(H_t,J_t,\eta).\]
We consider curves
$u : \R \times S^1 \longrightarrow \widehat{N}$ satisfying
the perturbed Cauchy-Riemann equations:
\[ \partial_s u + J_t(u(s,t)) \partial_t u = \nabla^{g_t} H_t\]
where $\nabla^{g_t}$ is the gradient associated to the $S^1$ family of metrics
$g_t := \omega(\cdot,J_t(\cdot))$.
For two $1$ periodic orbits $o_{-},o_{+}$ let
$\widetilde{U}(o_{-},o_{+})$ denote the set of all curves $u$ satisfying
the Cauchy-Riemann equations such that $u(s,\cdot)$ converges
to $o_{\pm}$ as $s \rightarrow \pm \infty$. This has a natural
$\R$ action given by translation in the $s$ coordinate.
Let $U(o_{-},o_{+})$ be equal to $\widetilde{U}(o_{-},o_{+}) / \R$.
For a $C^{\infty}$ generic admissible complex structure
we have that $U(o_{-},o_{+})$ is an $i(o_{-})-i(o_{+}) -1$  dimensional
oriented manifold (see \cite{FHS:transversalitysymplectic}).
Because $(H,J)$ satisfies the maximum principle with respect to $U^H$,
$(H_t,J_t) = (H-\epsilon,J)$ for small $\epsilon>0$
outside a closed subset of $U^H$ and the closure of $U^H$ is compact,
we have that all elements of $U(o_{-},o_{+})$ stay inside
a compact set $K$.
Hence we can use a compactness theorem (see for instance
\cite{BEHWZ:compactnessfieldtheory}) which ensures that
if $i(o_{-}) - i(o_{+}) = 1$, then
$U(o_{-},o_{+})$ is a compact zero dimensional manifold.
The maximum principle for $(H,J)$ is crucial here as it
ensures that $U(o_{-},o_{+})$ is compact.
The class $b \in H^2(N,\Z / 2\Z)$ enables us to orient this manifold
(see \cite[Section 3.1]{Abouzaid:contangentgenerate}).
Let $\# U(o_{-},o_{+})$ denote the number of
positively oriented points of $U(o_{-},o_{+})$
minus the number of negatively oriented points. Then we have a differential:
\[\partial : CF_k^d(H_t,J_t,\eta) \longrightarrow CF_{k-1}^d(H_t,J,\eta) \]
\[\partial \langle o_{-} \rangle := \displaystyle \sum_{i(o_{-}) - i(o_{+})=1 } \# U(o_{-},o_{+}) \langle o_{+} \rangle\]
By analyzing the structure of 1-dimensional moduli spaces, one shows
$\partial^2=0$ and defines
$SH_*(H_t,J_t,\eta)$ as the homology of the above chain complex.
As a $\K$ vector space $CF_k^d(H_t,J_t,\eta)$ is independent of $J_t$ and $b$,
but its boundary operator does depend on $J_t$ and $b$.
We define $SH_*^{(c,d]}(H_t,J_t,\eta)$ as
the homology of the chain complex $CF_*^d(H_t,J_t,\eta) / CF_*^c(H_t,J_t,\eta)$.

Suppose we have two approximating pairs $(H_t^1,J_t^1),(H_t^2,J_t^2)$ of $(H,J)$
such that $H_t^1 < H_t^2$ for all $t$.
Then there is a natural
map:
\[SH_*^{(c,d]}(H_t^1,J_t^1,\eta) \longrightarrow SH_*^{(c,d]}(H_t^2,J_t^2,\eta)\]
This map is called a {\bf continuation map}.
This map is defined from a map $C$ on the chain level as follows:
\[C : CF_k^d(H_t^1,J_t^1,\eta) \longrightarrow CF_k^d(H_t^2,J_t^2,\eta) \]
\[\partial \langle o_{-} \rangle := \displaystyle
\sum_{i(o_{-}) = i(o_{-}) } \# P(o_{-},o_{+}) \langle o_{+} \rangle\]
where $P(o_{-},o_{+})$ is a compact oriented zero dimensional manifold
of solutions of the following equations:
Let $(K^s_t,Y^s_t)$ $(s,t) \in \R \times S^1$ be a smooth family of pairs such that
\begin{enumerate}
\item $(K^s_t,Y^s_t) = (H_t^1,J_t^1)$ for $s \ll 0$.
\item $(K^s_t,Y^s_t) = (H_t^2,J_t^2)$ for $s \gg 0$.
\item $(K^s_t,Y^s_t) = (H-\epsilon_s,J)$ outside some closed subset of $U^H$
where $\epsilon_s >0$ is a smooth family of constants.
\item $K^s_t$ is non-decreasing with respect to $s$.
\end{enumerate}
The set  $P(o_{-},o_{+})$ is the
set of solutions to the parameterized Floer equations
\[ \partial_s u + Y^s_t(u(s,t)) \partial_t u = \nabla^{g_t} K^s_t\]
such that $u(s,\cdot)$ converges
to $o_{\pm}$ as $s \rightarrow \pm \infty$.
For a $C^{\infty}$ generic family $(K^s_t,Y^s_t)$ this is a compact zero dimensional
manifold (if $o_-,o_+$ have the same index).
Again the the class $b \in H^2(N,\Z / 2\Z)$ enables us to orient this manifold.
If we have another family of pairs joining $(H_t^1,J_t^1)$ and $(H_t^2,J_t^2)$
then the continuation map induced by this second family is the same as
the map induced by $(K^s_t,Y^s_t)$.
The composition of two continuation maps is a continuation map.
This means that we can define $SH_*^{\#}(H,J)$
as the direct limit of $SH_*^{[0,\infty)}(H_t,J_t)$
over all approximating pairs $(H_t,J_t)$
with respect to the ordering $<$.

Suppose that we have another pair $(H',J')$ satisfying the maximum principle
with respect to an open set $U^{H'}$ such that 
\begin{enumerate}
\item $U^H \subset U^{H'}$.
\item $H' \geq H$.
\item $(H',J') = (\lambda H + \kappa,J)$ outside a closed subset of $U^{H'}$
for some $\lambda \geq 1$, $\kappa \in \R$.
\end{enumerate}
then we have a natural map
\[SH_*^{\#}(H,J) \rightarrow SH_*^{\#}(H',J').\]
The reason why we have such a map is that if $(H_t,J_t)$
is an approximating pair for $(H,J)$
then we can choose a non-decreasing family $(H_t^s,J_t^s)$
which is equal to $\lambda_s H + \kappa_s$ outside a closed subset of $U^{H'}$
such that $(H_t^s,J_t^s)$ is equal to the approximating pair
$(H_t,J_t)$ for $s \ll 0$ and some approximating pair
$(H'_t,J'_t)$ of $(H',J')$ for $s \gg 0$.
This gives us a morphism
\[SH_*^{\#}(H_t,J_t) \rightarrow SH_*^{\#}(H'_t,J'_t).\]
Because this morphism is induced by a continuation map,
it induces a morphism of directed systems
defining 
$SH_*^{\#}(H,J)$ and $SH_*^{\#}(H',J')$ respectively.
Hence it induces a morphism
\[SH_*^{\#}(H,J) \rightarrow SH_*^{\#}(H',J').\]

If we use orbits of all actions we can define the group
$SH_*^{(-\infty,\infty)}(H,J)$.
We will write $SH_*$ instead of $SH_*^{(-\infty,\infty)}$.
If all of the $1$-periodic orbits of $(H,J)$ have non-negative action
then $SH_*^\#(H,J) = SH_*(H,J)$.
If we wish to stress which coefficient field we are using,
we will write $SH_*^{\#}(H,J,\K)$ if the field is $\K$ for instance.

\subsection{Growth rates} \label{section:growthrates}

In order to define growth rates, we will need some linear algebra first.
Let $(V_x)_{x \in [1,\infty)}$ be a family of vector spaces indexed by $[1,\infty)$.
For each $x_1 \leq x_2$ we will assume that there is a homomorphism
$\phi_{x_1,x_2}$ from $V_{x_1}$ to $V_{x_2}$ with the property that
for all $x_1 \leq x_2 \leq x_3$, $\phi_{x_2,x_3} \circ \phi_{x_1,x_2} = \phi_{x_1,x_3}$
and $\phi_{x_1,x_1} = \text{id}$.
We call such a family of vector spaces a {\it filtered directed system}.
Because these vector spaces form a directed system, we can take the direct
limit $V := \varinjlim_x V_x$.
From now on we will assume that the dimension of $V_x$ is finite dimensional.
For each $x \in [1,\infty)$ there is a natural map:
\[q_x : V_x \rightarrow \varinjlim_x V_x. \]
Let $a : [1,\infty) \rightarrow [0,\infty)$ be a function such that
$a(x)$ is the rank of the image of the above map $q_x$.
We define the growth rate as:
\[\Gamma( (V_x) ) : =\varlimsup_x \frac{\log{a(x)}}{\log{x}} \in \{-\infty\} \cup [0,\infty].\]
If $a(x)$ is $0$ then we just define $\log{a(x)}$ as $-\infty$.
If $a(x)$ was some polynomial of degree $n$ with positive leading coefficient, then the growth rate
would be equal to $n$.
If $a(x)$ was an exponential function with positive exponent, then the growth rate is $\infty$.
The good thing about growth rates is that
if we had some additional vector spaces $(V'_x)_{x \in [1,\infty)}$
such that the associated function $a'(x) := \text{rank}(V'_x \rightarrow \varinjlim_x V'_x)$
satisfies $a'(x) = A a(B x)$ for some constants $A,B>0$ then
$\Gamma(V'_x)  = \Gamma(V_x)$.
The notation we use for filtered directed systems
is usually of the form $(V_x)$ or $(V_*)$, and we will usually write
$V_x$ without brackets if we mean the vector space indexed by $x$.
There is a notion of isomorphism of filtered directed systems.
We do not need to know the exact definition in this section
(it is defined in Section \ref{section:growthratelinearalgebra}).
The only property we need to know is that if two filtered directed
systems are isomorphic then they have the same growth rate
by \cite[Lemma 3.1]{McLean:affinegrowth}.

Let $N$ be a Liouville domain and $\widehat{N}$ its completion.
Let $\theta_N$ be the respective Liouville form.
An {\it $SH_*$ admissible pair} $(H,J)$ on $\widehat{N}$ is a pair satisfying:
\begin{enumerate}
\item For all $\lambda \geq 1$ outside some discrete subset $A^H$,
$(\lambda H,J)$ satisfies the maximum principle
with respect to an open set $U^H_\lambda$.
\item $U^H_{\lambda_1} \subset U^H_{\lambda_2}$ for $\lambda_1 \leq \lambda_2$.
\end{enumerate}
A {\it growth rate admissible pair} $(H,J)$ is an
$SH_*$ admissible pair such that:
\begin{enumerate}
\item (bounded below property)

The Hamiltonian $H$ is greater than or equal to zero,
and there exists a compact set $K$ and a constant $\delta_H > 0$
such that:
$H > \delta_H$ outside $K$.

\item (Liouville vector field property)

There exists an exhausting function $f_H$, and $1$-form $\theta_H$
such that:
\begin{enumerate}
\item $\theta_N - \theta_H$ is exact.
\item There exists a small $\epsilon_H > 0$ such that
$dH(X_{\theta_H}) > 0$ in the region $H^{-1}(0,\epsilon_H]$
where $X_{\theta_H}$ is the $\omega_N$-dual of $\theta_H$.
\item
There is a constant $C$ such that
$df_H(X_{\theta_H}) > 0$
in the region $f_H^{-1}[C,\infty)$ and
$f_H^{-1}(-\infty,C]$ is non-empty and is contained in the interior of $H^{-1}(0)$.
\end{enumerate}
\item (action bound property)

There is a constant $C_H$ such that
the function
$-\theta(X_H) - H$ is bounded above by $C_H$
where $X_H$ is the Hamiltonian vector field associated to $H$.
Here $\theta$ is some $1$-form such that $\theta - \theta_N$
is exact.

\end{enumerate}

So for all $\lambda \notin A_H$
we can define $SH_*^{\#}(H,J)$ as in the previous section.
Also for $\lambda_1 \leq \lambda_2$
we have a natural map (induced by continuation maps)
from $SH_*^{\#}(\lambda_1 H,J)$ to $SH_*^{\#}(\lambda_2 H,J)$.
If $\lambda \in A_H$ then we define $SH_*^{\#}(\lambda H,J)$
as the direct limit of $SH_*^{\#}(\lambda' H,J)$
over all $\lambda' < \lambda$ and not in $A_H$.
If $\lambda_1 = \lambda_2$ then the respective continuation
map is an isomorphism.
Hence we have a filtered directed system
$(SH_*^{\#}(\lambda H,J))$.

From \cite[Corollary 4.3]{McLean:affinegrowth},
we have the following theorem:
\begin{theorem} \label{theorem:filtereddirectedsysteminvariance}
Suppose that $\widehat{N}$ is symplectomorphic to
$\widehat{N'}$ and $(H,J)$ is growth rate admissible
on $\widehat{N}$ and $(H',J')$ is growth rate admissible on
$\widehat{N'}$.
This symplectomorphism must preserve our choice of trivialization of
the canonical bundle and also our choice of $b \in H^2(\widehat{N},\Z / 2\Z)$.
Then the filtered directed system 
$(SH_*^{\#}(\lambda H,J))$ is isomorphic to
$(SH_*^{\#}(\lambda H',J'))$.
\end{theorem}
Hence we have an isomorphism class of filtered directed systems
which is an invariant of $\widehat{N}$ up to symplectomorphism
(preserving our choice of trivialization of the canonical bundle and $b$).
We will write
$(SH_*(\widehat{N},d\theta_N,\lambda))$ for any
filtered directed system in this isomorphism class.
If the context is clear we will just write
$(SH_*(\widehat{N},\lambda))$.

\begin{defn} \label{defn:growthrate}
We define the growth rate $\Gamma(\widehat{N},d\theta_N)$ as:
\[\Gamma(\widehat{N},d\theta_N) := \Gamma(SH_*(\widehat{N},d\theta_N,\lambda)).\]
\end{defn}
Again we suppress $d\theta_N$ from the notation if the context is clear.
Theorem \ref{theorem:filtereddirectedsysteminvariance} combined with the fact that growth rate is an invariant
of filtered directed systems up to isomorphism tells us that $\Gamma(\widehat{N},d\theta_N)$
is an invariant up to symplectomorphism preserving our choice of trivialization of
the canonical bundle and also our choice of $b \in H^2(\widehat{N},\Z / 2\Z)$.

\subsection{A Floer homology group for symplectomorphisms}

We will use coefficients in a field $\K$.
Here we define the Floer cohomology groups $HF_*(\phi,k)$ for each $k \in \N$,
where $\phi : \widehat{F} \rightarrow \widehat{F}$ is a compactly supported symplectomorphism
and where $\widehat{F}$ is the completion of a Liouville domain $(F,\theta_F)$.
We write $\omega_F := d\theta_F$.
The boundary $\partial F$ has a natural contact form $\alpha_F := \theta_F|_{\partial F}$
and $\theta_F = r_F\alpha_F$ on the cylindrical end $[1,\infty) \times \partial F$
where $r_F$ parameterizes $[1,\infty)$.
For simplicity we assume that the first Chern class of the
symplectic manifold $F$ is trivial.
We will assume that $\phi$ is an exact symplectomorphism.
An exact symplectomorphism is a map that satisfies $\phi^* \theta_F = \theta_F + df$
for some function $f : \widehat{F} \rightarrow \R$.
Any compactly supported symplectomorphism is isotopic
through compactly supported symplectomorphisms to an exact
symplectomorphism anyway so this does not really
put any constraint on $\phi$
(see the proof of \cite[Lemma 1.1]{BEE:legendriansurgery}).

By enlarging $F$ we may as well assume that the support of $\phi$ is
contained in $F$.
From \cite[Section 2.1]{McLean:spectralsequence} we have that the mapping
torus $M_\phi$ has a natural contact form $\alpha_\phi$ satisfying:
\begin{enumerate}
\item $d\alpha_\phi$ restricted to each fiber is a symplectic
form. This means we have a connection on this fibration
coming from the line field that is $d\alpha_\phi$ orthogonal
to the fibers.
\item The monodromy map going positively around $S^1$ is Hamiltonian isotopic to $\phi$.
This means that it is equal to $\phi$ composed with a compactly supported Hamiltonian symplectomorphism.
\item \label{item:propertynearinfinity}
Near infinity, the fibration is equal to
the product fibration
\[ [R,\infty) \times \partial F \times S^1 \twoheadrightarrow S^1\]
where $\alpha = d\theta + \theta_F$.
Here $\theta$ is the angle coordinate.
We can enlarge $F$ so that $R = 1$.
\end{enumerate}
The contact plane distribution has a natural symplectic form.
We make this hyperplane distribution into a complex bundle
by putting a compatible almost complex structure on it.
Choose a trivialization $\tau$ of the highest exterior power
of this complex bundle
and a class
$b \in H^2(\widehat{F},\Z / 2\Z)$ which is invariant under the
symplectomorphism.
This class $b$ can be viewed as a class $\widetilde{b} \in H^2(M_\phi,\Z/2\Z)$
which restricts to $b$ on any particular fiber.
The group $HF_*(\phi,k)$
depends on these choices but we suppress them from the notation
when the context is clear.

This group is defined in \cite{McLean:spectralsequence}.
Let $S$ be the surface $(0,\infty) \times S^1$
parameterized by $(r_S,t)$. Here we identify $t \in \R/\Z = S^1$.
From \cite[Section 2.3]{McLean:spectralsequence},
there exists a fibration
$\pi_\phi : W_\phi = (0,\infty) \times M_\phi \rightarrow (0,\infty) \times S^1$ such that:
\begin{enumerate}
\item $\pi_\phi$ splits up as a product $\text{id} \times p_\phi$
where
$p_\phi : M_\phi \twoheadrightarrow S^1$ is a fibration
whose fiber is $\widehat{F}$.
\item 
Its Liouville form $\theta_\phi$ is equal to $r_S\pi_\phi^* dt + \alpha_\phi$
where $\alpha_\phi$ is a  contact form on $M_\phi$.
\item
The monodromy map of $M_\phi$ is the symplectomorphism $\phi$.
\end{enumerate}
Near infinity (in the fiberwise direction), $W_\phi$ looks like $[1,\infty) \times \partial F \times S^1$
where $\pi_\phi$ is the projection map to $S^1$ and $\alpha_\phi = (r_S+1)dt + r_F\alpha_F$.
The coordinate $r_S$ can be viewed as a coordinate on $W_\phi$ by pulling it back via $\pi_\phi$.
The cylindrical coordinate $r_F$ of $\widehat{F}$ can be viewed as a well
defined coordinate in the region $[1,\infty) \times \partial F \times S^1$ inside $W_\phi$
parameterizing $[1,\infty)$.
We will write this region as $\{r_F \geq 1\}$.
Let $\epsilon>0$ be a constant smaller than the smallest
Reeb orbit of $\partial F$.
We can put a trivialization $\tilde{\tau}$ of the canonical
bundle of $W_\phi$ using the trivialization $\tau$
and the fact that the base $S$ has a canonical trivialization.

Let $h : [1,\infty) \rightarrow \R$ be a function with $h(x) = 0$
near $x=1$ and $h = \epsilon r_F$ near infinity.
We can view $h(r_F)$ as a function on $W_\phi$
by extending it by zero in the region where $r_F$ is ill defined as a function.
We say that a Hamiltonian $H : S^1 \times W_\phi \rightarrow \R$ is {\it admissible}
in this context if it is equal to $g(r_S) + h(r_F)$
outside a large compact set where $g : (0,\infty) \rightarrow \R$
is a function satisfying:
\begin{enumerate}
\item $g',g'' \geq 0$.
\item $g'(s)$ is constant for $s$ near $0$ or near $\infty$.
\item $0<g'(s) < 1$ for $s$ near $0$.
\end{enumerate}
The value of $g'$ near infinity is called the {\it slope} of $H$.
%
%
Let $j$ be the complex structure on $S$ where
we identify $S$ with ${\mathbb H} / \Z$ where ${\mathbb H}$ is the
upper half plane in $\C$ and the $\Z$ action is translation.
We also choose an $S^1$ family of almost complex structures $J$
on $W_\phi$ making $\pi_\phi$ $(J,j)$ holomorphic
and such that in the region $\left( [1,\infty) \times \partial F \right) \times S$
it splits up as a product $J_F + j$ where $J_F$
is convex on the cylindrical end $[1,\infty) \times \partial F$.
Here convex means that $dr_F \circ J_F = -\theta_F$.
Then it turns out by maximum principles \cite[Lemma 7.2]{SeidelAbouzaid:viterbo}
and \cite[Lemma 5.2]{McLean:symhomlef}
that $SH_*(H,J)$ is well defined for generic such $(H,J)$.
If we have some subset $A$ of $H^1(W_\phi)$
then we can consider only those $1$-periodic orbits
whose $H^1$ class lies inside $A$.
This is a subgroup $SH_*^A(H,J)$ of $SH_*(H,J)$.
Again this group depends on $\tilde{b}$ and $\tilde{\tau}$
but we suppress this from the notation.
Let $\beta_k$ be the subset of $H^1$ represented
by loops in $W_\phi$ that wrap around the $S^1$
factor of $(0,\infty) \times S^1$ $k$ times after
projecting the loop down by $\pi_\phi$.
We define $HF_*(\phi,k)$ as the direct
limit over all admissible pairs $(H,J)$
with $H|_{\bracket{\pi_\phi^* s < 1}} < 0$ of
$SH_*^{\beta_k}(H,J)$.
Note that this Floer homology group depends on $\tilde{\tau}$
and $\tilde{b}$.
The ordering of this direct limit is the ordering
where $(H_1,J_1)$ is less than $(H_2,J_2)$
if and only if $H_1 < H_2$.
This turns out to be a finite dimensional group as a consequence
of \cite[Lemma 2.9]{McLean:spectralsequence}.
If we have two such exact symplectomorphisms $\phi_1$
and $\phi_2$ that can be joined together by a smooth
family of compactly supported exact symplectomorphisms
then $HF_*(\phi_1,k) = HF_*(\phi_2,k)$ hence
this is an invariant up to isotopy.
Sometimes we just have a symplectomorphism $\phi : F \rightarrow F$
which fixes the boundary $\partial F$.
This has a Floer homology group $HF_*(\phi,k)$
as well as we extend $\phi$ by the identity map
giving us a $C^1$ function which we smooth
to a map $\widehat{\phi} : \widehat{F} \rightarrow \widehat{F}$
and we define $HF_*(\phi,k)$ as $HF_*(\widehat{\phi},k)$.

\section{Products} \label{subsection:products}

The goal of this section is to prove Theorem \ref{theorem:growthrateproduct}.
Let $N$ and $N'$ be Liouville domains.
We have that $N \times N'$ is not a Liouville domain,
but we can smooth the corners slightly so that it becomes a Liouville domain
whose completion is symplectomorphic to $\widehat{N} \times \widehat{N'}$.
The statement of this theorem is:
\[\Gamma(\widehat{N} \times \widehat{N'},(\tau \otimes \tau',b \otimes b'))
 = \Gamma(\widehat{N},(\tau,b)) + \Gamma(\widehat{N}',(\tau',b')).\]
From now on we will suppress our choice of $\tau,\tau',b,b'$ from the notation
unless it is unclear which choices to make.
Theorem \ref{theorem:growthrateproduct} is a consequence of the following Theorem:
Suppose that $(V_\lambda),(V'_\lambda)$ are filtered directed systems.
Then we can form a new filtered directed system
$(V_\lambda \otimes V'_\lambda)$.
The growth rate of $(V_\lambda \otimes V'_\lambda)$
is equal to the sum of the growth rates of $(V_\lambda)$
and $(V'_\lambda)$.
\begin{theorem} \label{theorem:kunnethformula}
The filtered directed system
$(SH_*(\widehat{N} \otimes \widehat{N'},\theta_N + \theta_{N'},\lambda))$
is isomorphic to the tensor product:
$(SH_*(\widehat{N},\theta_N,\lambda) \otimes SH_*(\widehat{N'},\theta_{N'},\lambda))$.
\end{theorem}
This proves Theorem \ref{theorem:growthrateproduct}.

\proof of Theorem  \ref{theorem:kunnethformula}.
Let $(H,J)$ (resp. $(H',J')$) be growth rate admissible for
$\widehat{N}$ (resp. $\widehat{N}'$).
We will now show that the pair $(H+H',J \oplus J')$ is growth rate admissible.

{\it $(H+H',J \oplus J')$ is $SH_*$ admissible:}
Let $A_H$ (resp. $A_{H'}$) be a discrete subset
of $(0,\infty)$ such that $\lambda H$ (resp. $\lambda H'$)
has all of its $1$-periodic orbits inside the relatively
compact open set $U^H_\lambda$ (resp. $U^{H'}_\lambda$)
where $\lambda \in (0,\infty) \setminus A_H$ (resp. $(0,\infty) \setminus A_{H'}$).
This implies that $\lambda(H+H')$
has all of its $1$-periodic orbits inside $U^H_\lambda \times U^{H'}_\lambda$
for $\lambda \in (0,\infty) \setminus (A_H \cup A_{H'})$.
We also have that $(H+H',J+J')$ is $SH_*$-admissible as any solution
$u : S \rightarrow \widehat{N} \times \widehat{N'}$
of the perturbed Cauchy-Riemann equations whose boundary maps
into $U^H_\lambda \times U^{H'}_\lambda$
must be contained in a compact set of the form $K \times K'$.
This is because $u$
projects to a solution of the perturbed Cauchy-Riemann equations for the pair
$(H,J)$ or $(H',J')$ under the projection to $\widehat{N}$
or $\widehat{N'}$ whose boundary maps to $U^H_\lambda$
or $U^{H'}_\lambda$ and hence by the maximum principle
must be contained in a compact set $K \subset \widehat{N}$
or $K' \subset \widehat{N'}$.

{\it $(H+H',J \oplus J')$ satisfies the bounded below property and the action bound property}
because $H$ and $H'$ satisfies the bounded below property and the action bound property.

{\it $(H+H',J \oplus J')$ satisfies the Liouville vector field property:}
Let $V_H,f_H$ (resp. $V_{H'}, f_{H'}$) be the respective
Liouville vector field and function as stated in the Liouville
vector field property for $H$ (resp $H'$).
We also have a constant $C$ so that $N_C := f_H^{-1}(-\infty,C]$
is contained in the interior of $H^{-1}(0)$ and
$V_H(f_H) >0$ for $f_H \geq C$.
We have a similar constant $C'$ for $f_{H'}$.
Write $N'_{C'} := f_{N'}^{-1}(-\infty,C]$.
We have that $N_C \times N'_{C'}$ is contained in the
interior of $(H+H')^{-1}(0)$
and $d(f_H + f_{H'})(V_H + V_{H'}) > 0$ on $\partial \left(N_C \times N'_{C'}\right)$
and outside $N_C \times N'_{C'}$. 
Smooth the corners of $N_C \times N'_{C'}$ slightly
to give a new manifold $A$ with boundary so that
\begin{enumerate}
\item $V_H + V_{H'}$ is transverse to $\partial A$ and pointing outwards.
\item $d(f_H + f_{H'})(V_H + V_{H'}) > 0$ on $\partial A$
and outside $A$.
\item $A \subset N_C \times N'_{C'}$.
Hence $A$  is contained in the
interior of $(H+H')^{-1}(0)$.
\end{enumerate}
Choose some positive function $g : \widehat{N} \times \widehat{N}' \rightarrow (0,\infty)$
which is small so that $g.(V_H + V_{H'})$ is an integrable vector field.
Flow $\partial A$ along $g.(V_H + V_{H'})$ for all time
so that we get a diffeomorphism from $(\widehat{N} \times \widehat{N}') \setminus A^o$
to $[1,\infty) \times \partial A$ such that $g.(V_H + V_{H'})$
maps to $\frac{\partial}{\partial r}$ where $r$ parameterizes $[1,\infty)$.
We can extend this to $(0,\infty) \times \partial A$ by flowing $\partial A$
backwards along $g.(V_H + V_{H'})$.
So $g.(V_H + V_{H'})$ still maps to $\frac{\partial}{\partial r}$
and $\partial A$ is identified with $\{1\} \times \partial A$.
Let $h : (0,\infty) \rightarrow \R$ be a function with
\begin{enumerate}
\item $h'(x) \geq 0$, $h(x) = 0$ for $x \leq \frac{1}{2}$
\item $h'(x)>0$,$h>1$ for $x \geq \frac{3}{4}$
\item $h(x)$ tends to infinity as $x$ tends to infinity.
\end{enumerate}
We define $f_{H+H'} : \widehat{N} \times \widehat{N}' \rightarrow \R$ as $h(r)$ when $r$ is well defined and $0$ elsewhere.
This is exhausting and $f_{H+H'}^{-1}(\infty,1]$
is contained in the
interior of $(H+H')^{-1}(0)$.
Also $df_{H+H'}(V_H + V_{H'}) > 0$ for $f_{H+H'} \geq 1$.
Also for $\epsilon$ small enough we  have that
$d(H+H')(V_H+V_{H'}) > 0$ in the region
$(H+H')^{-1}(0,\epsilon)$. This is because
if $H+H'$ is small and positive then $H$ and $H'$
are small and at least one of them is positive
and because $H$ and $H'$ satisfy:
\begin{enumerate}
\item $dH(V_H) \geq 0$ for $H$ small.
\item $dH(V_H) > 0$ if and only if $H>0$ for $H$ small.
\item the same properties are true for $H'$.
\end{enumerate}
Hence $H+H'$ satisfies the Liouville vector field property
and so is growth rate admissible.

The directed system $(SH_*(\lambda(H+H'),J \oplus J'))$
is equal to $(SH_*(\lambda H,J) \otimes SH_*(\lambda H',J'))$.
This is because we can choose approximating pairs for $(\lambda(H+H'),J \oplus J')$
which respect the product structure.
\qed

\section{Growth rate linear algebra}
\label{section:growthratelinearalgebra}

Recall that a filtered directed system is a family of vector spaces $(V_x)$
parameterized by $[1,\infty)$ forming a category where for $x_1 \leq x_2$
there is a unique homomorphism from $V_{x_1}$ to $V_{x_2}$ and no homomorphism when $x_1 > x_2$.
A morphism of filtered directed systems $\phi : (V_x) \rightarrow (V'_x)$
consists of some constant $C_\phi \geq 1$ and a sequence
of maps \[a_x : V_x \rightarrow V'_{C_\phi}x\]
so that we have the following commutative diagram:
\[
\xy
(0,0)*{}="A"; (40,0)*{}="B";
(0,-20)*{}="C"; (40,-20)*{}="D";
%
"A" *{V_{x_1}};
"B" *{V'_{C_\phi x_1}};
"C" *{V_{x_2}};
"D" *{V'_{C_\phi x_2}};
%
{\ar@{->} "A"+(12,0)*{};"B"-(15,0)*{}};
{\ar@{->} "C"+(12,0)*{};"D"-(15,0)*{}};
{\ar@{->} "A"+(0,-4)*{};"C"+(0,4)*{}};
{\ar@{->} "B"+(0,-4)*{};"D"+(0,4)*{}};
%
"A"+(20,3) *{a_{x_1}};
"C"+(20,3) *{a_{x_2}};
\endxy
\]
for all $x_1 \leq x_2$
where the vertical arrows come from the filtered directed systems.

Let $\psi_{x_1,x_2}$ be the natural map from $V_{x_1}$ to $V_{x_2}$
in this filtered directed system for $x_1 \leq x_2$.
For each constant $C \geq 0$, we have an automorphism $C_V$ from
$(V_x)$ to $(V_x)$ given by the map $\psi_{x,Cx}$.
We say that $(V_x)$ and $(V'_x)$ are {\it isomorphic}
if there is a morphism $\phi$ from $(V_x)$ to $(V'_x)$
and another morphism $\phi'$ from $(V'_x)$ to $(V_x)$
such that $\phi' \circ \phi = C_V$ and $\phi \circ \phi' = C'_{V'}$
where $C,C'\geq 0$ are constants and $C_V : (V_x) \rightarrow (V_x)$,
$C'_{V'} : (V'_x) \rightarrow (V'_x)$ are the automorphisms described above.
From \cite[Lemma 3.1]{McLean:affinegrowth} we have that
if $(V_x),(V'_x)$ are two isomorphic filtered directed systems,
then $\Gamma(V_x) = \Gamma(V'_x)$.

\begin{lemma} \label{lemma:filteredlongexactisomorphism}
Suppose $(V_x),(V'_x),(V''_x)$ are filtered directed systems
such that for all $x_1 \leq x_2$, we have the following commutative diagram
where the horizontal arrows are long exact sequences between $V_x,V'_x$ and $V''_x$
and the vertical arrows are the natural directed system maps:
\[
\xy
(0,0)*{}="A1"; (20,0)*{}="A2"; (40,0)*{}="A3";(60,0)*{}="A4"; (80,0)*{}="A5";
(0,-20)*{}="B1"; (20,-20)*{}="B2"; (40,-20)*{}="B3";(60,-20)*{}="B4"; (80,-20)*{}="B5";
%
"A2" *{V_{x_1}};
"A3" *{V'_{x_1}};
"A4" *{V''_{x_1}};
"B2" *{V_{x_2}};
"B3" *{V'_{x_2}};
"B4" *{V''_{x_2}};
%
{\ar@{.>} "A1"+(6,0)*{};"A2"-(6,0)*{}};
{\ar@{->} "A2"+(6,0)*{};"A3"-(6,0)*{}};
{\ar@{->} "A3"+(6,0)*{};"A4"-(6,0)*{}};
{\ar@{.>} "A4"+(6,0)*{};"A5"-(6,0)*{}};
{\ar@{.>} "B1"+(6,0)*{};"B2"-(6,0)*{}};
{\ar@{->} "B2"+(6,0)*{};"B3"-(6,0)*{}};
{\ar@{->} "B3"+(6,0)*{};"B4"-(6,0)*{}};
{\ar@{.>} "B4"+(6,0)*{};"B5"-(6,0)*{}};
{\ar@{->} "A2"+(0,-4)*{};"B2"+(0,4)*{}};
{\ar@{->} "A3"+(0,-4)*{};"B3"+(0,4)*{}};
{\ar@{->} "A4"+(0,-4)*{};"B4"+(0,4)*{}};
%
\endxy
\]
Suppose also that $(V''_x)$ is isomorphic
to the filtered directed system $(0)$ (i.e. all the vector
spaces are $0$).
Then $(V_x)$ is filtered isomorphic to $(V'_x)$.
\end{lemma}
\proof of Lemma \ref{lemma:filteredlongexactisomorphism}.
Because $(V''_x)$ is isomorphic to $(0)$, there is a constant $C>0$ such that
the directed system map $V''_x \rightarrow V''_{Cx}$ is $0$.
We look at the following commutative diagram:
\[
\xy
(0,0)*{}="A1"; (20,0)*{}="A2"; (40,0)*{}="A3";(60,0)*{}="A4"; (80,0)*{}="A5";
(0,-20)*{}="B1"; (20,-20)*{}="B2"; (40,-20)*{}="B3";(60,-20)*{}="B4"; (80,-20)*{}="B5";
(0,-40)*{}="C1"; (20,-40)*{}="C2"; (40,-40)*{}="C3";(60,-40)*{}="C4"; (80,-40)*{}="C5";
%
"A2" *{V_{x}};
"A3" *{V'_{x}};
"A4" *{V''_{x}};
"B2" *{V_{Cx}};
"B3" *{V'_{Cx}};
"B4" *{V''_{Cx}};
"C2" *{V_{C^2x}};
"C3" *{V'_{C^2x}};
"C4" *{V''_{C^2x}};
%
{\ar@{.>} "A1"+(6,0)*{};"A2"-(6,0)*{}};
{\ar@{->} "A2"+(6,0)*{};"A3"-(6,0)*{}};
{\ar@{->} "A3"+(6,0)*{};"A4"-(6,0)*{}};
{\ar@{.>} "A4"+(6,0)*{};"A5"-(6,0)*{}};
{\ar@{.>} "B1"+(6,0)*{};"B2"-(6,0)*{}};
{\ar@{->} "B2"+(6,0)*{};"B3"-(6,0)*{}};
{\ar@{->} "B3"+(6,0)*{};"B4"-(6,0)*{}};
{\ar@{.>} "B4"+(6,0)*{};"B5"-(6,0)*{}};
{\ar@{.>} "C1"+(6,0)*{};"C2"-(6,0)*{}};
{\ar@{->} "C2"+(6,0)*{};"C3"-(6,0)*{}};
{\ar@{->} "C3"+(6,0)*{};"C4"-(6,0)*{}};
{\ar@{.>} "C4"+(6,0)*{};"C5"-(6,0)*{}};
{\ar@{->} "A2"+(0,-4)*{};"B2"+(0,4)*{}};
{\ar@{->} "A3"+(0,-4)*{};"B3"+(0,4)*{}};
{\ar@{->} "A4"+(0,-4)*{};"B4"+(0,4)*{}};
{\ar@{->} "B2"+(0,-4)*{};"C2"+(0,4)*{}};
{\ar@{->} "B3"+(0,-4)*{};"C3"+(0,4)*{}};
{\ar@{->} "B4"+(0,-4)*{};"C4"+(0,4)*{}};
"A1"+(10,3) *{a^x_{31}};
"A2"+(10,3) *{a^x_{12}};
"A3"+(10,3) *{a^x_{23}};
"A4"+(10,3) *{a^x_{31}};
"B1"+(10,3) *{a^{Cx}_{31}};
"B2"+(10,3) *{a^{Cx}_{12}};
"B3"+(10,3) *{a^{Cx}_{23}};
"B4"+(10,3) *{a^{Cx}_{31}};
"C1"+(10,3) *{a^{C^2x}_{31}};
"C2"+(10,3) *{a^{C^2x}_{12}};
"C3"+(10,3) *{a^{C^2x}_{23}};
"C4"+(10,3) *{a^{C^2x}_{31}};
"A2"+(5,-10) *{\psi_{x,Cx}};
"A3"+(5,-10) *{\psi'_{x,Cx}};
"A4"+(2,-10) *{0};
"B2"+(7,-10) *{\psi_{Cx,C^2x}};
"B3"+(7,-10) *{\psi'_{Cx,C^2x}};
"B4"+(2,-10) *{0};
\endxy
\]
We have a map $\phi_x := a^x_{12}$ from $V_x$ to $V'_x$ which induces a morphism of filtered directed systems.
We now wish to create an inverse morphism $\phi'$
so that $\phi \circ \phi'$ and $\phi' \circ \phi$ are
directed system maps respectively (defined earlier in this section).
Here is how we construct $\phi'$:
We will construct it so that $\phi'_x$ is a map from $V'_x$
to $V_{C^2x}$.
Let $q \in V'_x$.
We have that $a^{Cx}_{23} \circ \psi'_{x,Cx}(q) = 0$
by the commutativity of the diagram.
This implies that
$\psi'_{x,Cx}(q) = a^{Cx}_{12}(w)$
for some $w \in V_{Cx}$ by the fact that we have a long exact sequence.
Let $w' \in V_{Cx}$ be another element such that 
$\psi'_{x,Cx}(q) = a^{Cx}_{12}(w')$.
Then by the long exact sequence property we have
that $w - w' = a^{Cx}_{31}(u)$ for $u \in V''_{Cx}$.
Commutativity implies that $\psi_{Cx,C^2x}(w-w') = 0$.
Hence $\psi_{Cx,C^2x}(w)$ is independent of the choice of $w$.
We define $\phi'_x(q) := \psi_{Cx,C^2x}(w)$.

We have that $\phi' \circ \phi = \psi_{Cx,C^2x} \circ \psi_{x,Cx}$
by commutativity of the diagram.
Similarly $\phi \circ \phi' = \psi'_{Cx,C^2x} \circ \psi'_{x,Cx}$.
Hence $\phi$ and $\phi'$ give us our isomorphism between
$(V_x)$ and $(V'_x)$.
\qed

\begin{lemma} \label{lemma:changeisomorphismconstant}
Let $(V_x),(V'_x)$ be isomorphic filtered directed systems.
For any constant $C$ sufficiently large, there exists maps
$\phi : V_x \rightarrow V'_{Cx}$ and $\phi' : V'_x \rightarrow V_{Cx}$
such that $\phi \circ \phi'$ and $\phi' \circ \phi$ are directed system
maps.
\end{lemma}
\proof of Lemma \ref{lemma:changeisomorphismconstant}.
Let $p : V_x \rightarrow V_{C_p x}$, $p' : V'_x \rightarrow V_{C_{p'} x}$
be the isomorphisms.
We choose any $C$ greater than both $C_p$ and $C_{p'}$.
We define $\phi := \psi_{C_p x, C x} \circ p$
and $\phi' := \psi_{C_{p'} x, C x} \circ p'$.
Because $p$ and $p'$ are morphisms of filtered directed systems,
they commute with the directed system maps and hence
because $p \circ p'$ and $p' \circ p$ are equal to
directed system maps then so are
$\phi \circ \phi'$ and $\phi' \circ \phi$.
\qed

We need a criterion that is invariant under isomorphism
that tells us when the growth rate of one filtered directed system
is greater than or equal to another one.
Let $(V_x),(V'_x)$ be filtered directed systems.
Let $\psi_{x_1,x_2},\psi'_{x_1,x_2}$ be the respective
directed system maps.
\begin{defn} \label{defn:growthrateinequalitydefn}
We say that $(V_x)$ is {\it bigger than}
$(V'_x)$ if there exists constants $A,B,C \geq 1$ with $C \geq B$
such that for all $x \geq 1$, $y \geq 1$,
\[\text{rank im}(\psi_{Bx,Cyx}) \geq  \text{rank im}(\psi'_{x,Ayx})\]
\end{defn}

\begin{lemma} \label{lemma:inequalityinvariance}
Suppose that $(V_x)$ is bigger than $V'_x$, and
suppose that $(\tilde{V}_x)$ (resp. $(\tilde{V}'_x)$)
is isomorphic to $(V_x)$ (resp. $V'_x$), then
$(\tilde{V}_x)$ is bigger than $(\tilde{V}'_x)$.
\end{lemma}
\proof of Lemma \ref{lemma:inequalityinvariance}.
Choose $A,B,C$ as in Definition \ref{defn:growthrateinequalitydefn}.
Also let $\phi_x : V_x \rightarrow \tilde{V}_{C_\phi x}$,
$\tilde{\phi}_x : (\tilde{V}_{x}) \rightarrow V_{C_{\tilde{\phi}}x}$
be the maps giving the isomorphism between
$(V_x)$ and $(\tilde{V}_x)$.
Similarly let
\[\phi'_x :  V'_x \rightarrow \tilde{V}'_{C_{\phi'} x}\] and 
\[\tilde{\phi}'_x : \tilde{V}'_{x} \rightarrow V'_{C_{\tilde{\phi}'}x}\] be the isomorphisms
between $(V'_x)$ and $(\tilde{V}'_x)$.
By Lemma  \ref{lemma:changeisomorphismconstant},
we can assume that \[\kappa = C_\phi = C_{\tilde{\phi}} = C_{\phi'} = C_{\tilde{\phi}'} > 1.\]
Choose any $y > 1$.
We have that
\[\psi_{\kappa Bx,C \kappa^3 y x} =
\tilde{\phi}_{\kappa^2 C y x} \circ \tilde{\psi}_{\kappa^2 B x, \kappa^2 C y x} \circ \phi_{\kappa Bx}.\]
Hence the rank of the image of
$\psi_{\kappa Bx,C \kappa^3 y x}$ is less than or equal
to the rank of the image of $ \tilde{\psi}_{\kappa^2 B x, \kappa^2 C y x}$.
We also have:
\[\tilde{\psi}'_{x,\kappa^4 A y x} =
\phi'_{\kappa^3 y A x} \circ \psi'_{\kappa x, \kappa^3 y A x} \circ \tilde{\phi}'_x.\]
Hence
\[\text{rank im} (\tilde{\psi}'_{x,\kappa^4 A y x}) \leq
\text{rank im}(\psi'_{\kappa x, \kappa^3 y A x}).\]
Because $(V_x)$ is bigger than $(V'_x)$
we get
\[\text{rank im}(\psi'_{\kappa x, \kappa^3 y A x}) \leq
\text{rank im}(\psi_{\kappa Bx,C \kappa^3 y x}).\]
This implies using the other inequalities that
\[\text{rank im}(\tilde{\psi}'_{x,\kappa^4 A y x}) \leq
\text{rank im}( \tilde{\psi}_{\kappa^2 B x, \kappa^2 C y x}).\]
This implies that $(\tilde{V}_x)$ is bigger than $(\tilde{V}'_x)$
where our constants $A,B,C$ (as described in definition
\ref{defn:growthrateinequalitydefn})
are replaced with constants
$\kappa^4 A, \kappa^2 B, \kappa^2 C$.
\qed

\begin{lemma} \label{lemma:inequalitygrowthratedifference}
If $(V_x)$ is bigger than $(V'_x)$ then
$\Gamma(V_x) \geq \Gamma(V'_x)$.
\end{lemma}
\proof of Lemma \ref{lemma:inequalitygrowthratedifference}.
Let $A,B,C \geq 1$ be the constants so that
for all $x,y \geq 1$,
\[\text{rank im}(\psi_{Bx,Cyx}) \geq  \text{rank im}(\psi'_{x,Ayx})\]
where $\psi,\psi'$ are the directed system maps for
$(V_x)$ and $(V'_x)$ respectively.
Let $f(x)$ (resp. $g(x)$) be the rank of
the image of the natural map
$a(x) : V_x \rightarrow \varinjlim_y V_y$
(resp. $b(x) : V'_x \rightarrow \varinjlim_y V'_y$).
We have that
$f(Bx) \geq g(x)$ by using the above inequality for large enough $y$.
The point is that because the rank of $V_x$ (resp. $V'_x$) is finite,
we have for large enough $y$ that $f(Bx)$ (resp. $g(x)$)
is the rank of the image of $\psi_{Bx,Cyx}$ (resp. $\psi'_{x,Ayx}$)
for large enough $y$.
So,
\[\varlimsup_x \frac{\log{f(x)}}{\log{x}} =
\varlimsup_x \frac{\log{f(Bx)}}{\log{Bx}} =\]
\[\varlimsup_x \frac{\log{f(Bx)}}{\log{x}} \geq
\varlimsup_x \frac{\log{g(x)}}{\log{x}}.\]
Hence $\Gamma(V_x) \geq \Gamma(V'_x)$.
\qed

\begin{lemma} \label{lemma:biggerthanrank}
Suppose that $(V_x),(V'_x)$ are filtered directed systems such that
$(V_x)$ is bigger than $(V'_x)$.
Then $(V'_x)$ is isomorphic to a filtered directed system
$(\tilde{V'_x})$ such that $|\tilde{V'_x}| \leq |V_{x}|$ for all $x \in [1,\infty)$.
\end{lemma}
\proof of Lemma \ref{lemma:biggerthanrank}.
Let $A,B,C \geq 1$ be the constants so that
for all $x,y \geq 1$,
\[\text{rank im}(\psi_{Bx,Cyx}) \geq  \text{rank im}(\psi'_{x,Ayx})\]
where $\psi,\psi'$ are the directed system maps for
$(V_x)$ and $(V'_x)$ respectively.
We have that $(\text{im}(\psi'_{x,Ax}))$ is a filtered directed system
where the directed system maps are the ones induced by $\psi'$.
This is because $\psi'_{x,Ax}$ is a morphism of filtered directed
systems and so its image is also a filtered directed system.
It is filtered isomorphic to $(V'_x)$ because we have maps:
$\phi : V_x \rightarrow \text{im}(\psi'_{x,Ax})$
given by $\phi = \psi'_{x,Ax}$ and
$\phi' :  \text{im}(\psi'_{x,Ax}) \rightarrow V_{Ax}$
given by the natural inclusion.
So $\phi \circ \phi'$ and $\phi' \circ \phi$
are the natural directed system maps.
Hence we have an isomorphism.
We define $\tilde{V'_x} :=  \text{im}(\psi'_{x/B,(A/B)x})$.
If $x < B$ then we define $\tilde{V'_x} := 0$.
This is a filtered directed system which (by using directed system maps)
is filtered isomorphic to $(\text{im}(\psi'_{x,Ax}))$.
We have
\[|\tilde{V'_x}| \leq \text{rank im}(\psi'_{x/B,(A/B)x}) \leq
\text{rank im}(\psi_{x,(C/B)x}) \leq |V_x|.\]
This proves the Lemma.
\qed

We will now state a technical lemma which will be used to give an upper bound for growth rates
in terms of open books (see Section \ref{section:lefschetzfibrationsbounds}).
Let $(C_\lambda)$ be a filtered directed system and
let $a_{\lambda_1,\lambda_2} : C_{\lambda_1} \rightarrow C_{\lambda_2}$
be the filtered directed system maps.
Let $\partial_\lambda$ be a differential on $C_\lambda$
so that the filtered directed system maps are chain maps.
Let $F^\lambda_0 \subset F^\lambda_1 \subset \cdots$ be a filtration
on $C_\lambda$ whose union is $C_\lambda$ so that:
$\partial_\lambda(F^\lambda_i) \subset F^\lambda_i$
and $a_{\lambda_1,\lambda_2}(F^{\lambda_1}_i) \subset F^{\lambda_2}_i$.
We define $F^\lambda_{-1} = 0$.
We have two filtered directed systems:
$(H_*(C_\lambda,\partial_\lambda))$
and $(\bigoplus_{i=0}^\infty H_*(F^\lambda_i / F^\lambda_{i-1},\partial_\lambda))$.
Let \[\bar{a}^i_{x,y} : H_*(F^x_i / F^x_{i-1},\partial_x)
\rightarrow H_*(F^y_i / F^y_{i-1},\partial_y)\]
be the induced directed system maps.
Let \[H(a_{x,y}) : H_*(C_x,\partial_x) \rightarrow H_*(C_y,\partial_y)\]
be the other induced directed system maps.
\begin{lemma} \label{lemma:biggerthanspectralsequence}
Suppose that there is a constant $M > 1$  such that
$\bar{a}^i_{x,y}$ are isomorphisms for all $y>x \geq Mi$.
Suppose also that there is a constant $N>1$ such that
$F^{\lambda}_i = F^\lambda_j$ for all $j \geq i \geq \lfloor N \lambda \rfloor$.
Then the filtered directed system $(\bigoplus_{i=0}^\infty H_*(F^\lambda_i / F^\lambda_{i-1},\partial_\lambda))$
is bigger than $(H_*(C_\lambda,\partial_\lambda))$.
\end{lemma}
\proof of Lemma \ref{lemma:biggerthanspectralsequence}.
We build a new filtered directed system $C'_\lambda$
with differential $\partial'_\lambda$ as follows:
We define $C'_\lambda$ to be equal to $F^\lambda_i$
in the region $Mi \leq \lambda < M(i+1)$.
A morphism from $C'_{\lambda_1}$ to $C'_{\lambda_2}$
where $Mi_1 \leq \lambda_1 < M(i_1+1)$
and $Mi_2 \leq \lambda_2 < M(i_2+1)$ ($i_1 \leq i_2$)
is the natural morphism from
$F^{\lambda_1}_{i_1}$ to $F^{\lambda_2}_{i_1}$
composed with the inclusion into $F^{\lambda_2}_{i_2}$.
We also define $\partial'_\lambda$ to be the induced differential.
This is  isomorphic to the filtered directed system $C_\lambda$,
where the isomorphism is built as follows:
We have a map $\phi : C_\lambda \rightarrow C'_{MN\lambda}$
given by the filtered directed system map $a_{\lambda,MN\lambda}$
because this map respects the filtration structure
and $C_\lambda \subset F^\lambda_{\lfloor N\lambda \rfloor}$.
Also we have a map
$\phi' : C'_\lambda \rightarrow C_{\lambda}$ given by inclusion
as $C'_\lambda$ is a subcomplex.
We have that $\phi \circ \phi'$ and $\phi' \circ \phi$
are directed system maps.
We also have an induced filtration structure
${F'}^\lambda_i := F^\lambda_i \cap C'_\lambda$.

This isomorphism also commutes with the differentials
and respects the filtration structure hence it induces filtered directed system isomorphisms
$(H_*(C_\lambda,\partial_\lambda)) \cong (H_*(C'_\lambda,\partial'_\lambda))$
and
 \[\left(\bigoplus_{i=0}^\infty H_*(F^\lambda_i / F^\lambda_{i-1},\partial_\lambda)\right)
\cong \left(\bigoplus_{i=0}^\infty H_*({F'}^\lambda_i / {F'}^\lambda_{i-1},\partial'_\lambda)\right).\]
All the filtered directed system maps
in $\left(\bigoplus_{i=0}^\infty H_*({F'}^\lambda_i / {F'}^\lambda_{i-1},\partial'_\lambda)\right)$
are injections because $ H_*({F'}^\lambda_i / {F'}^\lambda_{i-1},\partial'_\lambda))$
is non-trivial only when $Mi \leq \lambda$.
Let $\bar{a}'_{x,y}$ be the respective filtered directed system maps.
Let $H(a')_{x,y}$ be the directed system maps
for $(H_*(C'_\lambda,\partial'_\lambda))$.
We have that the rank of
$(H_*(C'_\lambda,\partial'_\lambda))$ is less than or equal
to the rank of $\left(\bigoplus_{i=0}^\infty H_*({F'}^\lambda_i / {F'}^\lambda_{i-1},\partial'_\lambda)\right)$
by a spectral sequence argument.
Also for all $y \geq 1$, the rank of the image of
$\bar{a}'_{\lambda,y\lambda}$
is equal to the rank of
 $(\bigoplus_{i=0}^\infty H_*({F'}^\lambda_i / {F'}^\lambda_{i-1},\partial'_\lambda))$
as the filtered directed system maps are injective.
The rank of the image of $H(a')_{\lambda,y\lambda}$
is less than or equal to the rank of
$H_*(C'_\lambda,\partial'_\lambda)$
hence we get that
\[\text{rank im} (\bar{a}'_{\lambda,y\lambda}) \geq
\text{rank im} (H(a')_{\lambda,y\lambda}).\]
This implies that
 $(\bigoplus_{i=0}^\infty H_*({F'}^\lambda_i / {F'}^\lambda_{i-1},\partial'_\lambda))$
is bigger than
$(H_*(C'_\lambda,\partial'_\lambda))$.
Hence by Lemma \ref{lemma:inequalityinvariance},
$(\bigoplus_{i=0}^\infty H_*(F^\lambda_i / F^\lambda_{i-1},\partial_\lambda))$
is bigger than $(H_*(C_\lambda,\partial_\lambda))$.
\qed

The following technical lemma will be used in Section
\ref{section:fillingsofalgebraiclefschetzfibrations}
so that we have a bound on the growth rate.
Let $(Q_\lambda)$ be a filtered directed system and
let $q_{\lambda_1,\lambda_2}$ be the respective filtered directed system maps.
Suppose as a vector space $Q_\lambda = A_\lambda \oplus B_\lambda$.
The filtered directed system maps $q_{\lambda_1,\lambda_2}$
can be viewed as a matrix
\[ \left( \begin{array}{cc}
a_{\lambda_1,\lambda_2} & b^a_{\lambda_1,\lambda_2} \\
a^b_{\lambda_1,\lambda_2} & b_{\lambda_1,\lambda_2}
\end{array} \right)\]
where $a_{\lambda_1,\lambda_2}$ is a map from $A$ to $A$,
$b_{\lambda_1,\lambda_2}$ is a map from $B$ to $B$,
$b^a_{\lambda_1,\lambda_2}$ is a map from $B$ to $A$
and
$a^b_{\lambda_1,\lambda_2}$ is a map from $A$ to $B$.
Suppose that we also have a differential $\partial_q$ on $Q_\lambda$
again of the form
\[ \left( \begin{array}{cc}
\partial_a & \partial_{b^a} \\
\partial_{a^b} & \partial_b
\end{array} \right)\]
which commutes with the filtered directed system maps $q_{\lambda_1,\lambda_2}$.
We will assume that $B_\lambda$ is a filtered directed system with filtered
directed system maps $b_{\lambda_1,\lambda_2}$.
We will also assume that the map $\partial_b$ is a differential
on $B_\lambda$ that commutes with the filtered directed system maps
$b_{\lambda_1,\lambda_2}$.

\begin{lemma} \label{lemma:chaincomplexhomologyupperbound}
Suppose that the rank of $A_\lambda$ is bounded above by some
function $P(\lambda)$ and that the filtered directed system
$(H_*(Q_\lambda,\partial_q))$ is isomorphic to some filtered
directed system $(V_\lambda)$ satisfying $|V_\lambda| \leq R(\lambda)$.
Then $(H_*(B_\lambda,\partial_b))$ is isomorphic as a filtered directed system
to $(W_\lambda)$ satisfying $|W_\lambda| \leq 6P(\lambda) + R(C'\lambda)$
for some constant $C'>1$.
\end{lemma}
\proof of Lemma \ref{lemma:chaincomplexhomologyupperbound}.
Because $(H_*(Q_\lambda,\partial_q))$ is isomorphic to $(V_\lambda)$,
there are constants $C,C'>1$ such that the map
$H_*(q_{\lambda,C\lambda})$ factors through $V_{C'\lambda'}$
for all $\lambda$.
We define $(B'_\lambda)$ to be the filtered directed system
$\text{im}(b_{\lambda,C\lambda})$.
We have a differential $\partial'_b$ on $B'_\lambda$ induced by $\partial_b$
because it commutes with the filtered directed system maps.
We have that $(H_*(B_\lambda,\partial_b))$ and
$(H_*(B'_\lambda,\partial'_b))$ are isomorphic as filtered directed systems
where the isomorphism is induced by the map $b_{\lambda,C\lambda}$
and the inclusion map of $B'_\lambda$ into $B_\lambda$.

We have that the rank of $H_*(B'_\lambda,\partial'_b)$
is equal by the first isomorphism theorem to
\[|B'_\lambda| - 2|\text{im} \partial'_b|.\]
Let $Q'_\lambda$ be equal to the image of $q_{\lambda,C\lambda}$.
Let $\partial'_q$ be the differential on $Q'_\lambda$.
We will view $\partial_a$ as a map from $A_\lambda \oplus B_\lambda$
to $A_\lambda \oplus B_\lambda$ by first projecting to $A_\lambda$ and then composing with $\partial_a$
and including the result into $A_\lambda \oplus B_\lambda$.
We will also view the maps $\partial_b$ , $\partial_{a^b}$ and $\partial_{b^a}$ in a similar way.
This means \[\partial_q = \partial_a + \partial_b + \partial_{a^b} + \partial_{b^a}.\]
The rank of the image of $\partial'_q$ is equal to the rank
of the image of $q_{\lambda,C\lambda} \circ \partial_q$
and the rank of the image of $\partial'_b$
is also equal to the rank of the image of
$q_{\lambda,C\lambda} \circ \partial_b$.
We have
\[|\text{im}(\partial'_b)| =
|\text{im}(q_{\lambda,C\lambda} \circ \partial_b)| \geq
|\text{im}(q_{\lambda,C\lambda} \circ \partial_q)| -
|\text{im}(q_{\lambda,C\lambda} \circ \partial_a)|\]
\[- |\text{im}(q_{\lambda,C\lambda} \circ \partial_{a^b})|
- |\text{im}(q_{\lambda,C\lambda} \circ \partial_{b^a})|
\geq |\text{im}(\partial'_q)| - 3|A_\lambda|.\]
Hence
\[|H_*(B',\partial'_b)| \leq |B'_\lambda| - 2(|\text{im}(\partial'_q)| - 3|A_\lambda|)
\leq |Q'_\lambda| - 2|\text{im}(\partial'_q)| + 6|A_\lambda|\]
\[ = |H_*(Q'_\lambda,\partial'_q)| + 6|A_\lambda|\]
We have that $|A_\lambda| \leq P(\lambda)$
and because $H_*(q_{\lambda,C\lambda})$ factors through $V_{C'\lambda}$
we have that $|H_*(Q'_\lambda,\partial'_q)| \leq R(C'\lambda)$.
Hence
\[|H_*(B',\partial'_b)| \leq R(C'\lambda) + 6P(\lambda).\]
Hence we have proven our lemma with $W_\lambda = H_*(B',\partial'_b)$.
\qed

\section{Lefschetz fibrations} \label{section:lefschetzfibrations}

\subsection{Partial Lefschetz fibrations} \label{subsection:partialefschetzfibrations}

Instead of dealing directly with open books, we will deal
with something which is basically equivalent to an open book called
a partial Lefschetz fibration.

We will now define partial Lefschetz fibrations.
A {\it partial Lefschetz fibration} $\pi : E \setminus K \rightarrow S$ is defined as follows:
The manifold $S$ is a compact surface with boundary, and
$E$ is a manifold with boundary and corners.
The manifold $E$ consists of two codimension $1$-boundary components
$\partial_h E$ and $\partial_v E$ meeting in a codimension $2$ component.
The set $K$ is a compact subset of the interior of $E$.
There is a $1$-form $\theta_E$ on $E$ making $E$ into a Liouville
domain after smoothing the corners.
The map $\pi$ must satisfy the following properties:
\begin{enumerate}
\item A neighbourhood of $\partial_h E$ is diffeomorphic
to $S \times (1-\epsilon,1] \times \partial F$ where
$F$ is some Liouville domain called the fiber of $\pi$.
Here $\theta_E = \theta_S + r_F \alpha_F$ where
$\theta_S$ is a Liouville form on $S$ and $r_F$ parameterizes the interval.
The $1$-form $\alpha_F$ is the contact form on $\partial F$.
The map $\pi$ is the projection map to $S$ here.
\item  $\theta_E$ restricted to the fibers of $\pi$
is non-degenerate away from the singularities of $\pi$.
\item
We have that $\pi|_{\partial_v E}$ is a fibration whose
fibers are exact symplectomorphic to $F$ and such the fibers of $\pi$ are either
disjoint or entirely contained in $\partial_v E$.
\item There are only finitely many singularities of $\pi$
and they are all disjoint from the boundary $\partial E$.
They are modelled on non-degenerate holomorphic
singularities.
\end{enumerate}
The Liouville domain $F$ is called the {\it fiber} of this
partial Lefschetz fibration.
For the purposes of this paper it does not matter
too much what the singularities of $\pi$ are.
In fact by enlarging the set $K$, we can assume
that $\pi$ has no singularities.
We call $\partial_h E$ the {\it horizontal boundary}
and $\partial_v E$ the {\it vertical boundary}.
Near the boundary $\partial_v E$, we have a connection
given by the $\omega_E$-orthogonal plane field
to the fibers.
Because the fibration is a product near $\partial_h E$,
the parallel transport maps associated to this connection
are well defined and are compactly supported if we
transport around a loop.
We call the symplectomorphism $\phi : F \rightarrow F$
given by parallel transporting around a loop on $\partial S$
the monodromy symplectomorphism around this boundary component.
If $S$ has a single boundary component then $\phi$ is called the monodromy
symplectomorphism of $\pi$.

We define $\alpha_\partial$ to be $\theta_E|_{\partial_v E}$.
This is a contact form on $\partial_v E$.
In the region $E_h = S \times (1-\epsilon,1] \times \partial F$,
we have that $\alpha_\partial = r_F \alpha_F + \alpha_S$
where $\alpha_S$ is $\theta_S|_{\partial_S}$ pulled back to $\partial_v E$
via $\pi$.
A partial Lefschetz fibration $\pi$ is said to be in {\it standard form}
if there is a neighbourhood $(1-\epsilon_S,1] \times \partial_v E$
of $\partial_v E$ where $\theta_E = (r_S-1) \alpha_S + \alpha_\partial$ where $r_S$ parameterizes the interval.
Also there is a neighbourhood $(1-\epsilon_S,1] \times \partial S$ of $\partial S$
where $\pi$ is the map $(\text{id},\pi|_{\partial_v E})$.
The good thing about partial Lefschetz fibrations in standard form is that we can
form their {\it completion}
\[\pi : \widehat{E} \rightarrow \widehat{S} \]
as follows:
we first glue on $S \times [1,\infty) \times \partial F$
to the horizontal boundary and extend $\pi$ as the projection map
to $S$. We also extend $\theta_E$ as $\theta_S + r_F \alpha_F$
over this region.
Let $\tilde{\pi} : \tilde{E} \setminus K \rightarrow S$
be the resulting map.
The region $\partial \tilde{E} = \tilde{\pi}^{-1}(\partial S)$
is a contact manifold with contact form $\alpha_\partial := \theta_{\tilde{E}}|_{\partial \tilde{E}}$.
This is in fact a union of mapping tori
(as described in \cite{McLean:spectralsequence}).
Let $\pi_\partial : \partial \tilde{E} \twoheadrightarrow \partial S$
be equal to $\tilde{\pi}|_{\partial \tilde{E}}$.
We call $\tilde{E}$ a {\it vertically completed partial Lefschetz fibration}.
Because $E$ is in standard form,
a neighbourhood of $\partial \tilde{E}$ is diffeomorphic to
$(1-\epsilon_S] \times \partial \tilde{E}$ with
$\theta_E = (r_S - 1)\alpha_S + \alpha_\partial$.
Here by abuse of  notation
we write $\alpha_S$ as the pullback of
$\alpha_S$ via $\pi_\partial$ to $\partial \tilde{E}$.
Hence we can glue on $[1,\infty) \times \partial \tilde{E}$
to $\partial \tilde{E}$ and extend $\theta_E$
by $(r_S - 1) \alpha_S + \alpha_\partial$.
We also extend $\pi$ to
\[ (\text{id}, \pi_\partial) : [1,\infty) \times \partial \tilde{E} \rightarrow [1,\infty) \times \partial S\]
in this region where $[1,\infty) \times \partial S$ is the cylindrical end of
$\widehat{S}$.
We write $\widehat{E}$ as the resulting manifold.
By abuse of notation we write $\pi$, $\theta_E$ for the associated
projection map and $1$-form on this manifold.
Here $\widehat{E}$ is called the {\it completion}.
If we smooth the boundary of $E$ slightly to create a manifold $E' \subset E$
then $(E',\theta_E)$ can be made into a Liouville domain and $\widehat{E}$
is in fact exact symplectomorphic to the completion $\widehat{E}'$.

The problem is that not every partial Lefschetz fibration is in standard form.
A {\it deformation of partial Lefschetz fibrations} is a smooth family of $1$-forms
$\theta_t$ on $E$ making $E$ into a Lefschetz fibration.
We require that the trivialization $S \times (1-\epsilon,1] \times \partial F$
and $\pi$ are fixed.
\begin{lemma} \label{lemma:partiallefschetzfibrationinstandardform}
Let $\pi : E \rightarrow S$ be a partial Lefschetz fibration.
Then this partial Lefschetz fibration is deformation equivalent to one which is in standard form.
\end{lemma}
\proof of Lemma \ref{lemma:partiallefschetzfibrationinstandardform}.
We have a natural connection on $E$ given by the planes that are $\omega_E$
orthogonal to the fibers (away from the singularities and the region $K$).
Let $X_{\theta_S}$ be the $\omega_S$-dual of $\theta_S$.
Let $\widetilde{X}_{\theta_S}$ be its lift.
We have $-\widetilde{X}_{\theta_S}$ is integrable near $\partial_v E$ because $\widetilde{X}_{\theta_S} = X_{\theta_S}$
in the region $E_h := S \times (1-\epsilon,1] \times \partial F$.
We first flow back $\partial_v E$ along $\widetilde{X}_{\theta_S}$
so that we get a region diffeomorphic to
\[ (-\epsilon',0] \times \partial_v E \]
with 
\[\theta_E = \alpha_\partial + \rho + dg \]
where $g$ is a function which vanishes inside
$E_h$, $\rho$ is a $1$-form which vanishes when restricted
to the fibers. Here $\epsilon'>0$ is a small constant.
Let $t$ parameterize the interval $(-\epsilon',0]$.
Because $X_{\theta_S}$ is the Liouville flow of $\theta_S$
we get that $r_S = e^t$.
Hence $\rho = (r_S-1) \alpha_S$ in the region $E_h$
as $\widetilde{X}_{\theta_S} = X_{\theta_S}$ inside $E_h$.
So we get that a neighbourhood of $\partial_v E$
is diffeomorphic to $(1-\epsilon'_S,1] \times \partial_v E$
with $\theta_E = \alpha_\partial + \rho + dg$
where $g = 0$ and $\rho = (r_S -1) \alpha_S$ inside $E_h$.
Here $\epsilon'_S = 1- e^{-\epsilon'}$.

We have a non-decreasing bump function $\nu : (1-\epsilon'_S,1] \rightarrow [0,1]$
which is $1$ near $1$ and $0$ near $1-\epsilon'_S$.
Note that for $\epsilon$ small we can define a new bump function
$\mu(x)=\nu(x + \epsilon)$ when $x \leq 1 - \epsilon$ and $\mu = 1$
otherwise.

\begin{figure}[H]
\centerline{
 \scalebox{1.0}{
   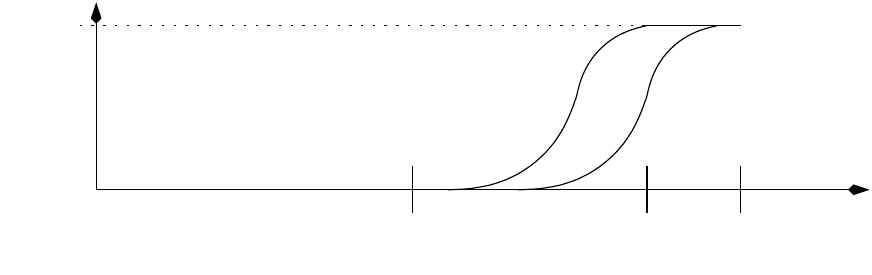
   }
   }
   \end{figure}

Let $C_t\geq 0$ be a smooth family of constants with $C_0 = 0$ such that
\[ \theta^t_E := (1-t\nu(r_S))\theta_E + t\nu(r_S) ( (r_S-1) \alpha_S + \alpha_\partial) + C_t \mu(r_S) r_S \alpha_S\]
makes $E$ into a Lefschetz fibration for all $t$
(basically by \cite[Theorem 2.15]{McLean:symhomlef}).
This is a deformation of Lefschetz fibrations.

For some small enough $\epsilon''>0$ we get that
\[ \theta^1_E = (r_S - 1) \alpha_S + \alpha_\partial + C_1 r_S \alpha_S\]
inside the region $\{r_S > 1-\epsilon''\}$.
We define $\alpha^1_\partial := \theta^1_E|_{\partial_v E}$,
\[\theta^1_S := \theta_S + C_1 \mu(r_S) r_S \alpha_S\]
and $\alpha^1_S := \theta^1_S|_{\partial S}$.
We have that $(S,\theta^1_S)$ is a Liouville domain.
For $\epsilon_S>0$ small enough we get that:
\[ \theta^1_E = (r_S - 1) \alpha^1_S + \alpha^1_\partial\]
in the region $\{r_S > 1-\epsilon_S\}$.
Hence $(E,\theta^1_E)$ is a partial Lefschetz fibration which is in standard form.
\qed

If we have some partial Lefschetz fibration $E$ then we define its {\it completion}
$\widehat{E}$ as follows:
We first deform $E$ so that it is in standard form and then we complete it as before.
From now on we will assume that all partial Lefschetz fibrations are in standard form
unless stated otherwise.
A {\it Lefschetz fibration} is defined as a partial Lefschetz fibration
$\pi : E \rightarrow S$ which is well defined everywhere
(i.e. the set $K$ is empty).

\subsection{An upper bound for growth rate} \label{section:lefschetzfibrationsbounds}

Let $(E,\pi)$ be a partial Lefschetz fibration whose fiber is $F$ and whose base
is some disk $\D$ in $\C$.
Let $\phi : F \rightarrow F$ be the monodromy map around $\partial \D$.
We have Floer homology groups $HF_*(\phi,k)$.
These groups have finite rank, so we can define the following
function:
\[b(x) := 1+\text{rank} \bigoplus_{k \leq x} HF_*(\phi,k).\]
This has an associated growth rate:
$\Gamma(\phi) := \varlimsup_x \frac{\log{b(x)}}{\log{x}}$.
The aim of this section is prove the following theorem
\begin{theorem} \label{theorem:lefschetzgrowthratebound}
$\Gamma(\widehat{E}) \leq \Gamma(\phi)$.
\end{theorem}
The reason why we want to do this is because it is a crucial ingredient
for proving Theorem \ref{theorem:subcriticalhandleattachaffinevariety}.
This theorem might be interesting  in its own right.
For instance if we combine it with Lemma \ref{lemma:partiallefschetzfibrationexistence}
we get an upper bound for $\Gamma(\widehat{E})$ in terms of
$HF_*(\psi,k)$ where $\psi$ is now the monodromy of any open book
supporting the contact structure on a smoothing of $\partial E$.

The growth rate of $\widehat{E}$ depends on a choice
of trivialization $\widetilde{\tau}$ of the canonical bundle on $\widehat{E}$
and a choice of some class $\widetilde{b} \in H_2(\widehat{E},\Z/2\Z)$.
The growth rate of $\phi$ also depends on a family of trivializations
of the canonical bundle on $\widehat{F}$ and a class $b \in H^2(F,\Z / 2\Z)$
which is fixed by the symplectomorphism.
The family of trivializations comes from the trivialization
$\widetilde{\tau}$ restricted to $\pi^{-1}(\partial \D)$ as follows:
Because we have a canonical trivialization of $T\D$
and a trivialization $\widetilde{\tau}$ of the canonical bundle of $E$,
we get a trivialization $\tau_\partial$ of the top exterior power of the vertical
bundle of $\pi$ restricted to $\pi^{-1}(\partial \D)$.
Because $\pi^{-1}(\partial \D) = F \times [0,1] / \sim$
where $\sim$ identifies $(0,f)$ with $(1,\phi(f))$,
we can lift this trivialization $\tau_\partial$
to $F \times [0,1]$. This can be viewed as a family
of trivializations $\tau_s$ of the canonical bundle of $F$
parameterized by $[0,1]$ such that $\phi^* \tau_1 = \tau_0$.
The choice of our class $b \in (F, \Z/ 2\Z)$ is just $\widetilde{b}$
restricted to the fiber $F$ of $\pi$.
From now all these Floer homology groups are defined
with respect to these choices.

We will prove Theorem \ref{theorem:lefschetzgrowthratebound} by proving a stronger Theorem:
We have a filtered directed system as follows:
we define \[V^\phi_x := H^{n-*}(E) \oplus \bigoplus_{i=1}^{\lfloor x
\rfloor} HF_*(\phi,k).\]
The morphism between $V^\phi_x$ and $V^\phi_y$
for $x \leq y$ is the natural inclusion map.
\begin{theorem} \label{theorem:lefschetzbiggerthan}
We have that $(V^\phi_x)$ is bigger than the filtered directed system
$(SH_*^+(\widehat{E},\theta_E,\lambda))$.
\end{theorem}
This Theorem combined with
Lemma \ref{lemma:inequalitygrowthratedifference}
implies Theorem \ref{theorem:lefschetzgrowthratebound}.
Here we also used the fact that the growth rate of a strictly
positive function does not change if we add any non-negative constant.
Hence all we need to do is prove Theorem \ref{theorem:lefschetzbiggerthan}.

We will now construct a pair $(H_\pi,J_\pi)$
for $\widehat{E}$ such that $(SH_*(\lambda H_\pi,J_\pi))$
is isomorphic to $(SH_*^+(\widehat{E},\theta_E,\lambda))$.
This pair can be defined for partial Lefschetz fibrations whose
base is any surface we like although we are usually interested in the disk.
Later on we will define families that only work when the base is a disk.

The base surface $S$ is a Liouville domain, and its completion $\widehat{S}$
has a cylindrical end $[1,\infty) \times \partial S$.
Let $\overline{r}_S$ be the cylindrical coordinate for this cylindrical end.
Let $\pi : \widehat{E} \twoheadrightarrow \widehat{S}$ be the completion
of the Lefschetz fibration above.
We will write $r_S$ for $\pi^* \overline{r}_S$.
Also the set $\{r_S \leq C\}$ really means $\pi^{-1}(S) \cup \{r_S \leq C\}$.
Similarly we have the coordinate $r_F$ which parameterizes the interval
in the region
\[\widehat{E}_h := \widehat{S} \times [1,\infty) \times \partial F.\]
In this region $\pi$ is the natural projection to $\widehat{S}$
and $\theta_E = \theta_S + r_F \alpha_F$.
In the region $\{r_S \geq 1\}$ we have that
$\theta_E = (r_S - 1)\alpha_S + \alpha_\partial$
where $\alpha_\partial = \theta_E|_{\pi^{-1}(\partial S)}$.

An almost complex structure $J$ is called {\it Lefschetz admissible}
if there exists a constant $C \geq 1$ such that in the region $\{r_F \geq C\} \subset \widehat{E}_h$,
$J = J_S \bigoplus J_F$ where $J_S$ (resp. $J_F$) is an almost complex structure
on $\widehat{S}$ (resp. $\widehat{F}$) compatible with the symplectic form which is cylindrical at infinity.
We also require that $\pi$ is $(J,J_S)$ holomorphic in the region $\{r_S \geq C\}$.
Let $h : [1,\infty) \rightarrow \R$ be a function
such that $h(x) = 0$ for $x$ near $1$, $h'(x),h''(x) \geq 0$
and $h'(x)=1$ near infinity.
We write $h(r_F)$ (resp. $h(r_S)$) as a function on $\widehat{E}$
which we extend by $0$ over the region where $r_F$ (resp. $r_S$)
is ill defined.
We define $H_\pi := h(r_F) + h(r_S)$
and $J_\pi$ to be any Lefschetz admissible almost complex structure.
After perturbing $E$ very slightly we may assume that
the period spectrum of $\partial F$ (where $F$ is the fiber)
is discrete.

\begin{lemma} \label{lemma:lefschetzpairadmissible}
The pair $(H_\pi,J_\pi)$ is growth rate admissible.
Also ${\mathcal A}_H := -\theta_E(X_{H_\pi}) - H_\pi \geq 0$.
\end{lemma}

Instead of proving this lemma, we will prove a more general
lemma which will be used later on in this paper.
Let $(H_F,J_F)$ be a growth rate admissible pair
on $\widehat{F}$ where $F$ is a subset of the interior of $H_F^{-1}(0)$.
We let $h : [1,\infty) \rightarrow [0,\infty)$ be a function equal to
$0$ near $1$ with $h'(r) = 1$ for $r \geq 2$.
We also assume that $h',h'' \geq 0$.
We write $h(r_S)$ for the function $\widehat{E}$ which is
equal to $h(r_S)$ when the coordinate $r_S$ is well defined
and zero elsewhere.
Let $J_E$ be an almost complex structure on $\widehat{E}$
making $\pi$ $(J_E,j_S)$ holomorphic where $j_S$ is a complex structure compatible with the symplectic form
on $\widehat{S}$ which is cylindrical on its cylindrical end.
We also assume that $J_E = j_S \oplus J_F$
in the region 
\[\widehat{E}_h := \widehat{S} \times [1,\infty) \times \partial F \subset \widehat{E}.\]
Let $\pi_2$ be the natural projection from $\widehat{E}_h$ to $[1,\infty) \times \partial F \subset \widehat{F}$.
By abuse of notation we write $\pi_2^* H_F$ as the Hamiltonian on $\widehat{E}$
defined by $\pi_2^* H_F$ inside $\widehat{E}_h$ and $0$ elsewhere.

\begin{lemma} \label{lemma:welldefinedgrowthrateadmissiblehamiltonian}
The pair $(H_E := \pi_2^* H_F + h(r_S),J_E)$ is growth rate admissible.
If $-\theta_F(X_{\theta_{H_F}}) - H_F \geq 0$ then
${\mathcal A}_H := -\theta_E(X_{\theta_{H_E}}) - H_E \geq 0$
(Here $\theta_{H_F}$ and $\theta_{H_E}$ are the $1$-forms
enabling $H_F$ and $H_E$ to satisfy the Liouville vector field property).
\end{lemma}

This lemma proves Lemma \ref{lemma:lefschetzpairadmissible}
when $H_F = h(r_F)$ and $J_F$ is cylindrical at infinity.
This is because such a pair is growth rate admissible
(see the second example in \cite[Section 4]{McLean:affinegrowth}).
Also because $h''(r_F) \geq 0$ we get that
$-\theta_F(X_{H_F}) - H_F \geq 0$ which implies that
${\mathcal A}_H \geq 0$.

\proof of Lemma \ref{lemma:welldefinedgrowthrateadmissiblehamiltonian}.

{\it $(H_E,J_E)$ is $SH_*$ admissible}: 
Let ${\mathcal P}_S$ the period spectrum of the contact boundary $\partial S$.
The Hamiltonian vector field $X_{h(r_S)}$ is equal to the
horizontal lift of the Hamiltonian vector field associated to $h(r_S)$
on the base $S$.
This means that for $\lambda$ not in ${\mathcal P}_S$ we have that
all the $1$-periodic orbits of $\lambda h(r_S)$ are contained
in the region where $r_S < 2$.
Outside $\widehat{E}_h$ we have that $H_E = h(r_S)$
which means that for $\lambda$ not in ${\mathcal P}_S$,
the orbits of $\lambda H_E$ outside $\widehat{E}_h$
are contained inside a fixed compact subset.

In the region $\widehat{E}_h$ we have that $H_E$
splits up as a product. Because $(H_F,J_F)$
is growth rate admissible there is a discrete
subset $A_{H_F}$ of $(0,\infty)$
and a sequence of relatively open compact subsets
$U^{H_F}_\lambda$ of $\widehat{F}$ satisfying
\begin{enumerate}
\item $U^{H_F}_{\lambda_1} \subset U^{H_F}_{\lambda_2}$ for
$\lambda_1 \leq \lambda_2$
\item
All the $1$-periodic orbits of $\lambda H_F$ are contained
in some closed subset of $U^{H_F}_\lambda$.
\end{enumerate}
This means that for $\lambda$ not in $A_{H_F}$
or ${\mathcal A}_S$ we have that all the $1$-periodic
orbits of $H_E$ that are contained in $\widehat{E}_h$
are also contained in $\{r_S < 2\} \times U^{H_F}_\lambda \subset
\widehat{E}_h$.
So we define $U^{H_E}_\lambda$ to be equal to
$\{\overline{r}_S < 2\} \times U^{H_F}_\lambda$ inside $\widehat{E}_h$
and equal to $\{r_S < 2\}$ outside $\widehat{E}_h$.

Suppose we have a map $u : S \rightarrow \widehat{E}$
satisfying the perturbed Cauchy-Riemann equations
whose boundary is contained in $U^{H_E}_\lambda$.
Let $\bar{S} \subset S$ be a subsurface with boundary such that
$u(\bar{S})$ is contained in $\widehat{E}_h$.
We can assume that $\pi_2 \partial \bar{S} \subset U^{H_F}_\lambda$.
We have that $\pi_2 \circ u|_{\bar{S}}$ satisfies the perturbed Cauchy-Riemann
equations with respect to $(H_F,J_F)$.
The maximum principle for $(H_F,J_F)$ ensures that $u(\bar{S}) \subset U^{H_F}_\lambda$.
Now let $\tilde{S} \subset S$ be a subsurface of $S$ mapped
by $u$ to the region $\{r_S > 1\}$ we can assume that the boundary
is mapped to a subset of $\{r_S \leq 2\}$.
Lemma \ref{lemma:lefschetzmaximumprinciple2} ensures that
$\tilde{S}$ is mapped to the region $\{r_S \leq 2\}$.
Hence $u(S)$ is contained in $U^{H_E}_\lambda$.
Hence $(H_E,J_E)$ satisfies the maximum principle.

{\it $H_E$ satisfies the bounded below property} because
$h(r) \geq 0$ and is greater than $1$ for $r$ sufficiently large
and because $H_F$ satisfies the bounded below property.

{\it $H_E$ satisfies the Liouville vector field property}: 
We add an exact $1$-form $d\psi$ to $\theta_F$ so that it is equal to
$\theta_{H_F}$ where $\theta_{H_F}$ is the $1$-form enabling the pair
$(H_F,J_F)$ to satisfy the Liouville vector field property.
We assume that $\psi = 0$ inside $F$ which is a closed subset of the interior of $H_F^{-1}(0)$.
Hence we can define $g_F$ to be equal to $\psi \circ \pi_2$ in $\widehat{E}_h$
and $0$ elsewhere.
We let $\theta_{H_E} := \theta_E + dg_F$.
The $\omega_E$-dual of $\theta_{H_E}$ in $\widehat{E}_h \cap \{r_S \geq 1\}$ is
$X_{\theta_{H_F}} + r_S \frac{\partial}{\partial r_S}$.
We have
\[dH_E(X_{\theta_E + dg_F}) - H_E = dH_F(X_{\theta_{H_F}+ dg_F}) - H_F + r_Sh'(r_S) - h(r_S)\]
inside $\widehat{E}_h \cap \{r_S \geq 1\}$.
In the region $\widehat{E}_h \cap \{r_S \leq 1\}$, $h(r_F)=0$
so:
\[dH_E(X_{\theta_E + dg_F}) - H_E = dH_F(X_{\theta_{H_F}+dg_F}) - H_F.\]
Using the above action calculations and the fact that
$H_F$, $h$, $h'$ and $h''$  are non-negative, we get the inequality
$dH_E(X_{\theta_E + dg_F}) - H_E > 0$ inside 
$\widehat{E}_h \cap \{H_E^{-1}(0,\epsilon_{H_F})\}$.
Because $\alpha_\partial \left( \frac{\partial}{\partial r_S}\right)=0$
and $\alpha_\partial \left( \widetilde{ \frac{\partial}{\partial \vartheta}}\right)>0$
where $\widetilde{ \frac{\partial}{\partial \vartheta}}$ is the horizontal lift of
$\frac{\partial}{\partial \vartheta}$ we get that the $\omega_E$-dual
of $\alpha_\partial$ is $X + \nu { \frac{\partial}{\partial r_S}}$
where $\nu>0$ is a function on $\pi^{-1}(\partial S)$
which is constant near infinity
and $X$ is tangent to the fibers of $\pi$.
Hence the $\omega_E$-dual of $\theta_E$ is 
\[(r_S-1) \frac{\partial}{\partial r_S} + X + \nu{ \frac{\partial}{\partial r_S}}\]
in the region $\{r_S \geq 1\}$.
Because $dh(r_S)$ is trivial when restricted to the fibers,
we have that inside $\{r_S \geq 1\} \setminus \widehat{E}_h$,
\[dH_E(X_{\theta_E + dg_F}) - H_E = (r_S-1)h'(r_S) + \nu h'(r_S) - h(r_S)\]
\[ = \int_1^{r_S} (t-1)h''(t) dt + \nu h'(r_S) \]
which is greater than $0$ when $h(r_S) > 0$.
Hence
$dH_E(X_{\theta_E + dg_F}) -H_E > 0$
in the region $\{H_E^{-1}(0,\epsilon_{H_F})\}$.

Let $f_F : \widehat{F} \rightarrow \R$ be an exhausting function
such that $df_F(X_{\theta_{H_F}}) > 0$ outside some closed subset of $H_F^{-1}(0)$.
We can ensure that $\pi_2^* f_F$ is zero near $\partial F$ and we can extend it by
zero outside $\widehat{E}_h$.
Let $h_1 = h(1+ \frac{r_S - 1}{2})$.
Then $(d(\pi_2^* f_F) + dh_1(r_S))(X_{\theta_E + dg_F}) > 0$ outside a closed subset of
$H_E^{-1}(0)$.
All of this means that $(H_E,J_E)$ satisfies the Liouville vector field property.

{\it $H_E$ satisfies the action bound property}:
The function $\nu$ above is bounded because $\alpha_\partial = \alpha_S + r_F \alpha_F$ in the region
$r_F \geq 1$ and $\nu$ is invariant under translations in the $r_S$ coordinate.
Hence $(r_S-1) h'(r_S) + \nu h'(r_S) - h(r_S)$ is bounded as $h$ is linear near infinity.
Because $H_F$ satisfies the action bound property,
we can add an exact $1$-form $d(q_F \circ \pi_2)$ to $\theta_E$ which we define to be zero outside $\widehat{E}_h$
so that $d(H_F \circ \pi_2)(X_{\theta_E + dq_E}) - H_F \circ \pi_2$ is bounded.
So
\[{\mathcal A}_H = (r_S-1) h'(r_S) + \nu(r_S)h'(r_S) - h(r_S) + d(H_F \circ \pi_2)(X_{\theta_E + dq_E}) - H_F \circ \pi_2\]
is bounded inside $\widehat{E}_h$
and
\[{\mathcal A}_H = (r_S-1)h'(r_S) + \nu(r_S)h'(r_S) - h(r_S)\]
away from $\widehat{E}_h$.
Hence $dH_E(X_{\theta_E + dq_E})$ is bounded.
This means that $(H_E,J_E)$ satisfies the action bound property.
Hence $(H_E,J_E)$ is growth rate admissible.

{\it ${\mathcal A}_H \geq 0$} :
Suppose that $-\theta_F(X_{H_F}) - H_F \geq 0$.
Outside the region $\widehat{E}_h$ we have that
\[{\mathcal A}_H = (r_S-1)h'(r_S) + \nu(r_S)h'(r_S) - h(r_S) =
 \nu h'(r_S) + \int_1^{r_S} (t-1) h''(t)dt \geq 0.\]
We have that ${\mathcal A}_H$ is equal to:
\[(r_S-1)h'(r_S) + \nu h'(r_S) - h(r_S) - \theta_F(X_{H_F}) - H_F \geq 0\]
inside $\widehat{E}_h$.
Hence ${\mathcal A}_H \geq 0$ everywhere.
\qed

Fix some $\epsilon_F > 0$ smaller than the length
of the smallest Reeb orbit of $\partial F$.
Let $(H_\lambda,J_\lambda)$ be smooth family of Hamiltonians
such that $H_\lambda = \epsilon_F h(r_F) + \lambda \kappa_1 h(r_S) + \kappa_2$
where $\kappa_1>0,\kappa_2$ are constants
and $J_\lambda$ is Lefschetz admissible.
We say that $(H_\lambda,J_\lambda)$
is a smooth family of {\it half Lefschetz admissible Hamiltonians
of slope} $\lambda$.
If we just have a family of Hamiltonians $H_\lambda^\pi$ then we say that they are Lefschetz admissible too.
The following Lemma is a technical Lemma. It will be used in this section
and in section \ref{section:fillingsofalgebraiclefschetzfibrations}.

\begin{lemma} \label{lemma:monodromyfiltereddirectedsystems}
Suppose that the base surface $S$ is connected and not contractible and let $A_S$
be a boundary component. Let $\alpha \subset H_1(E)$ be represented
by loops which project to loops which wrap around $A_S$ a non-zero number of times and let
$\phi$ be the associated monodromy map around $A_S$.
We also suppose that $0 \notin \alpha$.
The filtered directed system
$(SH_*^{\#,\alpha}(H_\lambda,J_\lambda))$
is isomorphic to $\left(\bigoplus_{i=1}^{\lfloor \lambda
\rfloor}(HF_*(\phi,i))\right)$
where the directed system maps are the natural inclusion maps.
\end{lemma}
\proof of Lemma \ref{lemma:monodromyfiltereddirectedsystems}.
We have that the component of the Lefschetz cylindrical end
corresponding to $A_S$ is equal to
$[1,\infty) \times M_\phi$ where $M_\phi$ is the mapping torus of $\phi$.
This mapping torus has a natural contact form $\alpha_\phi$
and $\theta_E = (r_S-1) \alpha_S + \alpha_\phi$ inside this cylindrical end
where $r_S$ parameterizes the interval
and $\alpha_S$ is the contact form on $\partial S$ which we pull back via $\pi$.

In order to calculate $SH_*^{\#,\alpha}(H_\lambda,J_\lambda)$
we need to perturb $H_\lambda$ slightly
(to create an approximating pair).
We have for some small $\epsilon_h>0$ that
$h(x) = 0$ for $x \in [1,1+\epsilon_h)$.
We perturb $(H_\lambda,J_\lambda)$
to an approximating pair $(H'_\lambda,J'_\lambda)$ so that
near $\{r_S = 1\}$ it is equal to
$l(r_S) + h(r_F)$ where $l'(1) < 0$.
We also assume that near $\{r_S = 1+\epsilon_h\}$,
$H'_\lambda$ is equal to $g(r_S) + h(r_F)$
where $g' > 0$ is very small.
By Lemma
\ref{lemma:lefschetzmaximumprinciple2}, all
$1$-periodic orbits representing some class in $\alpha$
and all the Floer trajectories connecting them stay inside
$[1,\infty) \times M_\phi$.
This is where our cohomological condition
$0 \notin \alpha$ is used so that we can construct our closed $1$-form $\beta$
as stated in Lemma \ref{lemma:lefschetzmaximumprinciple2}.
The point is that $\alpha$ contains a non-torsion class so represents a non-trivial
class in $H_1(E,\R)$, hence by the universal coefficient theorem we get
a non-trivial closed $1$-form $\beta$ extending $\pi^* \alpha_\phi$.
Also because these periodic orbits representing classes in $\alpha$ are actually contained
in the region $\{r_S > 1+\epsilon_h\}$ we have by
the maximum principle \cite[Lemma 5.2]{McLean:symhomlef}
that Floer trajectories connecting them are in fact contained in
$\{r_S > 1+\epsilon_h\}$.

We define $W_\phi := \{r_S > 1\}$
with $\theta_\phi := \theta_E|_{W_\phi}$.
We define $\pi_\phi := \pi$.
Let $\vartheta$ be the angle coordinate on $A_S$, then there is some constant
$\kappa_\phi>0$ such that $\alpha_S = \kappa_\phi d\vartheta$.
We let $r_\phi := \kappa_\phi(r_S - 1)$.
So $W_\phi = (0,\infty) \times M_\phi$
with $\theta_\phi = r_\phi d\vartheta + \alpha_\partial$.
We define an almost complex structure $J_\phi$
on $W_\phi$ making $\pi_\phi$ $(J_\phi,j)$-holomorphic with respect
to a natural complex structure on $A_S \times (0,\infty) \cong {\mathcal H}
/ \Z$ (${\mathcal H}$ is the upper half plane)
and such that it is equal to $J_\lambda$
inside the region $r_S \geq 1+ \frac{\epsilon_h}{2}$.
Here we might have to adjust $J_\lambda$ so that this $(J_\phi,j)$
holomorphic condition is satisfied.

We define
$\tilde{H}^\phi_\lambda$
to be equal to $H'_\lambda$ in the region
$\{r_S \geq 1+ \epsilon_h\}$
and outside this region we define it to be
equal to some function $\pi_\phi^* \tilde{h}(r_\phi) + h_F(r_F)$
where $\tilde{h} : (0,\infty) \rightarrow \R$ is some
function with small positive derivative so that
it is equal to $g(r_S)$ near $\{r_S = 1+ \epsilon_h\}$.
All the orbits of $\tilde{H}^\phi_\lambda$ and $H'_\lambda$
are identical inside $W_\phi$.
Also by \cite[Lemma 5.2]{McLean:symhomlef} we have that
their Floer trajectories are identical as well
as they must be contained in $\{r \geq 1+ \epsilon_h\}$.
Hence
$(SH_*^{\#,\alpha}(H_\lambda,J_\lambda))$ and
$(SH_*^{\#}(\tilde{H}^\phi_\lambda,J_\phi))$
are isomorphic as filtered directed systems.
Let $\alpha_k \subset H^1(W_\phi)$ be the subset of $H^1$
classes represented by loops which project to loops wrapping
$k$-times around the base $(0,\infty) \times A_S$.

By \cite[Lemma 2.9]{McLean:spectralsequence}
we have an isomorphism
$SH_*^{\#,\alpha_k}(\tilde{H}^\phi_\lambda,J_\phi) \cong HF_*(\phi,k)$
for $\lambda \geq 2 \kappa_\phi \pi k$.
Also the natural filtered directed system
maps
\[SH_*^{\#,\alpha_k}(\tilde{H}^\phi_{\lambda_1},J_\phi)
\rightarrow SH_*^{\#,\alpha_k}(\tilde{H}^\phi_{\lambda_2},J_\phi)\]
commute with this isomorphism for $2 \kappa_\phi \pi k \leq \lambda_1 \leq \lambda_2$.
Let $\beta_k := \cup_{i = 1}^k \alpha_i$.
We have the following isomorphism (which commutes with the filtered directed system maps):
\[SH_*^{\#,\beta_{\lfloor \lambda \rfloor}}(\tilde{H}^\phi_{2 \kappa_\phi \pi \lambda},J_\phi) \cong
\bigoplus_{k=1}^{\lfloor \lambda \rfloor} HF_*(\phi,k)\]
where
$\lfloor x \rfloor$ is the largest integer $\leq x$.
This means there is a natural morphism $\Phi$ of filtered directed systems
from $\left(\bigoplus_{k=1}^{\lfloor \lambda \rfloor} HF_*(\phi,k)\right)$ to
$(SH_*^{\#}(\tilde{H}^\phi_\lambda,J_\phi))$ induced by the above isomorphism
to the subgroup $SH_*^{\#,\beta_{\lfloor \lambda \rfloor}}(\tilde{H}^\phi_{2 \kappa \pi \lambda},J_\phi)$.

We also have a constant $K > 0$ such that
all orbits of $SH_*^{\#}(\tilde{H}^\phi_\lambda,J_\phi)$
wrap around the base at most $\lfloor K \lambda \rfloor$ times.
This means that we have a natural map
from:
$SH_*^{\#}(\tilde{H}^\phi_\lambda,J_\phi)$
into
$SH_*^{\#,\beta_{\lfloor K \lambda \rfloor}}(\tilde{H}^\phi_{2 \kappa_\phi \pi K \lambda},J_\phi)$
which is isomorphic to
$\bigoplus_{k=1}^{\lfloor K\lambda \rfloor} HF_*(\phi,k)$.

This induces the following morphism $\Phi'$ of filtered directed systems:
\[(SH_*^{\#}(\tilde{H}^\phi_\lambda,J_\phi)) \rightarrow
\left(\bigoplus_{k=1}^{\lfloor \lambda \rfloor} HF_*(\phi,k)\right).\]
Because these morphisms are compositions of continuation maps
and inclusion maps induced by classes in $H_1(W_\phi)$ we have that $\Phi \circ \Phi'$
and $\Phi' \circ \Phi$ are filtered directed system maps.
Hence we have our isomorphism from $(SH_*^{\#,\alpha}(H_\lambda,J_\lambda))$
to
$\left(\bigoplus_{k=1}^{\lfloor \lambda \rfloor} HF_*(\phi,k)\right)$.
\qed

If we have Hamiltonians $H_\lambda^\pi = \lambda (\kappa'_1 h(r_F) + \kappa'_2 h(r_S)) + \kappa'_3$
where $\kappa'_1,\kappa'_2>0$ and $\kappa'_3$ are constants then a smooth family of pairs
$(K_\lambda,Y_\lambda)$ is {\it Lefschetz admissible} if $K_\lambda = H_\lambda^\pi$ near infinity
and $Y_\lambda$ is a Lefschetz admissible almost complex structure.
\begin{lemma} \label{lemma:halflefschetzlefschetzequivalence}
Suppose the base $S$ is the disk $\D$
and $(H_\lambda,J_\lambda)$ are half Lefschetz admissible Hamiltonians of slope $\lambda$
then we have that the filtered directed system
$(SH_*(H_\lambda,J_\lambda))$ is isomorphic to
$(SH_*(\widehat{E},\theta_E,\lambda))$.
\end{lemma}

\proof of Lemma \ref{lemma:halflefschetzlefschetzequivalence}.
First of all if we have another choice $(H'_\lambda,J'_\lambda)$
with the same properties as stated above for $(H_\lambda,J_\lambda)$
then $(SH_*(H_\lambda,J_\lambda))$ is isomorphic to
$(SH_*(H'_\lambda,J'_\lambda))$.
In fact $H'_\lambda$ can be of the form $\lambda g(r_S) + \epsilon_F h(r_F)$
near infinity for any function $g$ where $g'$ is constant for $r_S$ large.
Suppose also we have any smooth family of pairs
$(K_\lambda,Y_\lambda)$
that are Lefschetz admissible.
Then by Lemmas \ref{lemma:lefschetzpairadmissible} and
\cite[Lemma 4.7]{McLean:affinegrowth}
we have that the filtered directed system
$(SH_*(K_\lambda,Y_\lambda))$ is isomorphic to
$(SH_*(\widehat{E},\theta_E,\lambda))$.
So in order to prove our Lemma, we need to construct
two families $(H'_\lambda,J'_\lambda)$, $(K_\lambda,Y_\lambda)$
as described above so that
$(SH_*(H'_\lambda,J'_\lambda)) = (SH_*(K_\lambda,Y_\lambda))$.
We will construct a long exact sequence
and then use Lemma \ref{lemma:filteredlongexactisomorphism}
to give us our isomorphism.
The proof is similar to the proof of \cite[Theorem 2.24]{McLean:symhomlef}.

We will first deform the partial Lefschetz fibration
$\pi : \widehat{E} \twoheadrightarrow \C = \widehat{\D}$ inside a compact set.
This does not affect our result.
We have a small neighbourhood of $\pi^{-1}(\partial \D)$
diffeomorphic to $(1-\epsilon_S,1+\epsilon_S) \times \pi^{-1}(\partial\D)$
where $\{1\} \times \pi^{-1}(\partial \D)$ is identified with
$\pi^{-1}(\partial \D)$.
The map $\pi$ here is the map $(\text{id},\pi|_{\pi^{-1}(\partial \D)})$
to $(1-\epsilon_S,1+\epsilon_S) \times \partial \D$.
Here $r_S$ parameterizes the interval $(1-\epsilon_S,1+\epsilon_S)$.
Let $p \in \{r_S = 1-\frac{\epsilon_S}{2}\} \subset S$.
The fibration $\pi$ is well defined in this region,
so let $U$ be a small neighbourhood of $p$ in $\D \subset \widehat{\D}$
diffeomorphic to a disk.
Let $(r,\vartheta)$ be polar coordinates for this disk $U$
such that $\{r \leq l\} = U$ and $dr^2 \wedge d\theta = d\theta_S$.
Basically by \cite[Lemma 3.1]{McLean:spectralsequence}
we can deform $\pi$ to a new partial Lefschetz fibration
so that $\pi$ is a trivial Lefschetz fibration
$U \times \widehat{F} \twoheadrightarrow U$
around $U$ with $\theta_E = \theta_S + \theta_F$.
Here $\theta_S$ is the Liouville form on the base $\C$.
We can add an exact $1$-form to $\theta_S$ so that $\theta_S = rd\theta$ in $U$, which means that
$\theta_E = rd\vartheta + \theta_F$ in $U \times \widehat{F}$.
This trivialization also extends to the trivialization
$(1-\epsilon_F, 1] \times \partial F \times S$ near $\partial_h E$.

We will now construct a family of Lefschetz admissible Hamiltonians
parameterized by three variables $m,\epsilon$ and $\lambda$.
Because the base $S$ is $\C$ and $U$ is disjoint from the cylindrical end of $\widehat{\D}$,
we can extend $r$ and $\vartheta$
so that $dr^2 \wedge d\vartheta = d\theta_S$ everywhere,
$d\vartheta = \kappa \alpha_S$ inside the cylindrical end of $S$ where $\kappa>0$ is a constant.
We can also ensure that $r = f(r_S)$ for some function $f$ inside the cylindrical end of $S$.
Because the integral of $r^2 d\vartheta$ around $\{r_S = C\}$
is equal to the integral of $r_S \alpha_S$ around this same curve for every $C \geq 1$
we get that $f(r_S) = \sqrt{\frac{r_S}{\kappa}}$.
We write (by abuse of notation) $r$ for the pullback of $r$
via $\pi$.
%
%
We define Hamiltonians $B_\lambda$ smoothly parameterized by
$\lambda$ such that in the region $\{r_S \geq 2\} \cup \{r_F \geq 1\}$ they are equal
to $\lambda r^2$.
In the region $\{r_S \leq 2\} \cap \{r_F \leq 1\}$, the Hamiltonian
$B_\lambda$ can be anything we like.

Let $h_{m,\lambda,\epsilon} : [1,\infty) \rightarrow \R$ be a function such that
$h_{m,\lambda,\epsilon}(t)$
is equal to zero near $t = 1$ and equal to $\epsilon (t-1)$
for $2 \leq t \leq m$.
For $t \geq m+1$ we require $h_{m,\lambda,\epsilon}(t) := \lambda(t - m)$.
We also require $h_{m,\lambda,\epsilon}',h_{m,\lambda,\epsilon}'' \geq 0$
(this works as long as $\epsilon$ is small
enough with respect to $m$).
We also require that $h_{m,\lambda,\epsilon}(t) < 1$ in the region $2 \leq t \leq m$
and that in the region where $h_{m,\lambda,\epsilon}(t) \geq 1$, $h_{m,\lambda,\epsilon}(t) = \lambda(t - m)$.
Here is a graph of this function:

\begin{figure}[H]
\centerline{
 \scalebox{1.0}{
  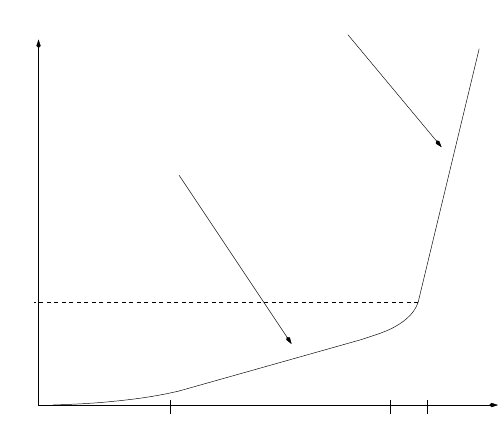
}
}
\end{figure}

The function $h_{m,\lambda,\epsilon}(r_F)$ can be viewed as a function on the whole
of $\widehat{E}$ by defining it to be zero over the region $\{r_F \leq 1\}$.
We define $B_{m,\lambda,\epsilon} := B_\lambda + h_{m,\lambda,\epsilon}(r_F)$.
We let $J$ be a Lefschetz admissible almost complex structure.
Then $(B_{m,\lambda,\epsilon},J)$ is a Lefschetz admissible family.
Here we should view $\lambda$ as any number $\geq 1$,
$m$ as large with respect to $\lambda$ and $\epsilon$ as small
with respect to $\frac{1}{m}$.
If $\lambda$ is not a multiple of $2\pi$ then all the $1$-periodic orbits
in the region $\{r_F \geq 1\}$ come in pairs
$(o_1,o_2)$ where $o_1$ is the orbit corresponding
to the minimum of $\lambda r^2$ at $p$ and
$o_2$ is any orbit of $h_{m,\lambda,\epsilon}(r_F)$ in the region
$\{r_F \geq m\}$.
There are no orbits in the region $2 \leq r_F \leq m$
because $\epsilon$ is small.
The action of the orbit $(o_1,o_2)$
described above is equal to the action of $o_2$
in $\widehat{F}$.
The action of the orbit $o_2$ is equal to $r_F h'_{m,\lambda,\epsilon}(r_F) - h_{m,\lambda,\epsilon}$.
The orbit lies in the region $\{m \leq r_F \leq m + 1\}$ (for $\lambda$ not in the action spectrum), and the length
of the smallest Reeb orbit is greater than some $\mu > 0$
so the action is greater than
$m \mu - 1$.
Choose two functions $\xi_1,\xi_2 : (0,\infty) \rightarrow \R$
so that $\xi_1(t)$ is large and $\xi_2(t)$ is large with respect to $\xi_1(t)$.
The point here is that we want $m = \xi_1(\lambda)$ and $\epsilon = \frac{1}{\xi_2(\lambda)}$,
and so for $b_{m,\lambda,\epsilon}$ to be well defined we need $\xi_1,\xi_2$
to be sufficiently large.
We define
$H'_\lambda := B_{\xi_1(\lambda),\lambda,\frac{1}{\xi_2(\lambda)}}$.
We have a short exact sequence of chain complexes:
\[0 \rightarrow CF_*^{\leq m\mu-1}(H'_\lambda,J) \rightarrow CF_*(H'_\lambda,J) \rightarrow
\frac{CF_*(H'_\lambda,J)}{CF_*^{\leq m\mu-1}(H'_\lambda,J)} \rightarrow 0.\]
We choose $m$ (and hence the function $\xi_1$) large enough with respect to $\lambda$ so that
$-\theta_E(X_{H'_\lambda}) - H'_\lambda \leq m \mu -1$ in the region $\{r_F \leq 1\}$.
This can be done because $B_{m,\lambda,\epsilon}$ does not depend on $m$ in the region $\{r_F \leq 1\}$.
This means that any orbit of $H'_\lambda$ inside $\{r_F \leq 1\}$
has action less than $m \mu -1$.
We can ensure that in the region $\{r_F \leq 1\}$
that $\frac{\partial H'_\lambda}{\partial \lambda} \geq 0$
but we cannot guarantee this in the region $\{r_F \geq 1\}$.
We wish that the continuation map joining $H'_{\lambda_1}$ to
$H'_{\lambda_2}$ for $\lambda_1 \leq \lambda_2$ through
the family $H'_\lambda$ preserves the subcomplex
$ CF_*^{\leq m\mu-1}(H'_\lambda,J)$.
All the orbits of this subcomplex are contained in
$\{r_F \leq 1\}$ (after possibly $C^1$ perturbing the function $H'_\lambda$ slightly)
and the maximum principle \cite[Lemma 7.2]{SeidelAbouzaid:viterbo} ensures that Floer trajectories
connecting these orbits stay inside this region.
Also $\frac{\partial H'_\lambda}{\partial \lambda} \geq 0$
in this region which ensures that the continuation map
sends orbits of action $\leq m \mu - 1$ to orbits
of action $\leq m \mu - 1$.
Hence the continuation map induces a morphism of short exact sequences:
\[
\xy
(0,0)*{}="A1"; (25,0)*{}="A2"; (60,0)*{}="A3";(95,0)*{}="A4"; (120,0)*{}="A5";
(0,-15)*{}="B1"; (25,-15)*{}="B2"; (60,-15)*{}="B3";(95,-15)*{}="B4"; (120,-15)*{}="B5";
%
"A1" *{0};
"A2" *{CF_*^{\leq \xi_1(\lambda_1)\mu-1}(H'_{\lambda_1},J)};
"A3" *{CF_*(H'_{\lambda_1},J)};
"A4" *{\frac{CF_*(H'_{\lambda_1},J)}{CF_*^{\leq \xi_1(\lambda_1)\mu-1}(H'_{\lambda_2},J)}};
"A5" *{0};
"B1" *{0};
"B2" *{CF_*^{\leq \xi_1(\lambda_2)\mu-1}(H'_{\lambda_2},J)};
"B3" *{CF_*(H'_{\lambda_2},J)};
"B4" *{\frac{CF_*(H'_{\lambda_2},J)}{CF_*^{\leq \xi_1(\lambda_2)\mu-1}(H'_{\lambda_2},J)}};
"B5" *{0};
%
{\ar@{->} "A1"+(2,0)*{};"A2"-(19,0)*{}};
{\ar@{->} "A2"+(19,0)*{};"A3"-(11,0)*{}};
{\ar@{->} "A3"+(11,0)*{};"A4"-(16,0)*{}};
{\ar@{->} "A4"+(16,0)*{};"A5"-(2,0)*{}};
{\ar@{->} "B1"+(1,0)*{};"B2"-(19,0)*{}};
{\ar@{->} "B2"+(19,0)*{};"B3"-(11,0)*{}};
{\ar@{->} "B3"+(11,0)*{};"B4"-(16,0)*{}};
{\ar@{->} "B4"+(16,0)*{};"B5"-(2,0)*{}};
{\ar@{->} "A2"+(0,-4)*{};"B2"+(0,4)*{}};
{\ar@{->} "A3"+(0,-4)*{};"B3"+(0,4)*{}};
{\ar@{->} "A4"+(0,-4)*{};"B4"+(0,6)*{}};
%
\endxy
\]
We also have another Hamiltonian $K_\lambda$
which is equal to $H'_\lambda$ in the region $\{r_F \leq m\}$
and equal to $\epsilon (r_F-1) + B_\lambda$ outside this region
where $m = \xi_1(\lambda)$ and $\epsilon = \frac{1}{\xi_2(\lambda)}$ as before.
The maximum principle \cite[Lemma 7.2]{SeidelAbouzaid:viterbo}
ensures that any Floer trajectory connecting orbits
inside $\{r_F \leq 1\}$ stays inside this region (as long
as we choose an appropriate $J$).
Any orbit of $K_\lambda$ also has action less than $m \mu -1$.
Hence there is a chain isomorphism
$CF_*^{\leq \xi_1(\lambda)\mu -1}(H'_\lambda,J) \cong
CF_*(K_\lambda,J)$.
The continuation maps between $CF_*(K_{\lambda_1},J)$
and $CF_*(K_{\lambda_2},J)$ are the same
as the ones between
$CF_*^{\leq \xi_1(\lambda_1)\mu -1}(H'_{\lambda_1},J)$
and $CF_*^{\leq \xi_1(\lambda_2)\mu -1}(H'_{\lambda_2},J)$
under this isomorphism.
Hence the filtered directed system
$(SH_*^{\leq \xi_1(\lambda)\mu -1}(H'_\lambda,J_\pi))$
is isomorphic to
$(SH_*(K_\lambda,J_\pi))$.
So in order to show that
$(SH_*(H'_\lambda,J_\pi))$
is isomorphic to
$(SH_*(K_\lambda,J_\pi))$
we need to show that
$\left(H_*\left(\frac{CF_*(H'_\lambda,J)}{CF_*^{\leq m\mu-1}(H'_\lambda,J)}\right)\right)$
is isomorphic to the trivial filtered directed system $(0)$
by Lemma \ref{lemma:filteredlongexactisomorphism}.

This is done as follows:
We will show that for $|\lambda_1 - \lambda_2| > 2\pi$,
the continuation map
\[\frac{CF_*(H'_{\lambda_1},J)}{CF_*^{\leq \xi_1(\lambda_1)\mu-1}(H'_{\lambda_1},J)}
\rightarrow \frac{CF_*(H'_{\lambda_2},J)}{CF_*^{\leq \xi_1(\lambda_2)\mu-1}(H'_{\lambda_2},J)}\]
is trivial.
The vector space $\frac{CF_*(H'_\lambda,J)}{CF_*^{\leq m\mu-1}(H'_\lambda,J)}$
has a basis given by orbits in the region
$\{r_F \geq m\}$.
All these orbits are contained in the region
$\pi^{-1}(U) = U \times \widehat{F}$ described earlier in this proof.
We can choose $J$ so that it splits as $j \oplus J_F$
in $U \times \widehat{F}$ where $\widehat{F} \cong \pi^{-1}(p)$.
A maximum principle ensures that any Floer trajectory
or continuation map trajectory joining orbits in
this region must stay inside this region.
If the partial Lefschetz fibration was in fact a Lefschetz fibration
(i.e. $\pi$ is well defined everywhere) then the maximum principle
needed to achieve this would be \cite[Lemma 5.2]{McLean:symhomlef}.
In general, we use Lemma \ref{lemma:lefschetzmaximumprinciple2}.
%
Because the almost complex structure
splits as a product as well as the Hamiltonian inside $U \times \widehat{F}$,
we have that $\frac{CF_*(H'_\lambda,J)}{CF_*^{\leq m\mu-1}(H'_\lambda,J)}$
is isomorphic to a tensor product:
\[CF_*(\lambda r^2,j) \otimes \frac{CF_*(H'_\lambda|_{\pi^{-1}(p)},J_F)}{CF_*^{\leq \xi_1(\lambda)\mu-1}(H'_\lambda|_{\pi^{-1}(p)},J_F)}.\]
We have that the chain map $CF_*(\lambda_1 r^2,j) \rightarrow CF_*(\lambda_2 r^2,j)$
is $0$ for $|\lambda_1 - \lambda_2| > 2\pi$ for index reasons
(see \cite[Section 3.2]{Oancea:survey}).
This implies that the continuation map
\[\frac{CF_*(H'_{\lambda_1},J)}{CF_*^{\leq \xi_1(\lambda_1)\mu-1}(H'_{\lambda_1},J)}
\rightarrow \frac{CF_*(H'_{\lambda_2},J)}{CF_*^{\leq \xi_1(\lambda_2)\mu-1}(H'_{\lambda_2},J)}\]
is zero.

Hence we have shown  $(SH_*(H_\lambda,J_\lambda))$
is isomorphic to
$(SH_*(\widehat{E},\theta_E,\lambda))$.
\qed

\bigskip
From now on the base of our partial Lefschetz fibration is the disk $\D$
so the base of $\widehat{E}$ is $\C = \widehat{\D}$.
Let $p$ be a point on $\C \setminus \D$. This is a regular value of $\pi$.
We will now put a filtration on the Floer chain complex of some half Lefschetz admissible pairs.
By Lemma 3.1 of \cite{McLean:spectralsequence}
we can deform $E$ so that
there is a point $p \in \C$ and a small disk neighbourhood
$U$ such that $\pi^{-1}(U) = \widehat{F} \times U$
with $\theta_E = r^2 d\vartheta + \theta_F$ in this region.
Here $(r,d\vartheta)$ are polar coordinates.
Let $(H_\lambda,J_\lambda)$ be a smooth family
of Lefschetz admissible pairs where $H_\lambda$ is of slope
$\lambda$ such that:
$H_\lambda = \epsilon_H h(r_F) + \epsilon r^2$ in the region $\widehat{F} \times U$
($\epsilon>0$ is some constant).
When we calculate $SH_*(H_\lambda,J_\lambda)$
we perturb $H_\lambda$ so that it is non-degenerate and of the form
$H'_F + \epsilon r^2$ on $\pi^{-1}(U)$ where $H'_F$ is also non-degenerate.
Also we split $J_\lambda$ up as a product in this region
so that $\pi^{-1}(p)$ is a holomorphic hypersurface.
We have a natural filtration on the chain complex
for  $SH_*(H_\lambda,J_\lambda)$.
It is described as follows:

Choose a closed $1$-form $\beta$ on $\widehat{E} \setminus \pi^{-1}(p)$
which is equal to $\frac{1}{2\pi}d\vartheta$ inside $(U \setminus p) \times \pi^{-1}(p)$.
This $1$-form exists because we can choose such a closed $1$-form inside
$\pi^{-1}( \C \setminus (\{p\} \cup \{\D\}) )$
so that it is exact on a small neighbourhood of $\pi^{-1}(\partial \D)$.
Hence we can extend it to the whole of $\widehat{E} \setminus \pi^{-1}(p)$.
We have a filtration
$F_0^\lambda \subset F_1^\lambda \subset \cdots$ on the chain complex $CF_*(H_\lambda,J)$
where $F_0$ is generated by orbits which either project to $p$
or such that the integral of $\beta$ over the orbit is zero or positive.
The group $F_i$ for $i > 0$ is generated by orbits $\gamma$ where $\int_\gamma \beta \geq -i$
and by elements of $F_0$.
Any Floer trajectory of $H_\lambda$ or continuation map Floer
trajectory associated to a family of $H_\lambda$'s must intersect
the fiber $\pi^{-1}(p)$ positively or must be contained
inside the fiber.
This means that the differential respects the filtration structure.
Also the continuation maps between
$CF_*(H_{\lambda_1},J) \rightarrow CF_*(H_{\lambda_2},J)$
for $\lambda_1 \leq \lambda_2$ respect this filtration structure
(i.e send elements of $F^{\lambda_1}_i$ to elements of $F^{\lambda_2}_i$).
We define $F_{-1}^\lambda := 0$.
Hence we have a filtered directed system
$\bigoplus_{i=0}^\infty H_*(F_i^\lambda / F_{i-1}^\lambda)$.

\begin{lemma} \label{lemma:lefschetzfiltereddirectedsystemequivalence}
The filtered directed systems $(V^\phi_\lambda)$
and 
$\bigoplus_{i=0}^\infty H_*(F_i^\lambda / F_{i-1}^\lambda)$ are isomorphic.
\end{lemma}
\proof of Lemma \ref{lemma:lefschetzfiltereddirectedsystemequivalence}.
The groups $H_*(F_i^\lambda / F_{i-1}^\lambda)$ along with the continuation maps between them
do not change if we choose another
pair $(H_\lambda,J_\lambda)$ as described above.
Hence we will define $H_\lambda$ as follows:
Let $h : [1,\infty) \rightarrow \R$ be a function such that $h(x) = 0$ near $x=1$
and $h,h',h'' \geq 0$ and near infinity $h'$ is constant and less than the length of the smallest
Reeb orbit of $\partial F$.
We set $H_\lambda = h(r_F) + \lambda h(r_S)$ where $r_S$ is the cylindrical coordinate
on the base $\widehat{\D}$ and we choose $h$ so that $H_\lambda^{-1}(0)$ contains $\pi^{-1}(U)$.
We have that $H_*(F_0^\lambda / F_{-1}^\lambda)$ is isomorphic to
$H^{n-*}(E)$ because this complex is generated by the constant orbits of $H_\lambda$
(all non-trivial orbits wrap around $\pi^{-1}(p)$ a non-zero number of times).
If we look at $H_*(F_i^\lambda / F_{i-1}^\lambda)$  for $i > 1$
we see that all the generators for this chain complex lie in the region $r_S \geq 1$.
The closed $1$-form $\beta$ constructed in the paragraph before this lemma satisfies
$\beta = \kappa \alpha_S$ inside the region $\{r_S \geq 1\}$ and
$\beta = \frac{1}{2\pi}d\vartheta$ inside $U$ where $\kappa > 0$ is some constant.
Let $f : (0,1) \rightarrow (0,\infty)$ be a function equal to $1$ outside a neighbourhood of $0$
and which tends to infinity as we reach $0$.
We write $f(r)$ for the function equal to $f(r)$ inside $U \times \pi^{-1}(p)$
and equal to $1$ outside this region and we define $f$ so that this function is smooth.
We also assume that $f' < 0$ near zero and $f' \leq 0$ everywhere.
We can adjust $U$ and add some small positive multiple $\tau f(r)\beta$
to $\theta_E$ so that $(\widehat{E} \setminus \pi^{-1}(p), \theta_E + f(r)\beta)$
is the completion of a Lefschetz fibration whose base is the annulus.
Here the new Lefschetz fibration map is $\pi |_{\widehat{E} \setminus \pi^{-1}(p)}$,
and we have a new cylindrical coordinate $r'_S$ where $r'_S = \frac{1}{C}r_S$
in the region $r_S \geq \frac{1}{C}$ where $C\geq 1$ is a constant so that $\{r_S \leq C\}$
contains $U$. In this region $\theta_E = (r'_S-1) (C\alpha_s) + C\alpha_S + \tau \beta + \alpha_\partial$
($\alpha_\partial = \theta_E|_{\partial_v E}$).
Also $r'_S = (f(r) + 2\pi r^2)$
in the region $U$.
We can define $H'_\lambda$ so that it is half Lefschetz admissible
for $(\widehat{E} \setminus \pi^{-1}(p), \theta_E + f(r)\beta)$
and so that $H'_\lambda = H_\lambda$ outside $U \times \pi^{-1}(p)$.
Let $\alpha \subset H_1(\widehat{E} \setminus \pi^{-1}(p))$
be the set of $H_1$ classes represented by loops which wrap a strictly negative number
of times around $p$ when we project to the base $\widehat{\D} \setminus p$.
We can ensure all orbits of $H'_\lambda$ representing classes in $\alpha$ are contained in the region
$\{r_S \geq 1\}$.
Because $H_*(F_i^\lambda / F_{i-1}^\lambda)$  
does not change for $i > 1$
 if we change $\theta_E$ inside a closed subset of the region $U \times \pi^{-1}(p)$
(assuming $\pi$ still can be made holomorphic in this region), we have that:
\[H_*(F_i^\lambda / F_{i-1}^\lambda) = SH_*^{\#,\alpha}(H'_\lambda,J')\]
where $J'$ is some Lefschetz admissible almost complex structure.
Also the associated continuation maps induced by increasing $\lambda$ commute with this isomorphism.
By Lemma \ref{lemma:monodromyfiltereddirectedsystems},
$(SH_*^{\#,\alpha}(H'_\lambda,J'))$
is isomorphic to the filtered directed system
$\left(\bigoplus_{i=0}^{\lfloor \lambda \rfloor}(HF_*(\phi,i))\right)$
where the directed system maps are the natural inclusion maps.
Summing the above discussion up we have the following isomorphism
of filtered directed systems:
\[\left(\bigoplus_{i=1}^{\infty} H_*(F_i^\lambda / F_{i-1}^\lambda)\right) \cong
\left(\bigoplus_{i=1}^{\lfloor \lambda \rfloor}(HF_*(\phi,i))\right)\]
which in turn implies
that:
$\left(\bigoplus_{i=0}^{\infty} H_*(F_i^\lambda / F_{i-1}^\lambda)\right)$
is isomorphic to
\[(V^\phi_\lambda) := (H^{n-1}(E) \oplus \bigoplus_{i=1}^{\infty} H_*(F_i^\lambda / F_{i-1}^\lambda)).\]
\qed

We now have enough ingredients to prove Theorem \ref{theorem:lefschetzbiggerthan}

\proof of Theorem \ref{theorem:lefschetzbiggerthan}.
By the proof of \cite[Lemma 2.9]{McLean:spectralsequence},
we have a constant $R$ such that the filtered directed system maps
$H_*(F^{\lambda_1}_i / F^{\lambda_1}_{i-1}) \rightarrow H_*(F^{\lambda_2}_i / F^{\lambda_2}_{i-1})$
are isomorphisms for all $\lambda_2 \geq \lambda_1 \geq Ri$.
There is also a constant $\nu > 0$
so that any $1$-periodic orbit of $H_\lambda$
wraps around $\pi^{-1}(p)$ at most $\nu \lambda$ times
(for an appropriate choice of $H_\lambda$).
Hence we have by Lemma \ref{lemma:biggerthanspectralsequence},
\[\left( H^{n-*}(E) \oplus \bigoplus_{i=1}^{\infty} H_*(F^\lambda_i / F^\lambda_{i-1})\right)\]
is bigger than $(SH_*(H_\lambda,J))$.
Hence by Lemmas \ref{lemma:inequalityinvariance} and \ref{lemma:lefschetzfiltereddirectedsystemequivalence},
we have that
$(V^\phi_{\lambda})$ is bigger than $(SH_*(H_\lambda,J))$.

By Lemmas 
\ref{lemma:halflefschetzlefschetzequivalence} and
\ref{lemma:inequalityinvariance} we get that
$(V^\phi_{\lambda})$ is bigger than
$(SH_*(\widehat{E},\theta_E,\lambda))$.
This proves the Theorem.
\qed

\subsection{Growth rates of fillings and algebraic Lefschetz fibrations}
\label{section:fillingsofalgebraiclefschetzfibrations}

The aim of this section is to prove the following Theorem:
\begin{theorem} \label{theorem:affinelefschetzupperbound}.
Any smooth affine variety $A$ is exact symplectomorphic to
some Lefschetz fibration $\pi : \widehat{E} \twoheadrightarrow \C$
with monodromy map $\phi$ such that
$\Gamma(\phi) \leq \text{dim}_\C A$.
\end{theorem}

This theorem is in fact true if we replace $A$ by $\widehat{M}$
where $\partial M$ is contactomorphic to $\partial \overline{A}$
and where we replace `Lefschetz fibration' with `partial Lefschetz fibration'.
The proof is basically the same but we will deal with $A$ for simplicity.
The choice of trivialization of the canonical bundle and of our class $b$
is the same as in section \ref{section:lefschetzfibrationsbounds}.
It turns out (using similar geometrical ideas from the proof
of Lemma \ref{lemma:partiallefschetzfibrationexistence}) that if
we have a partial Lefschetz fibration over the disk with monodromy $\psi$
then we also have an open book on the contact boundary $\partial \overline{A}$
whose monodromy is $\psi$. Hence the above theorem tells
us that $\partial \overline{A}$ admits an open book whose
monodromy map $\phi$ satisfies $\Gamma(\phi) \leq \text{dim}_\C A$.
We will prove the following stronger theorem:
\begin{theorem} \label{theorem:Lefschetzfibrationpolynomialbound}
Any smooth affine variety $A$ is exact symplectomorphic to
some Lefschetz fibration $\pi : \widehat{E} \twoheadrightarrow \C$
with monodromy map $\phi$ such that
there exists a polynomial $R$ of degree $n := \text{dim}_\C A$ such that
$|R(x)| \geq |V^\phi_x|$ for all $x$.
\end{theorem}
Here $V^\phi_x$ is the filtered directed system associated to $\phi$
as described in Section \ref{section:lefschetzfibrationsbounds}.
Before we prove Theorem \ref{theorem:Lefschetzfibrationpolynomialbound}
we need to know about convex symplectic manifolds and
algebraic Lefschetz fibrations and also prove some preliminary Lemmas.
A {\it convex symplectic manifold} is a manifold $M$ with a $1$-form $\theta_M$
such that
\begin{enumerate}
\item $\omega_M := d\theta_M$ is a symplectic form.
\item There is an exhausting function $f_M : M \rightarrow \R$
and a sequence $c_1 < c_2 < \cdots$ tending to infinity
such that $(f_M^{-1}(-\infty,c_i],\theta_M)$ is a Liouville domain
for each $i$.
\end{enumerate}
We say that $M$ is of {\it finite type} if there is a $C \in \R$
such that $(f_M^{-1}(-\infty,c],\theta_M)$ is a Liouville domain
for all $c \geq C$.
A convex symplectic manifold is said to be {\it complete}
if the $\omega_M$-dual of $\theta_M$ is an integrable vector field.

Let $(M,\theta^t_M)$ be a smooth family of convex symplectic manifolds
parameterized by $t \in [0,1]$.
This is said to be a {\it convex deformation} if
for every $t \in [0,1]$ there is a $\delta_t>0$
and an exhausting function $f_M^t$ and a sequence
of constants $c_1^t < c_2^t < \cdots$ tending to infinity
such that $( (f_M^t)^{-1}(-\infty,c_i^t],\theta^s_M)$ is a Liouville domain
for each $s \in [t- \delta_t,t+\delta_t]$ and each $i$.
We do not require that $f_M^t$,$c_i^t$,$\delta_t$
smoothly varies with $t$. In fact it can vary in
a discontinuous way with $t$.

\begin{lemma} \label{lemma:convexsubmanifolds}
If we have two finite type convex symplectic manifolds
$B,B'$ such that
\begin{enumerate}
\item
$B$ is a codimension $0$ exact symplectic submanifold of $B'$.
\item $df_{B'}(X_{\theta_{B'}}) > 0$ outside some closed subset of $B$.
\end{enumerate}
Then $B$ is convex deformation equivalent to $B'$.
Here $f_{B'}$ is an exhausting function and
$X_{\theta_{B'}}$ is the $d\theta_{B'}$-dual of $\theta_{B'}$.
\end{lemma}
\proof of Lemma \ref{lemma:convexsubmanifolds}.
Let $P'$ be the regular hypersurface $f_{B'}^{-1}(C)$ for some $C$ where
$(f_{B'}^{-1}(-\infty,c],\theta_{B'})$ is a Liouville domain for all $c \geq C$.
Flow $P'$ backwards along $X_{\theta_{B'}}$ to a new hypersurface $P$
contained in $B$.
Let $f_t : B' \rightarrow [0,\infty)$, $t \in [0,1]$
be a smooth family of functions
with the following properties:
\begin{enumerate}
\item $f_0 > 0$ inside $B$ and $f_0 = 0$ outside $B$.
\item $f_1 > 0$ everywhere.
\item $f_t X_{\theta_{B'}}$ is integrable for all $t$
\item $f_t = f_0$ on a neighbourhood of $P$.
\end{enumerate}
We can construct a smooth family of embeddings $\iota_t : B \rightarrow B'$
as follows:
First $\iota_t$ is the identity map on
the compact submanifold $Q$
whose boundary is $P$.
We define $\iota_t(x)$ for $x$ outside $Q$ as follows:
we first flow $x$ via $-f_0 X_{\theta_{B'}}$
for some time $t_x$ until it hits $P$.
Then we flow $x$ forward via $f_t X_{\theta_{B'}}$
for time $t_x$.
The result of this flow is our definition of $\iota_t(x)$.

The vector field $f_t X_{\theta_{B'}}$ restricted to the image
$\iota_t(B)$ is complete.
We can flow $P$ along this vector field so that
the closure of the complement of $Q$ is diffeomorphic to
$[1,\infty) \times P$ with $X_{\theta_{B'}}$
equal to $\frac{1}{f_t} \frac{\partial}{\partial r_t}$
where $r_t$ parameterizes $[1,\infty)$.
We define an exhausting function $g_t : \iota_t(B) \rightarrow \R$
to be equal to our coordinate $r_t$ near infinity and anything
we like elsewhere.
Hence $dg_t(X_{\theta_{B'}}) > 0$ outside a compact subset of $\iota_t B$.
This means that $(\iota_t(B),\theta_{B'})$ is a smooth family of finite type convex
symplectic manifolds such that $\iota_0(B)$ is exact symplectomorphic
to $B$ and $\iota_1(B)$ is exact symplectomorphic to $B'$.
Hence $(B,\theta_{B'})$ is convex deformation equivalent
to $(B',\theta_{B'})$.
Also $(B,\theta_B)$ is exact symplectomorphic to $(B,\theta_{B'})$
and hence by applying \cite[Lemma 8.3]{McLean:symhomlef} 
we get that $(B,\theta_B)$ is convex deformation equivalent to
$(B,\theta_{B'})$.
Hence $(B,\theta_B)$ is convex deformation equivalent to
$(B',\theta_{B'})$.
\qed

From Section \ref{section:liouvilldomaindefinition}
we know that a smooth affine variety $A$ has
a natural symplectic structure induced by some
embedding into $\C^N$.
Because $A$ is the completion of some natural Liouville domain
$\overline{A}$ we have that $A$ 
has the structure of a finite type convex symplectic manifold.
%
Let $X$ be a compactification of $A$ by a smooth normal
crossing divisor $D$ and let $t$ be a section of $L$ such that $t^{-1}(0)$ is smooth, reduced
and transverse to all the strata of $D$. Then $p : A \rightarrow \C$,
$p := t / s$ is called an {\it algebraic Lefschetz fibration}
if all the singularities of $p$ are non-degenerate and on distinct fibers.
The map $p$ is an algebraic Lefschetz fibration for a generic section $t$.

Let $\pi : \widehat{E} \rightarrow \widehat{S}$ be a partial Lefschetz fibration.
The region $\pi^{-1}(\{r_S \geq 1\})$ is called the {\it Lefschetz cylindrical end}
of $\widehat{E}$. Here $r_S$ is the cylindrical coordinate of $\widehat{S}$.
This region is diffeomorphic to $[1,\infty) \times (\pi^{-1}(\partial S))$ and
$\pi$ maps this region to $[1,\infty) \times \partial S$ via the map
$(\text{id},\pi|_{\pi^{-1}(\partial S)})$.
The $1$-form $\theta_E$ is equal to $(r_S - 1) \alpha_S + \alpha_\partial$
where $\alpha_\partial = \theta_E|_{\pi^{-1}(\partial S)}$, $\alpha_S = \theta_S|_{\partial S}$.
If we have some connected component $\partial'S$ of $\partial S$
then the subset $\pi^{-1}([1,\infty) \times \partial'S)$
is called a {\it component of the Lefschetz cylindrical end}.
We also have a region in $\widehat{E}$ diffeomorphic to
\[\widehat{E}_h := \widehat{S} \times [1,\infty) \times \partial F\]
where $\theta_E = \theta_S + r_F \alpha_F$.
Here $r_F$ parameterizes $[1,\infty)$, $\alpha_F$
is the natural contact form on the boundary of $F$ and $\pi$ is the projection map
to $\widehat{S}$.

Let $\pi_1 : \widehat{E}_1 \rightarrow \widehat{S}_1$,
$\pi_2 : \widehat{E}_2 \rightarrow \widehat{S}_2$ be two
partial Lefschetz fibrations.
Let $\partial' S_1$ be a connected component of $\partial S_1$
and $\partial' S_2$ a connected component of $\partial S_2$.
We say that {\it $\widehat{E}_1$ and $\widehat{E}_2$
have two identical Lefschetz cylindrical end components}
if there exist $\partial' S_1$ and $\partial' S_2$ as above
and a diffeomorphism $\Phi$
from $\pi_1^{-1}(\partial' S_1 \times [1,\infty))$ to
$\pi_2^{-1}(\partial' S_2 \times [1,\infty))$ which is also a
map of fibrations from $\pi_1$ to $\pi_2$ covering an orientation
preserving diffeomorphism from $\partial' S_1 \times [1,\infty)$
to $\partial' S_2 \times [1,\infty)$.
This diffeomorphism $\Phi$ must pull back
$\theta_{E_1}$ to $\theta_{E_2} + \pi_{E_2}^* \beta$
where $\beta$ is a $1$-form on the base $\widehat{S}_2$.
This means that these regions have identical parallel transport
maps.

\begin{theorem} \label{theorem:algebraiclefschetzfibration}
For every algebraic Lefschetz fibration $p : A \rightarrow \C$,
where $q \in \C$ is a regular value,
there is a Lefschetz fibration $\pi : \widehat{E} \rightarrow \C$
such that $A$ is convex deformation equivalent to $\widehat{E}$
and $A \setminus p^{-1}(q)$ is convex deformation equivalent
to a Lefschetz fibration $\pi' : \widehat{E'} \rightarrow \C^*$
such that $\pi,\pi'$ have two identical Lefschetz cylindrical end components.
Also the fiber $\pi^{-1}(q)$ is convex deformation equivalent to $p^{-1}(q)$.
\end{theorem}
We will prove this Theorem in Section \ref{section:algebraicsymplecticlefschetzfibration}.
We recall the definition
of $P$-bounded from \cite[Section 6]{McLean:affinegrowth}.
Suppose we have a Hamiltonian $H : W \rightarrow \R$, a function $P : \R \rightarrow \R$
and a small open neighbourhood ${\mathcal N}$ of $H^{-1}(0)$ such that:
\begin{enumerate}
\item $H$ satisfies the Liouville vector field property.
\item For every $\lambda \in (0,\infty)$ outside some discrete subset,
there is a $C^2$ small perturbation $H_\lambda$
of $\lambda H$ such that all the $1$-periodic orbits
of $H_\lambda$ inside ${\mathcal N}$ are non-degenerate and the number of
such orbits is bounded above by $P(\lambda)$.
\end{enumerate}
If $H$,$P$ satisfy these properties then we say that $H$
is $P$-bounded.

We have the following Lemma:
\begin{lemma} \label{lemma:numberoforbitsbound}
Suppose that $H$ is $P$ bounded. Also we assume
$-X_H(\theta_H) - H \geq 0$ such that this function is greater than some constant
$\delta_H > 0$ outside some compact set.
Here $\theta_H$ is a $1$-form such that $\theta_H - \theta_M$ is exact.
Let ${\mathcal N}'$ be a small neighbourhood of $H^{-1}(0)$.
Then there exists a constant $\Delta_H > 0$ such that:
For every $\lambda \geq 1$ outside some discrete subset,
there is a $C^2$ small perturbation of $H_\lambda$
of $\lambda H$
such that $H_\lambda = H$ outside ${\mathcal N}'$ and the number of $1$-periodic orbits of $H_\lambda$
of action in $[0,\Delta_H \lambda]$ is $\leq P(\lambda)$.
Also all these orbits are non-degenerate.
\end{lemma}
We will omit the proof of this Lemma as the proof is contained
inside the proof of
\cite[Lemma 6.6]{McLean:affinegrowth}.

\begin{lemma} \label{lemma:extending0}
Suppose that $H$ is $P$ bounded.
Then for every compact set $K$, there exists a Hamiltonian
$H_K$ which is $P$ bounded with $H_K^{-1}(0)$
containing $K$.
\end{lemma}
\proof of Lemma \ref{lemma:extending0}.
We assume that $H$ is a Hamiltonian on $\widehat{N}$ for some Liouville domain $N$.
Let $\theta_H$ be the $1$-form on $\widehat{N}$
ensuring that $H$ satisfies the Liouville vector field property.
This means that $d\theta_H = \omega_N$ and there is an exhausting function $f_H$
with $df_H(X_{\theta_H}) > 0$ outside a closed
subset of the interior of $H^{-1}(0)$.
Here $X_{\theta_H}$ is the $\omega_N$-dual of $\theta_H$.
If $X_{\theta_H}$ was integrable then we could use the flow
of this vector field to make the zero set of $H$ as large as we
like and then we would have proven this Lemma.
The problem is that it may not be integrable.
We will modify $H$ and $\theta_H$ so that $\theta_H$ becomes integrable.
Let $C>1$ such that the interior of
$N' := f_H^{-1}(-\infty,C]$ contains $H^{-1}(0)$.
We have that $(N',\theta_H)$ is a Liouville domain.
Also the interior $(N')^0$ of $N'$ with $1$-form $\theta_H$
has the structure of a convex symplectic manifold
as follows:
we have a collar neighbourhood $(1-\epsilon,1] \times \partial N'$
with $\theta_H = r \alpha_{N'}$ where $\alpha_{N'}:=\theta_H|_{\partial N'}$
given by flowing $\partial N'$ backwards along $X_{\theta_H}$.
Let $g : (1-\epsilon,1)$ be a function equal
to $0$ near $1-\epsilon$ and which tends to infinity
as we approach $1$. We also require that its derivative
is positive near $1$.
We define $f_{(N')^0}$ to be equal to $g(r)$ on this collar neighbourhood
and zero elsewhere.
Then $((N')^0,\theta_H)$ has the structure of a finite type
convex symplectic manifold with exhausting function $f_{(N')^0}$.

Both $\widehat{N}$ and $\widehat{N'}$
are convex deformation equivalent to $((N')^0,\theta_H)$
by Lemma \ref{lemma:convexsubmanifolds}.
Hence $\widehat{N}$ and $\widehat{N'}$
are convex deformation equivalent and hence
by \cite[Corollary 8.6]{McLean:affinegrowth}
they are exact symplectomorphic.
Let $\Phi$ be this exact symplectomorphism.

Extend $H|_{(N')^0}$ to any positive Hamiltonian $H'$
on $\widehat{N'}$ which is bounded below by some positive constant
near infinity.
We write $\theta_{N'}$ for the natural $1$-form on $\widehat{N}'$
so that $\theta_{N'} = \theta_H$ on $N'$.
There is a function $f_{N'}$ 
such that $df_{N'}(X_{\theta_{N'}}) > 0$ outside a closed subset
of the interior of $(H')^{-1}(0)$.
Hence because $H = H'$ and $\theta_H = \theta_{N'}$ on a small neighbourhood of
$H^{-1}(0) = (H')^{-1}(0)$ we have that $H'$
must satisfy the Liouville vector field property
and be $P$ bounded.

Let $\phi_s$ be the flow of the vector field $X_{\theta_{N'}}$.
We define $H_s' := H' \circ \phi_{-s}$.
This still satisfies the Liouville vector field
property. Also the flow of any perturbation $K'$
of $H'$ is the same as the flow
of $e^{s}K' \circ \phi_{-s}$ pushed forward by $\phi_{-s}$.
This means that $H' \circ \phi_{-s}$ is $P$ bounded.
For $s$ large enough we have that $(H_s')^{-1}(0)$ contains
$\Phi^{-1}(K)$.

Because the property of being $P$-bounded is invariant under
exact symplectomorphisms we have that
$H_K := H'_s \circ \Phi$ is also $P$ bounded on $\widehat{N}$
and its zero set contains $K$.
\qed

Let $\pi : \widehat{E} \twoheadrightarrow \widehat{S}$ be a partial Lefschetz fibration
with fiber $\widehat{F}$ and base $\widehat{S}$.

\begin{lemma} \label{lemma:boundonlefschetzfibration}
Let $P,Q$ be functions from $[1,\infty)$ to $\R$
such that
\begin{enumerate}
\item  $H_F$ is a Hamiltonian on $\widehat{F}$ which is $P$ bounded.
\item $(SH_*^{\#}(\widehat{E}))$ is isomorphic to $(W_\lambda)$
with $|W_\lambda| \leq Q(\lambda)$.
\end{enumerate}
Let $\phi$ be the monodromy map around one of the components $A_S$ of $\partial S$.
If $\alpha \subset H_1(E)$ is the set of $H_1$ classes represented by loops
which wrap around $A_S$ non-trivially after projecting by $\pi$ then we assume that
$0 \notin \alpha$ and we also assume that the base $S$ is non-contractible.
Then there are constants $C,\kappa_1,\kappa_2 \geq 1$ such that
$|V^\phi_\lambda| \leq C \lambda P(\kappa_1 \lambda) + Q(\kappa_2 \lambda)$.
\end{lemma}
\proof of Lemma \ref{lemma:boundonlefschetzfibration}.
Our Lefschetz fibration is $\pi : \widehat{E} \rightarrow \widehat{S}$
and our monodromy map is $\phi$ and our fiber is $\widehat{F}$.
We will prove this lemma by using Lemma \ref{lemma:chaincomplexhomologyupperbound}
with the chain complex associated to a Hamiltonian roughly of the form
$h(r_S) + H_F$. The statement of Lemma \ref{lemma:chaincomplexhomologyupperbound}
has a chain complex of the form $A_\lambda \oplus B_\lambda$ where in our case
$A_\lambda$ correponds to orbits of this Hamiltonian inside $\widehat{E}_h$
and $B_\lambda$ corresponds to orbits away from this region.
We will prove this Lemma in $2$ steps. In Step $1$ we will construct the Hamiltonian
and also construct explicit perturbations of it
and in Step $2$ we will apply Lemma \ref{lemma:chaincomplexhomologyupperbound}
to this Hamiltonian.

{\it Step 1:}
Let $\text{supp}(\phi) \subset \widehat{F}$ be the relatively compact set where
$\phi$ is not the identity map.
By enlarging $F$ we can assume that it contains $\text{supp}(\phi)$.
By Lemma \ref{lemma:extending0} we can assume that $H_F^{-1}(0)$
contains $F$.
By \cite[Lemma 4.2]{McLean:affinegrowth} there is a growth rate admissible
pair $(H^p_F,J_F)$ such that:
\begin{enumerate}
\item \label{item:neighbourhoodproperty}
$H^p_F = H_F$ on a small neighbourhood of $H_F^{-1}(0)$ and $(H^p_F)^{-1}(0) = H_F^{-1}(0)$.
\item $-\theta_F(X_{H^p_F}) - H^p_F \geq 0$ everywhere and this function is positive
when $H^p_F$ is positive.
\item There is a constant $\delta^p_H > 0$ such that $-\theta_F(X_{H^p_F}) - H^p_F > \delta^p_H$
outside a small neighbourhood of $H_F^{-1}(0)$.
\end{enumerate}
By these properties we basically get that $H^p_F$ is $P$ bounded
because $H_F$ is $P$ bounded.
So from now on instead of writing $H^p_F$ we will assume that $H_F = H^p_F$.
Let $\pi_2 : \widehat{E}_h \rightarrow [1,\infty) \times \partial F$
be the natural projection map.
We view $\pi_2^* H_F$ as a function on $\widehat{E}$
by defining it to be zero away from $\widehat{E}_h$.
We let $h(x)$ be a function equal to $0$ near $x=1$
and equal to $x-1$ for $x \geq 2$ with $h,h',h'' \geq 0$.
We define $H_E := \pi^* h(r_S) + \pi_2^* H_F$.
Let $J_S$ be a complex structure on the base $\widehat{S}$ so that
it is admissible.
We let $J_E$ be an almost complex structure on $\widehat{E}$
compatible with the symplectic form so that
$J_E = J_F + J_S$ inside $\widehat{E}_h$
and so that $\pi$ is $(J_E,J_S)$ holomorphic outside some compact set.
This pair is growth rate admissible by Lemma \ref{lemma:welldefinedgrowthrateadmissiblehamiltonian}.
We will now construct an explicit perturbation of $(\lambda H_E,J_E)$.

Let $\Delta_{H_F}$ be the constant from Lemma \ref{lemma:numberoforbitsbound}
so that there is a small perturbation $H^\lambda_F$ of $\lambda H_F$ such that all the orbits
of action in $[0,\Delta_{H_F} \lambda]$
are non-degenerate and the number of them is bounded above by $P(\lambda)$.
Let $\delta_\lambda>0$ be a small constant such that
$lH_F$ has no orbits in the region $H_F^{-1}(0,\delta_\lambda)$
for all $0 < l \leq \lambda$.
By the bounded below property we can assume that this region is relatively compact
for small enough $\delta_\lambda$.
Let $\nu_\lambda : \R \rightarrow \R$ be a bump function with non-negative derivative
such that $\nu_\lambda(x) = 0$ near $x = 0$
and $\nu_\lambda(x) = 1$ for $x \geq \delta_\lambda$.
For a small enough perturbation $H^\lambda_F$
we have that all the orbits of $\nu(H_F) H^\lambda_F$ whose action is in $[0,\lambda \Delta_{H_F}]$
and that are disjoint from $\nu(H_F)^{-1}(0)$ are non-degenerate.
The point is that if there were any orbits intersecting the region where $d(\nu(H_F)) \neq 0$
then a compactness argument would imply that $lH_F$ would have an orbit in this region for $0 \leq l \leq \lambda$
which is impossible.
Let $\theta_{H_F}$ be the $1$-form enabling $H_F$ to satisfy the Liouville vector field property.
Let $H_F^{-1}(c)$ ($c>0$) be a small regular level set of $H_F$ contained in $\nu(H_F)^{-1}(0)$.
This is transverse to the Liouville vector field $X_{\theta_{H_F}}$,
so it has a small neighbourhood diffeomorphic to $(1-\epsilon,1+\epsilon) \times H_F^{-1}(c)$
where $\theta_{H_F} = r_{H_F} \theta_{H_F}|_{H_F^{-1}(c)}$ where $r_{H_F}$
parameterizes the interval $(1-\epsilon,1+\epsilon)$.
We can modify $J_F$ to $J_F^\lambda$ so that it is cylindrical in the region
$(1-\epsilon,1+\epsilon) \times H_F^{-1}(c)$
and so that $(\lambda H_F,J^\lambda_F)$ is still growth rate admissible
(i.e. still satisfies the maximum principle).
We also assume that $J^\lambda_F = J_F$ inside $F$
and also outside a small neighbourhood of $H_F^{-1}(-\infty,\delta_\lambda]$.

We assume that $h$ has the additional property that
there is some constant $C_h > 0$ so that
$h''(x) > 0$ if and only if $h(x)>0$ and $x < C_h$.
This ensures that the $1$-periodic orbits
of $\lambda h(r_S)$ on the base $S$ come in $S^1$ families
which are Morse-Bott non-degenerate for all $\lambda$ outside the period spectrum of $\partial S$.
Hence by \cite{CieliebakFloerHoferWysocki:SymhomIIApplications}
we can perturb these families into two orbits.
This means that there is a time dependent perturbation of $\lambda h(r_S)$
into $h^\lambda_S : S^1 \times S \rightarrow \R$ where the number of $1$-periodic orbits is bounded above
by $C_S \lambda$ where $C_S$ is some constant.
Let $K := \pi^* h^\lambda_S + \pi_2^* H^\lambda_F$.
All the orbits of $K$ in the region $\pi_2^* \left(\nu(H_F)^{-1}(0,\infty) \right)$
are non-degenerate because the Hamiltonian splits up as a product there and $H_F^\lambda$
is non-degenerate in this region.
The only degenerate orbits of $K$ inside $\widehat{E}_h$ are in the region $K^{-1}(0)$,
We can perturb $H^\lambda_F$ slightly again so that it is a very small increasing function of
$r_{H_F}$ near $H_F^{-1}(c)$ and so that it has no additional orbits
(i.e. this perturbation only removed orbits in the region where $\nu(H_F) = 0$).
The reason why we need to do this is that we wish to use a maximum principle later.
We also assume that all of its orbits are non-degenerate in the region $H_F^{-1}(c,\infty)$.
We can also ensure that the only constant orbits are in the region $\{H^\lambda_F = 0\}$ as well.
We now perturb $K$ outside the region $\pi_2^* H_F^{-1}(c,\infty)$
to a Hamiltonian $H^\lambda_E$ which only has non-degenerate orbits
(see \cite[Lemma 2.2]{McLean:affinegrowth}).
We define $J^\lambda_E$ to be equal to 
some Lefschetz admissible almost complex structure $J_E$ away from $\widehat{E}_h$ and equal to $J^\lambda_F \oplus J_S$
inside $\widehat{E}_h$.
By using the maximum principle \cite[Lemma 7.2]{SeidelAbouzaid:viterbo}
we have that any Floer trajectory connecting $1$-periodic orbits of
$H^\lambda_E$ away from $\pi_2^* H_F^{-1}(c,\infty)$ must stay away from $\pi_2^* H_F^{-1}(c,\infty)$.
Because $h,h',h'' \geq 0$ and $-\theta_{H_F}(X_{H_F}) - H_F \geq 0$
the perturbation $H^\lambda_E$ can be made so that all of its orbits
have action $\geq 0$ (we can do this by subtracting a small constant).

{\it Step 2:}
Let $B'$ be the $\K$ vector space generated by $1$-periodic orbits
of $H^\lambda_E$ away from $\pi_2^* H_F^{-1}(c,\infty)$
and $A'$ the vector space generated by orbits inside $\pi_2^* H_F^{-1}(c,\infty)$ all
of action inside $[0,\lambda \Delta_{H_F}]$.
The rank of $A'$ is bounded above by $(C_S \lambda)P(\lambda)$.
The chain complex for $H^\lambda_E$ as a vector space is
$A' \oplus B'$ and the differential is a matrix
\[\left( \begin{array}{cc}
\partial'_a & \partial'_{b^a} \\
\partial'_{a^b} & \partial'_b
\end{array} \right)\]
We have that $\partial_b^2 = 0$ as a maximum principle ensures that all Floer trajectories connecting orbits away from
$\pi_2^* H_F^{-1}(c,\infty)$ stay away from $\pi_2^* H_F^{-1}(c,\infty)$ so
$(B,\partial_b)$ is a chain complex.
Note that $SH_*^{[0,\lambda \Delta_{H_f}]}(\lambda H_E,J_E)$ has the structure of a filtered directed system
where the filtered directed system maps are induced by continuation maps and natural inclusion maps
(see the discussion before \cite[Lemma 6.5]{McLean:affinegrowth}).
By taking the direct limit over all approximating pairs $(H^\lambda_E,J^\lambda_E)$
we get that $SH_*^{[0,\lambda \Delta_{H_f}]}(\lambda H_E,J_E)$ is generated by a chain complex
$A \oplus B$ with differential
\[\left( \begin{array}{cc}
\partial_a & \partial_{b^a} \\
\partial_{a^b} & \partial_b
\end{array} \right)\]
such that the rank of $A$ is bounded above by
$(C_S \lambda) P(\lambda)$
and such that $\partial_b^2 = 0$.
The filtered directed system maps of $(SH_*^{[0,\lambda \Delta_{H_f}]}(\lambda H_E,J_E))$
are of the form
\[ \left( \begin{array}{cc}
a_{\lambda_1,\lambda_2} & b^a_{\lambda_1,\lambda_2} \\
a^b_{\lambda_1,\lambda_2} & b_{\lambda_1,\lambda_2}
\end{array} \right).\]
Similar reasoning and
the maximum principle \cite[Lemma 7.2]{SeidelAbouzaid:viterbo}
tells us that $b_{\lambda_1,\lambda_2}$
commutes with $\partial_b$.

Because $(H_E,J_E)$ is growth rate admissible we have by \cite[Lemma 6.5]{McLean:affinegrowth} that
$(SH_*^{[0,\lambda \Delta_{H_f}]}(\lambda H_E,J_E))$ is isomorphic as a filtered directed system to
$(W_\lambda)$
with $|W_\lambda| \leq Q(\lambda)$.
Also $(H_*(B,\partial_b))$ is a filtered directed system with filtered
directed system maps $b_{\lambda_1,\lambda_2}$.
Hence by Lemma \ref{lemma:chaincomplexhomologyupperbound},
we have that $(H_*(B,\partial_b))$ is isomorphic
to some filtered directed system
$W'_\lambda$ with $|W'_\lambda| \leq 6(C_S \lambda)P(\lambda) + Q(C'\lambda)$
where $C' \geq 1$ is some constant.

Let $b(r_F)$ be a function which is equal to $0$ near $r_F = 1$ and away from
$\widehat{E}_h$ and has very small slope when $b(r_F) > 0$ so that
it has no $1$-periodic orbits in this region.
Because $H_E = h(r_S)$ away from $\widehat{E}_h$ we can use similar reasoning
and the maximum principle \cite[Lemma 7.2]{SeidelAbouzaid:viterbo}
to ensure that $H_*(B,\partial_b)$
is equal to $SH_*^{[0,\lambda \Delta_{H_f}]}(\lambda h(r_S) + b(r_F), J_E)$
(we can perturb $(\lambda h(r_S) + b(r_F),J_E)$ and $(H_E,J_E)$
so that their chain complexes are identical away from $\pi_2^* H_F^{-1}(c,\infty)$).
Similarly the continuation maps of the respective directed systems
are identical (again only when considering Floer continuation map trajectories
away from $\pi_2^* H_F^{-1}(c,\infty)$).
Hence $(H_*(B,\partial_b))$ is isomorphic as a filtered directed system
to $(SH_*^{[0,\lambda \Delta_{H_f}]}(\lambda h(r_S) + b(r_F), J_E))$.
Basically by \cite[Lemma 6.5]{McLean:affinegrowth}
this is in turn isomorphic as a filtered directed system to
$(SH_*^{\#}(\lambda h(r_S) + b(r_F), J_E))$.
Technically  \cite[Lemma 6.5]{McLean:affinegrowth}
states that our Hamiltonians must be of the form $\lambda H$ but ours are of the form
$\lambda h(r_S) + b(r_F)$ but this does not matter as the same proof holds.

By Lemma \ref{lemma:monodromyfiltereddirectedsystems} we can ensure that
$(SH_*^{\#,\alpha}(\lambda h(r_S) + b(r_F), J_E))$
is isomorphic to $(V^\lambda_\phi)$ as filtered directed systems
where $\alpha$ is the subset of $H_1(\widehat{E})$
represented by loops which project under $\pi$
to loops wrapping around our chosen component of $\partial S$.
Here $SH_*^{\#,\alpha}$ means we restrict to orbits which represent
classes contained in $\alpha$.
Hence by Lemma \ref{lemma:inequalityinvariance} we have that $(W'_\lambda)$ is bigger than $V^\lambda_\phi$.
This means that $V^\lambda_\phi$ is isomorphic to a filtered directed system
$(W_\lambda)$ with $|W_\lambda| \leq |W'_\lambda|$ by Lemma
\ref{lemma:biggerthanrank}.
Because all the filtered directed system maps of $V^\lambda_\phi$
are injective, we have that $|V^\lambda_\phi| \leq |W'_{C_W \lambda}|$
for some constant $C_W$.
This implies that
\[|V^\lambda_\phi| \leq 6(C_S C_W \lambda)P(C_W \lambda) + Q(C'C_W\lambda)\]
This proves the Lemma where our constants satisfy
$C = 6C_SC_W$, $\kappa_1 = C_W$ and $\kappa_2 = C'C_W$.
\qed


\proof of Theorem \ref{theorem:Lefschetzfibrationpolynomialbound}.
Let $p : A \rightarrow \C$ be an algebraic Lefschetz fibration.
By Theorem \ref{theorem:algebraiclefschetzfibration},
there are Lefschetz fibrations
$\pi : \widehat{E} \rightarrow \C$ and
$\pi' : \widehat{E}' \rightarrow \C^*$
such that:
\begin{enumerate}
\item
$\widehat{E}$ is convex deformation equivalent to $A$.
\item
$\widehat{E}'$ is convex deformation equivalent to $A \setminus p^{-1}(0)$
where $0$ is without loss of generality as regular value of $p$.
\item
These Lefschetz fibrations have identical cylindrical ends and hence
their fibers are exact symplectomorphic to each other.
Let $\phi$ be the monodromy map around these cylindrical ends.
\item
Any regular fiber $\widehat{F}$ of $\pi$ and $\pi'$ is exact symplectomorphic to a smooth affine variety.
\end{enumerate}

By \cite[Theorem 6.3]{McLean:affinegrowth}
there is a polynomial $Q$ of degree $n :=  \text{dim}_\C A$
and a filtered directed system $(W_\lambda)$ such that
$(SH_*^{\#}(\widehat{E}',\lambda))$ is isomorphic to $(W_\lambda)$
and $|W_\lambda| \leq Q(\lambda)$.
Also by the proof of  \cite[Theorem 6.3]{McLean:affinegrowth}
we have a polynomial $P$ of degree
$n-1$ and a Hamiltonian $H_F$  on $\widehat{F}$ which is $P$-bounded.
By Lemma \ref{lemma:boundonlefschetzfibration}
we have that
$|V^\phi_\lambda|$ is bounded above by
a polynomial of degree $n$ in $\lambda$.
This proves the Theorem.
\qed

\section{Smooth affine varieties with subcritical handles attached} \label{section:affinesubcriticalhandles}

The aim of this section is to prove
Theorem \ref{theorem:subcriticalhandleattachaffinevariety}.
Here is the statement of this theorem:
{\it Suppose that we attach a series of subcritical handles to an algebraic Stein domain
$\overline{A}$ to create a Liouville domain $N'$.
Let $N''$ be any Liouville domain whose boundary is contactomorphic
to $\partial N'$.
If $M$ is a Liouville domain with $\widehat{M}$
symplectomorphic to $\widehat{N''}$
then $\Gamma(M) \leq \text{dim}_\C A$.}

We will assume that $\overline{A}$ is the algebraic Stein domain
obtained from a smooth affine variety $A$ and some embedding of it
into $\C^N$.
Let $F'$ be a Liouville domain and $F$ a Liouville subdomain.
Let $\phi' : F' \rightarrow F'$ be an exact symplectomorphism which
is the identity on the closure of $F' \setminus F$.
We define $\phi : F \rightarrow F$ to be equal to
$\phi'|_F$.

\begin{lemma} \label{lemma:trivialextensionoffloerhomology}
There exists a constant $\kappa$ so that
for all $k \in \N$,
\[\left|  \text{rank}(HF_*(\phi,k)) - \text{rank}(HF_*(\phi',k))\right| \leq \kappa.\]
\end{lemma}
\proof of Lemma \ref{lemma:trivialextensionoffloerhomology}.
Note that we can complete $F$ inside $\widehat{F}'$ by \cite[Lemma 3]{SS:rama},
hence we have a natural exact symplectic embedding of $\widehat{F}$
inside $\widehat{F}'$ extending the embedding of $F$ into $\widehat{F}'$.

We let $0<\delta<1$ be smaller than the length of the smallest Reeb orbit
of $\partial F$ and $\partial F'$.
We let $r_F$ and $r_{F'}$ be the cylindrical coordinates
for $\widehat{F}$ and $\widehat{F}'$ respectively.
We view $r_F$ as a function on a subset of $\widehat{F}'$ corresponding
to the cylindrical end of $\widehat{F}$.
Let $f_\delta : [0,\infty) \rightarrow [0,\infty)$ be a smooth function satisfying
$f_\delta',f_\delta'' \geq 0$.
We also assume that $f_\delta(x) = 0$ for $x \leq 1$ and $f_\delta(x) =  \delta (x-1)$
for $x \geq 2$.
Let $g_b : \widehat{F'} \rightarrow \R$ be a smooth family of functions
parameterized by $b \in (0,\infty)$ such that:
$g_b = 0$ inside $F$, $g = f_\delta(r_F)$ in the region $\{1 \leq r_F \leq b\}$,
and $g_b = \delta r_{F'} + \delta b$ for $r_{F'}$ large enough.
We assume that the derivatives of $g_b$ are so small everywhere else that
the only $1$-periodic orbits are critical points of $g_b$.
Note that $g_b$ has no $1$-periodic orbits in the region
$\{1 < r_F \leq b\}$.
We also assume that the number of critical points of $g_b$
in the region $\widehat{F}' \setminus F$
is bounded above by $\frac{1}{2} \kappa$ where $\kappa$ is some constant
independent of $b$
and all of these critical points are non-degenerate.
We also assume that $g_b \geq \delta (b-1)$ outside
$F \cup \{1 < r_F \leq b\}$ and $g_b=0$ inside $F$.

Let $M_{\phi}$ and $M_{\phi'}$ be the mapping tori
for $\phi$ and $\phi'$. Let $\alpha$ and $\alpha'$
be their respective contact forms.
We also have Lefschetz fibrations
$(0,\infty) \times M_\phi$
and
$(0,\infty) \times M_{\phi'}$
with associated Liouville forms
$r_\phi d\vartheta + \alpha$ and $r_{\phi'} d\vartheta +  \alpha'$
where $r_\phi$ and $r_{\phi'}$ parameterizes $(0,\infty)$.
Here $\vartheta$ by abuse of notation is the pullback of the angle
coordinate of $S^1$ by $\pi$ or $\pi'$.
The fibration $(0,\infty)\times M_\phi$
is naturally a subfibration of $(0,\infty) \times M_{\phi'}$
because $\phi'$ restricted to $F$ is $\phi$.
Let $p_\lambda : (0,\infty) \rightarrow \R$
be defined as follows:
The derivative $p'_\lambda$ is small and positive near zero,
and constant and equal to $\lambda$ near infinity.
Also we assume that $p'_\lambda,p''_\lambda \geq 0$.
If $0 < p_\lambda' < \lambda$ then we need $p''_\lambda > 0$.
This condition is useful later to ensure that certain periodic orbits
are Morse Bott non-degenerate.
We can pull back $p_\lambda$ to $(0,\infty) \times M_{\phi'}$
via the natural projection map to $(0,\infty)$.
By abuse of notation we will call this new function $p_\lambda$.
The Liouville form
$\alpha_{\phi'}$ on the region $(0,\infty) \times S^1 \times (\widehat{F'} \setminus F)$
is a product $r_{\phi'} d\vartheta + \theta_{F'}$.
Also we define $g_b$ as a function on this product (by pulling it
back via the projection map to $(\widehat{F'} \setminus F)$)
and then we extend it by $0$ to the whole of $(0,\infty) \times M_{\phi'}$.
We define $K_\lambda^b : (0,\infty) \times M_{\phi'} \rightarrow \R$
by $p_\lambda + g_b$.
The action spectrum of $p_\lambda$ forms a relatively
compact subset of $\R$. Hence
there exists a function $a(\lambda)$
such that all $1$-periodic orbits of $p_\lambda$ have action
greater than or equal to $a(\lambda)$.
We choose $b > |a(\lambda) / \delta|+1$.
This ensures that all the $1$-periodic orbits of $K_\lambda^b$
whose action is less than $a(\lambda)$ lie inside
$(0,\infty) \times S^1 \times (\widehat{F'} \setminus F)$.
The Hamiltonian $K_\lambda^b$
splits as a sum $p_\lambda + g_b$, so the
$1$-periodic orbits in this region are
pairs $(o_1,o_2)$ where $o_1$
is an orbit of $p_\lambda : (0,\infty) \times S^1 \rightarrow \R$
and $o_2$ is an orbit of $g_b|_{\widehat{F'} \setminus F}$.
All the orbits of $g_b|_{\widehat{F'} \setminus F}$
are critical points of this function.

We have a natural map
$q : (0,\infty) \times M_{\phi'} \twoheadrightarrow S^1$
defined as the projection to $M_{\phi'}$ composed with
$\pi'$.
Let $\beta_k \subset H_1((0,\infty) \times M_{\phi'})$
be subset represented by loops which project
under $q$ to loops which wrap around $S^1$ $k$ times.
The Hamiltonian $p_\lambda : (0,\infty) \times S^1 \rightarrow \R$
(after perturbing it using  work from \cite{CieliebakFloerHoferWysocki:SymhomIIApplications})
has exactly two orbits which wrap around $k$ times for $\lambda$ large enough.
Hence the number of $1$-periodic orbits of action less
than $a(\lambda)$ representing the class $\beta_k$
is $\kappa$.
We have that
$SH_*^{\beta_k}(K_\lambda^b,J)$ and
$SH_*^{(a(\lambda),\infty),\beta_k}(K_\lambda^b,J)$
are independent of $b$ for $b < a(\lambda)$.
The reason why this is true is because
$SH_*^{\beta_k}(K_\lambda^b,J)$
only depends on the slope $\lambda$ and not on $b$.
Also all the orbits of the Hamiltonian $K_\lambda^b$
of action greater than $a(\lambda)$ sit inside
$(0,\infty)\times M_\phi$ and this Hamiltonian only varies
with respect to $b$ outside $(0,\infty)\times M_\phi$.
A maximum principle (for an appropriate $J$) then ensures that no Floer trajectories
leave $(0,\infty)\times M_\phi$ which means that
$SH_*^{(a(\lambda),\infty),\beta_k}(K_\lambda^b,J)$
does not depend on $b$.
There is a long exact sequence between the groups
$SH_*^{\beta_k}(K_\lambda^b,J)$,
$SH_*^{(a(\lambda),\infty),\beta_k}(K_\lambda^b,J)$
and $SH_*^{(-\infty,a(\lambda)],\beta_k}(K_\lambda^b,J)$.
This ensures that $SH_*^{(-\infty,a(\lambda)],\beta_k}(K_\lambda^b,J)$
does not depend on $b$.
The rank of $SH_*^{(-\infty,a(\lambda)],\beta_k}(K_\lambda^b,J)$
is less than or equal to $\kappa$ because it is generated
by orbits in the region  $(0,\infty) \times S^1 \times (\widehat{F'} \setminus F)$.
Hence using the above long exact sequence we have
\[|\text{rank}(SH_*^{\beta_k}(K_\lambda^b,J)) - \text{rank}(SH_*^{(a(\lambda),\infty),\beta_k}(K_\lambda^b,J))| \leq \kappa.\]
Any Floer trajectory connecting orbits of action greater than $a(\lambda)$
is contained in $(0,\infty)\times M_\phi \subset (0,\infty) \times M_{\phi'}$.
Using this fact and a similar one for continuation maps we have
$HF_*(\phi,k) = \varinjlim_\lambda SH_*^{(a(\lambda),\infty),\beta_k}(K_\lambda^b,J))$.
We also have
$HF_*(\phi',k) = \varinjlim_\lambda SH_*^{\beta_k}(K_\lambda^b,J))$.
Hence
\[\left| \text{rank}(HF_*(\phi,k)) - \text{rank}(HF_*(\phi',k)) \right| \leq \kappa.\]
\qed

In Section \ref{section:lefschetzfibrations},
we defined a filtered directed system $(V^\phi_x)$.
We recall the definition here:
we define $V^\phi_x := H^{n-*}(E) \oplus \bigoplus_{i=1}^{\lfloor x \rfloor} HF_*(\phi,k)$.
The morphism between $V^\phi_x$ and $V^\phi_y$
for $x \leq y$ is the natural inclusion map.
Similarly we have a filtered directed system $(V^{\phi'}_x)$.
The previous Lemma tells us that
$\Gamma{(V^{\phi'})} \leq \text{max}( \Gamma{(V^\phi_x)},1)$.
We also have the inequality
$\Gamma{(V^\phi)} \leq \text{max}( \Gamma{(V^{\phi'}_x)},1)$
so these growth rates are equal if they are both greater than or equal to $1$.

\proof of Theorem  \ref{theorem:subcriticalhandleattachaffinevariety}.
By Lemma \ref{lemma:lefschetzsubcriticalhandleattachment}
and Theorem \ref{theorem:affinelefschetzupperbound}
there exists partial Lefschetz fibrations
$\pi : E \rightarrow \C$ and $\pi' : E' \rightarrow \C$
with fibers $F$ and $F'$ respectively
with the following properties:
\begin{enumerate}
\item
$\widehat{E}$ (resp. $\widehat{E'}$) is symplectomorphic to
$\widehat{\overline{A}}$
(resp. $\widehat{N'}$).
\item There is an exact symplectic embedding of $E$ into $E'$
so that $\pi'|_{E} = \pi$. Hence $F$ is an exact submanifold
of $F'$ as well.
\item
There is a neighbourhood ${\mathcal N}$ of the closure of $E' \setminus E$
diffeomorphic to
$\text{nhd}(F' \setminus F) \times \C$
where $\pi'$ is the natural projection map to $\C$.
The Liouville form $\theta_E$ is a product $\theta_{F'} + r^2 d\vartheta$ in
${\mathcal N}$.
\item
$\Gamma{(V^\phi_x)} \leq \text{dim}_\C A$ where $\phi$ is the monodromy map
of the partial Lefschetz fibration $\pi$.
\end{enumerate}
Let $\phi' : \widehat{F} \rightarrow \widehat{F}$ be the monodromy map
for $\pi'$.
We have that $\phi'$ is the identity outside $F$ and is equal to $\phi$
when restricted to $F$.
By the statement after Lemma \ref{lemma:trivialextensionoffloerhomology},
we have that $\Gamma{(V^{\phi'})} \leq \text{max}( \Gamma{(V^\phi_x)},1)$.
Hence $\Gamma(V^{\phi'}) \leq \text{dim}_\C A$.
The completion $\widehat{M}$ admits a partial Lefschetz fibration
with monodromy $\phi'$.
The reason for this is because $\widehat{N}'$ admits such a partial Lefschetz fibration
and there exists relatively compact open subsets $K_1 \subset \widehat{M}$, $K_2 \subset \widehat{N'}$
such that $\widehat{N'} \setminus K_2$ is symplectomorphic to $\widehat{N} \setminus K_1$.
We also have by Theorem \ref{theorem:lefschetzgrowthratebound}, that
$\Gamma(\phi') \geq \Gamma(M)$.
This implies that $\Gamma(M) \leq \text{dim}_\C A$.
There is a subtlety which is that $\Gamma(M)$ depends on the choice of $b \in H^2(M,\Z  / 2\Z)$.
This naturally induces a choice $b' \in H^2(\widehat{F}')$
for $\phi'$ as stated near the start of Section \ref{section:lefschetzfibrationsbounds}.
A choice of $b'' \in H^2(\overline{N}',\Z / 2\Z)$ also
induces a possibly different choice of $b'$ for $\phi'$
which in theory could give a different growth rate for $\phi'$.
But actually there is a choice of $b''$ which ensures the two possible values for
$b'$ are the same because
$H^2(N',\Z / 2\Z) \rightarrow H^2(\partial N',\Z / 2\Z)$ is surjective by assumption.
Hence $\Gamma(V^{\phi'}) \leq \text{dim}_\C A$
for any choice of $b' \in H^2(\widehat{F}')$ invariant under $H^2(\phi')$.
\qed

\section{Attaching subcritical handles} \label{section:attachingsubcriticalhandles}

We will prove Theorem \ref{theorem:subcriticalhandlelefschetzgrowth}.
Here is a Statement:
{\it
Let $M$ be a Liouville domain whose boundary supports
an open book whose pages are homotopic to $CW$-complexes of dimension less
than half the dimension of $M$ and let $M'$ be a Liouville domain
with a subcritical handle attached.
Then $\Gamma(M) = \Gamma(M')$.
}

This is a consequence of this following Theorem
combined with \cite[Lemma 3.1]{McLean:affinegrowth}:
\begin{theorem} \label{theorem:subcriticalfiltereddirectedsystem}
Let $M$, $M'$ be Liouville domains as in the statement of Theorem
\ref{theorem:subcriticalhandlelefschetzgrowth} above,
then
$(SH_*^{\#}(M,\theta_M,\lambda))$ is isomorphic to
$(SH_*^{\#}(M',\theta_{M'},\lambda))$
as a filtered directed system.
\end{theorem}

\proof of Theorem \ref{theorem:subcriticalfiltereddirectedsystem}.
By Lemma \ref{lemma:lefschetzsubcriticalhandleattachment}
and Lemma \ref{lemma:partiallefschetzfibrationexistence} 
there exists partial Lefschetz fibrations
$\pi : E \rightarrow \C$ and $\pi' : E' \rightarrow \C$
with fiber $F$ and $F'$ respectively
with the following properties:
\begin{enumerate}
\item
$\widehat{E}$ (resp. $\widehat{E'}$) is symplectomorphic to $\widehat{M}$
(resp. $\widehat{M'}$).
\item There is an exact symplectic embedding of $E$ into $E'$
so that $\pi'|_{E} = \pi$. Hence $F$ is an exact submanifold of $F'$ as well.
\item \label{item:trivialization}
There is a neighbourhood $N$ of the closure of $E' \setminus E$
diffeomorphic to
$\C \times \text{nhd}(F' \setminus F)$
where $\pi'$ is the natural projection map to $\C$.
The Liouville form $\theta_E$ is a product $r^2 d\vartheta + \theta_{F'}$ in this region.
\end{enumerate}

We can also assume without loss of generality that $0$
is a regular value of $\pi$ and $\pi'$ and that if $\text{bad}(\pi')$
is the region where $\pi'$ is ill defined then there is a small neighbourhood
$N$ of it so that $0$ and $z$ can be connected by a path inside
 \[\C \setminus \pi' \left (N \setminus \text{bad}(\pi')\right) \]
for any $z$ with $|z|$ sufficiently large.
The reason why we need this condition is that
we wish that $d\vartheta$
can be extended to a closed $1$-form on
$\widehat{E}' \setminus {\pi'}^{-1}(0)$ so that we can use a maximum principle
(Lemma \ref{lemma:lefschetzmaximumprinciple2}).

By \cite[Lemma 3.1]{McLean:spectralsequence},
we can assume that after deforming $\theta_{E'}$ slightly that
we have a trivialization $\D_\epsilon \times {\pi'}^{-1}(0)$ of $\pi$ around $0$
such that $\theta_{E'}$ splits as a product $r^2 d\vartheta + \theta_{F'}$
where $(r,\vartheta)$ are polar coordinates on the base $\C$.
Here $\D_\epsilon$ is a small disk of radius $\epsilon$ inside $\C$.
This trivialization matches the trivialization 
$\C \times \text{nhd}(F' \setminus F)$ mentioned above in the regions where they overlap.
By \cite[Lemma 3]{SS:rama}, we have an exact embedding
of $\widehat{F}$ into $\widehat{F'}$
extending the embedding of $F$ into $F'$.
Hence we can embed $\C \times [1,\infty) \times \partial F$ into $\C \times \overline{\widehat{F'} \setminus F}$
which in turn gives us an embedding of $\widehat{E}$ into $\widehat{E'}$.
This embedding has the property that $\pi'|_{\widehat{E}} = \pi$.
We view the cylindrical coordinate $r_F$ of $\widehat{F}$
as a function defined on $\C \times [1,\infty) \times \partial F \subset \widehat{E'}$
parameterizing $[1,\infty)$.
We also have a function $r_{F'}$
defined on $\C \times [1,\infty) \times \partial F' \subset \widehat{E'}$.
Let $\delta>0$ be smaller than the length of the smallest Reeb
orbits of both $\partial F$ and $\partial F'$.
Let $f_b : \widehat{E'} \rightarrow \R$ be a function such that:
$f_b = 0$ inside $E$, and $f = \delta r_F$ in the region $\{2 \leq r_F \leq b\}$,
and $f_b = \delta r_{F'} + \delta b$ for $r_{F'}$ large enough.
We assume that the derivatives of $f_b$ are so small everywhere else that
the only $1$-periodic orbits are critical points of $f_b$.
We also assume that $f_b$ has only finitely many non-degenerate
fixed points outside $E \subset \widehat{E'}$,
and that the value of $f_b$ in this region is greater than or equal to $\delta b$.
We define $K_\lambda$ to be equal to $(\lambda r^2) \circ \pi'$
in the region where $\pi'$ is well defined and anything else away from it.
We also need that $K_\lambda$ is smoothly parameterized by $\lambda$.
We define $H_\lambda^b := K_\lambda + f_b$.
Let $a : (0,\infty) \rightarrow \R$ be a function such that
$a(\lambda)$ is smaller than the action of all the $1$-periodic orbits of
of $K_\lambda$ (this exists because $K_\lambda$ is Lefschetz admissible).
We will assume that $b$ is greater than $a(\lambda) / \delta$.
These functions can be defined so that $\frac{\partial H_\lambda^b}{\partial b} \geq 0$
and $\frac{\partial H_\lambda^b}{\partial \lambda} \geq 0$.
Let $J$ be an almost complex structure compatible with the partial Lefschetz
fibration $\widehat{E'}$ such that in the region
$\{1 \leq r_F \leq b\}$ it splits up as a product $j \oplus J_F$
where $j$ is the complex structure on $\C$ and $J_F$ is an almost complex structure
on $[1,\infty) \times \partial F$ which is cylindrical.
We define $\bar{H}_\lambda : \widehat{E} \rightarrow \R$
to be equal to $H_\lambda^b$ in the region $\{r_F \leq 3\}$.
We then extend $\bar{H}_\lambda$ by the function $\delta r_F$.
We can define $H_\lambda^b$ so that in the region
$\{r_F \leq 3\}$, it does not change when we change $b$.
This means that $\bar{H}_\lambda$ is independent of $b$.
We define $\bar{J}$ to be equal $J$ in the region
$r_F \leq b$ and then extend it over the whole of $\widehat{E}$
so that it is compatible with this partial Lefschetz fibration.
We have a long exact sequence:
\[\rightarrow SH_*^{(-\infty,a(\lambda)]}(H_\lambda^b,J)
\rightarrow SH_*(H_\lambda^b,J)
\rightarrow SH_*^{(a(\lambda),\infty)}(H_\lambda^b,J) \rightarrow.\]
All the orbits of action greater than or equal to
$a(\lambda)$ are inside the region $\{r_F \leq b\}$.
Hence the maximum principle \cite[Lemma 7.2]{SeidelAbouzaid:viterbo},
tells us that $SH_*^{(a(\lambda),\infty)}(H_\lambda^b,J)$
is equal to $SH_*(\bar{H}_\lambda,\bar{J})$.
This also implies that
$SH_*^{(a(\lambda),\infty)}(H_\lambda^b,J)$
is independent of $b$ and so we can view it as a filtered directed system
only depending on $\lambda$.
This isomorphism commutes with the continuation
maps
\[SH_*^{(a(\lambda_1),\infty)}(H_{\lambda_1}^b,J)
\rightarrow 
SH_*^{(a(\lambda_2),\infty)}(H_{\lambda_2}^b,J)\]
and
\[SH_*(\bar{H}_{\lambda_1},\bar{J}) \rightarrow
SH_*(\bar{H}_{\lambda_2},\bar{J}).\]
We also have that $SH_*(H_\lambda^b,J)$
is independent of $b$ because symplectic homology of
these Hamiltonians only depends
on the slope (by using continuation maps
and the maximum principle mentioned above).
The five lemma then tells us that
$SH_*^{(-\infty,a(\lambda)]}(H_{\lambda}^b,J)$
is independent of $b$ and hence is a filtered directed system only depending on
$\lambda$.
The transfer maps between these groups are also independent
of $b$.
All the $1$-periodic orbits generating
$SH_*^{(-\infty,a(\lambda)]}(H_{\lambda}^b,J)$
lie inside the region $D_\epsilon \times {\pi'}^{-1}(0)$.
We can also assume that $H_\lambda^b$ is a product
$\lambda r^2 + f_b|_{ {\pi'}^{-1}(0) }$ in this region and that the almost complex structure
$J$ is a product $j + J_F$ where $j$ is the standard
complex structure on $D_\epsilon \subset \C$
and $J_F$ is some almost complex structure
on ${\pi'}^{-1}(0)$.
This ensures by a maximum principle
(Lemma \ref{lemma:lefschetzmaximumprinciple2}
applied to the closure of the complement of $\D_\epsilon \times \widehat{F}'$)
that any Floer trajectory connecting these orbits stays inside this region.
Hence the filtered directed system
$(SH_*^{(-\infty,a(\lambda)]}(H_{\lambda}^b,J))$
is equal to
$(SH_*(\lambda r^2,j) \otimes SH_*^{(-\infty,a(\lambda)]}(f_b,J_F))$ and the
directed system maps are induced by the ones
in $(SH_*(\lambda r^2))$ and $(SH_*^{(-\infty,a(\lambda)]}(f_b,J_F))$.
Hence the map:
\[SH_*^{(-\infty,a(\lambda_1)]}(H_{\lambda_1}^b,J)
\rightarrow
SH_*^{(-\infty,a(\lambda_2)]}(H_{\lambda_2}^b,J)\]
is zero for $|\lambda_2 - \lambda_1| > 2\pi$.
This implies that the above filtered
directed system is isomorphic to the trivial
one $0$.
Hence by Lemma \ref{lemma:filteredlongexactisomorphism},
we get that
$(SH_*(H_\lambda^b,J))$
is isomorphic to $(SH_*^{(a(\lambda),\infty)}(H_\lambda^b,J))$.
This is in turn isomorphic to
$(SH_*(\bar{H}_\lambda,\bar{J}))$.
Finally by
Lemma \ref{lemma:halflefschetzlefschetzequivalence}, the filtered directed system
$(SH_*(H_\lambda^b,J))$ (resp. $(SH_*(\bar{H}_\lambda,\bar{J}))$)
is isomorphic to
$(SH_*^{\#}(M,\theta_M,\lambda))$  (resp. $(SH_*^{\#}(M',\theta_{M'},\lambda))$).
Hence
$(SH_*^{\#}(M,\theta_M,\lambda))$
and $(SH_*^{\#}(M',\theta_{M'},\lambda))$ are isomorphic as filtered directed systems.
\qed

\section{Appendix A: Partial Lefschetz fibrations} \label{section:handlelefschetz}

\subsection{Relationship with open books} \label{subsection:openbookpartiallefschetz}
Let $N$ be a Liouville domain.
In this section we will prove the following Lemma:
\begin{lemma} \label{lemma:partiallefschetzfibrationexistence}
Suppose $\partial N$ admits an open book supporting the contact structure such that each page
is homotopic to an $n-1$ dimensional $CW$-complex where $n = \frac{1}{2}\text{dim}(N)$.
Then
$\widehat{N}$ is exact symplectomorphic to the completion of a partial Lefschetz fibration.
Also the monodromy map of the open book is the same as the
monodromy map of the associated partial Lefschetz fibration
up to isotopy through symplectomorphisms fixing the boundary of the fiber.
\end{lemma}

From \cite[Theorem 10]{Giroux:openbooks} we have that every contact manifold admits an open book
supporting the contact structure
whose pages are homotopic to an $n-1$-dimensional $CW$ complex.
We will now give a definition of an open book.
Let $F$ be a manifold with boundary and $\phi : F \rightarrow F$ a
diffeomorphism which is the identity near $\partial F$.
From this we can construct a manifold as follows:
Let $M_\phi$ be the mapping torus of $\phi$.
The boundary of $M_\phi$ is $\partial F \times S^1$.
Let $A$ be a manifold obtained by gluing $\partial F \times \D^2$ to $M_\phi$ 
by identifying the boundary $\partial F \times \partial \D^2$ with $\partial M_\phi$.
\begin{defn} \label{defn:openbook}
Such a manifold $A$ is said to have an {\it open book structure} with page $F$.
The submanifold $\partial F \times \{0\} \subset \partial F \times \D^2$
is called the binding $B$. There is a natural fibration
$\pi : A \setminus B \twoheadrightarrow S^1$ with fiber diffeomorphic to the
interior of $F$. The fibers of $\pi$ are called the fibers of the open book.
An open book on $A$ is said to {\it support a contact structure} given by the kernel of a $1$-form $\alpha$
if $\alpha$ is a contact form on the binding, and $d\alpha$ restricted
to each fiber is a symplectic manifold. We also assume that the orientation
of $B$ induced by the contact from $\alpha|_B$ is the same orientation as the
boundary of the symplectic fiber $F$ (note that the closure of a fiber in $A$ is
a manifold with boundary diffeomorphic to $F$).
\end{defn}

\begin{lemma} \label{lemma:standardcontactform}
Let $\alpha$ be a contact form supported by an open book $\pi : A \setminus B \twoheadrightarrow S^1$
as above. Then we can deform $\alpha$ so it is of the form $f(r)\alpha|_B + r^2d\theta$
on the neighbourhood $B \times \D^2$ of $B$
where $(r,\theta)$ are polar coordinates on $\D^2$ and $f$ is a positive function
with derivative $0$ at $0$ and strictly negative derivative elsewhere.
\end{lemma}
\proof of Lemma \ref{lemma:standardcontactform}.
The map $\pi$ restricted to $B \times (\D^2 \setminus \{0\})$ is the projection map from
$B \times (\D^2 \setminus \{0\}) \twoheadrightarrow \partial \D^2 = S^1$.
Let $\alpha' := f(r)\alpha|_B + r^2d\theta$.
This is a contact form if $f(r)$ has sufficiently small derivative.
The contact form $\alpha$
agrees with $\alpha'$ when restricted to $B \times \{0\}$
hence the orientations on $B \times \{x\}$ induced by
$\alpha$ and $\alpha'$ match for $x$ near $0$.
Also the contact forms $\alpha$ and $\alpha'$
induce the same orientations because
they induce the same orientations on the fibers
$\theta = \text{constant}$ and both Reeb vector fields
point in the direction where $\theta$ increases.
Using both these orientation conditions we get that the orientation of $d\alpha$ and $d\alpha'$
agree on $\{b\} \times \D^2$ for every $b \in B$.
We shrink the disk $\D^2$ so that for every $t \in [0,1]$,
we have that 
$t\alpha' + (1-t)\alpha|_{B \times \bracket{x}}$
is a contact form on $B \times \{x\}$ for all $x \in \D^2$.

Let $P_\D : B \times \D^2 \twoheadrightarrow \D^2$
be the natural projection map.
Let $\rho : \D^2 \rightarrow \R$ be a bump function
equal to $0$ near the boundary and $1$ near $0$.
We view $\rho$ now as the function $P_\D^* \rho$.
If I have any $1$-form $\beta$ on $B \times \D^2$
such that $\beta|_{B \times \bracket{x}}$ is a contact
form for all $x \in \D^2$, then pulling back
a $1$-form whose exterior derivative is a sufficiently large volume form on $\D^2$ and
adding it to $\beta$ gives us a contact form.
Also adding any pullback of a $1$-form
to $\alpha$ whose exterior derivative is a non-negative function multiplied by the volume form on $\D^2$
gives us a contact form.
This means that there is a family of $1$-forms $\nu_t$
equal to $\kappa d\vartheta$ ($\kappa > 0$) near $\partial \D^2$, and also
equal to zero near $t=0$ such that
\[P_\D^* \nu_t + t(\rho \alpha' + (1-\rho)\alpha) + (1-t)\alpha\]
is a contact form for all $t \in [0,1]$.

These contact forms are equal to $\alpha + \kappa d\vartheta$
near $B \times \partial \D^2$ hence we can extend them
by $\alpha + \kappa d\vartheta$ over the whole of our contact manifold $A$.
For $\nu_t$ large enough (in the region where $\rho \neq 0$),
we get that $P_\D^* \nu_t + \rho \alpha' + (1-\rho)\alpha$
supports the open book
decomposition.
If we choose $\nu_1$ appropriately then we get a contact form equal to
$f(r)\alpha|_B + g(r)d\theta$  near zero for some function $g$ such that $g'(r)>0$ for $r>0$.
We can deform $g$ near $0$ so that it is equal to $r^2$ near $0$ and
we can ensure that this is a deformation through contact forms.
Hence our new contact form has all the properties we need.
This proves the Lemma.
\qed

We will now prove Lemma \ref{lemma:partiallefschetzfibrationexistence}.
\proof
Let $\pi : \partial N \setminus B \twoheadrightarrow S^1$ be an open book with binding $B$
supporting a contact form $\alpha$ on $\partial N$ whose kernel is $\text{ker}(\theta_N |_{\partial N})$
where $\theta_N$ is the Liouville form on $N$. We assume $\alpha$ has the same coorientation as
$\theta_N|_{\partial_N}$.
We will use the same notation as in the proof of Lemma \ref{lemma:standardcontactform}.
We also assume that the pages are homotopic to $n-1$ dimensional $CW$-complexes
and that near the binding, $\alpha$ has the form as described in the above lemma.
After a Liouville deformation (which doesn't change the symplectomorphism type of $\widehat{N}$),
we can ensure that $\theta_N|_{\partial N} = \alpha$.
When we complete $N$, we extend $\theta_N$ so that is equal to
$r_N \alpha$ on the cylindrical end $\partial N \times [1, \infty)$.
A neighbourhood of the binding is diffeomorphic to
$B \times \D^2$ with a contact form $f(r)\alpha|_B + r^2d\theta$
where $(r,\vartheta)$ are polar coordinates for $\D^2$.
The part of the cylindrical end covering this neighbourhood
is diffeomorphic to $[1,\infty) \times B \times \D^2$ and
$\theta_N = r_N ( f(r)\alpha|_B + r^2d\vartheta )$.

Let $b : [1,\infty) \times B \times \D^2 \rightarrow \R$ be a function with
the following properties:
\begin{enumerate}
\item $b$ is a function of $r$ and $r_N$ only.
\item $b = r_N$ in the region $r \geq \frac{3}{4}$.
\item $b = r \sqrt{r_N}$ in the region $r \leq \frac{1}{2}$.
\item $\frac{\partial b}{\partial r} \geq 0$ and $\frac{\partial b}{\partial r_N} > 0$
in the region $r>0$.
\end{enumerate}
Because the function $b$ is equal to $r_N$ near $ [1,\infty) \times B \times \partial \D^2$,
we can extend $b$ to a function on $[1, \infty) \times \partial N$
by defining it to be $r_N$ outside $ [1,\infty) \times B \times \D^2$.
The function $b$ has no singular points.
Define a map $[1, \infty) \times \partial N \rightarrow \C$ by
$(b,\vartheta)$ where these are polar coordinates in $\C$.
First of all, outside the region $ [1,\infty) \times B \times \D^2$,
we have that the fibers of $b$ are symplectic.
In the region $ [1,\infty) \times B \times \D^2$, we have that
a fiber of $(b,\vartheta)$ is the set $\{b(r_N,r) = \text{const}\} \cap \{\vartheta = \text{const}\}$.
The Liouville form on this fiber is
$r_N f(r) \alpha|_B$.
The exterior derivative of this $1$-form is symplectic
because $\frac{\partial b}{\partial r}$ and $\frac{\partial b}{\partial r_N}$
are never negative and at least one of these derivatives is non-zero.
The point is that if we contract $d(r_N f(r) \alpha|_B)$ by the vector
(tangent to the fiber)
\[-\frac{\partial b}{\partial r_N}\frac{\partial}{\partial r} +
\frac{\partial b}{\partial r} \frac{\partial}{\partial r_N}\]
we get a positive multiple of $\alpha|_B$ inside
$[1,\infty) \times B \times \D^2$ which means that
$d(r_N f(r) \alpha|_B)$ restricted to each fiber is symplectic.
%
We have that the $\omega_N$-dual $X_{\vartheta_N}$
of $\vartheta_N$ is $\frac{\partial}{\partial r_N}$.
Because the $r_N$ derivative of $b(r_N,r)$
is greater than zero we have that it is transverse to the level sets of $b(r_N,r)$
and pointing in the direction where $b$ increases.
Let $a = r_N f(r)$ be a new function.
For $c$ large enough we have that
\[ E :=  (b^{-1}([0,2]) \cap (a^{-1}(-\infty,c]) \cup N \]
has the structure of a partial Lefschetz fibration with map $(b,\vartheta)$ and the vector
field $\frac{\partial}{\partial r_N}$ points outwards along its boundary.
The corner of the Lefschetz fibration is the region where $b = a$.
The reason why this partial Lefschetz fibration is trivial near the horizontal
boundary $a^{-1}(c)$ is because in the region $a^{-1}(c-\epsilon,c]$
for some small $\epsilon>0$,
we have that the Liouville form is
$a \alpha|_B + b^2 d\vartheta$ and so it is trivial at infinity.
Also the fibers of the completion of this partial Lefschetz fibration
are identical to the fibers of the open book (up to rescaling the symplectic form)
which means that they have identical monodromy maps.
Finally because $\frac{\partial}{\partial r_N}$ points outwards along the
boundary of the partial Lefschetz fibration and $r_N(\frac{\partial}{\partial r_N}) > 0$
outside a closed subset of the interior of $E$, we have that the completion
$\widehat{E}$ is exact symplectomorphic to $\widehat{N}$ by
Lemmas \ref{lemma:convexsubmanifolds} and
\cite[Corollary 8.6]{McLean:affinegrowth}.
\qed

\subsection{Attaching subcritical handles to partial Lefschetz fibrations}
\label{subsection:attachinghandlelefschetz}
We will show how partial Lefschetz fibrations change
when we add subcritical handles.
We will first describe handle attaching.
Let $1 \leq k \leq n$. We will describe a Weinstein $k$-handle.
Let $(\R^{2n},\omega_{\text{std}})$ be
the standard symplectic manifold.
We also assume that $z_j = p_j + i q_j$
are the standard complex coordinates form $\C^n$.
We have a Liouville vector field
\[V_k := \frac{1}{2}\sum_{i=1}^{n-k} \left(
q_i \frac{\partial}{\partial q_i} + p_i \frac{\partial}{\partial p_i} \right)
+ \sum_{i=n+1-k}^n \left( 2q_i \frac{\partial}{\partial q_i} -p_i \frac{\partial}{\partial p_i} \right).\]
This has exactly one singularity.
We define
\[\phi := \frac{1}{4} \sum_{i=1}^{n-k} \left( q_i^2 + p_i^2\right)
 + \sum_{i=n+1-k}^n \left(q_i^2-\frac{1}{2}p_i^2  \right).\]
Let $\tilde{F}_+,\tilde{F}_-$ be two embedded hypersurfaces
inside $\C^n$ such that:
\begin{enumerate}
\item outside a compact set
they coincide with $\phi^{-1}(1)$.
\item $V_k$ is transverse to $F_\pm$.
\item For every point $x \in F_+$, $\Phi^t_{V_k}(x)$ tends to infinity.
\item For every point $x \in F_-$, $\Phi^{-t}_{V_k}(x)$ tends to infinity.
\item The complement of $F_+ \cup F_-$ consists of three connected regions.
The only relatively compact region is the one containing $0$.
\end{enumerate}
Let ${\mathcal H}$ be the closure of relatively compact region.
This is a $k$-handle.
If we have two different choices of $F_\pm$, then
the associated $k$-handles are isotopic to each other through $k$-handles.
Such handles exist by work from \cite{Weinstein:contactsurgery}.
We say that $F_- \cap {\mathcal H}$ is the negative boundary
of the handle and $F_+ \cap {\mathcal H}$ the positive boundary.
The Liouville vector field is transverse to both boundaries
and it points inwards on the negative one and outwards on the positive one.
We write $\partial_\pm {\mathcal H}$ for
${\mathcal H} \cap F_\pm$.
The region $\partial_- {\mathcal H}$ is called the attaching region.
It is contactomorphic to the closure of some open subset of
the contact manifold ${\it J}^1 \left(S^{k-1} \times \R^{n-k}\right)$
which contains $S^{k-1} \times \{0\}$.
Here 
\[{\mathcal J}^1 \left(S^{k-1} \times \R^{n-k} \right) =
T^*\left(S^{k-1} \times \R^{n-k} \right) \times \R\]
is the first jet bundle of $S^{k-1}$.
A {\it framed isotropic $k-1$ sphere} on $\partial M$
consists of an open set $U \subset {\it J}^1 \left(S^{k-1} \times \R^{n-k}\right)$
containing the $S^{k-1} \times \{0\}$ and a contactomorphism
from $U$ to an open subset of $\partial M$.
This contactomorphism has to respect the coorientation
of the contact structures given by the $1$-form defining them.
We say that two framed isotropic $k-1$ spheres
$a_0 : U_0 \rightarrow \partial W$, $a_1 : U_1 \rightarrow \partial M$
are {\it isotopic} if there exists an open $U \subset U_0 \cap U_1$
containing the $S^{k-1} \times \{0\}$ and a family of framed isotropic
spheres of the form $b^t : U \rightarrow \partial M$
such that $a_0|_U = b^0$ and $a_1|_U = b^1$.

Let $M$ be a Liouville domain 
and let $a : U \rightarrow \partial M$ be an framed isotropic sphere of dimension $k-1$.
Then we find a $k$ handle whose attaching
region as a contact form is identical to the one on the closure of $U$
(possibly after shrinking $U$) and we can attach
this handle along this region.
This gives us a new Liouville domain $M'$.
If two framed isotropic spheres are isotopic then
the new Liouville domains obtained by attaching
handles along them are isotopic through Liouville domains.

We will define what a {\it handle attaching triple}
or a HAT is for a Liouville domain $M$.
This is a triple $(f,\beta,\gamma)$ where $f : S^{k-1} \rightarrow \partial M$
is a smooth embedding where $k$ is smaller than half the
dimension of $M$.
Also $\beta$ is a normal framing for $f$ inside $\partial M$
(i.e. a bundle isomorphism $\beta : S^{k-1} \times \R^{2n-k} \rightarrow \nu_f$
where $\nu_f$ is the normal bundle to $f$).
Here $\gamma : S^{k-1} \times \C^n \rightarrow f^* TM$
is a symplectic bundle isomorphism where we
give $\C^n$ the standard symplectic structure.
There is an injective bundle homomorphism
$df : TS^{k-1} \hookrightarrow f^*TM$.
Let $\underline{\R}$ be the trivial $\R$ bundle over $S^{k-1}$.
We also have a bundle morphism
$Df : TS^{k-1} \oplus \underline{\R} \hookrightarrow f^*TM$
given by $df + L$ where $L$ sends the positive unit
vector in $\R$ to an inward pointing vector along $\partial M$.
We say that $(f,\beta,\gamma)$ is a handle
attaching triple or HAT if the map $\gamma$ is isotopic
to $Df \oplus \beta$ through real bundle isomorphisms.
Here we view $Df \oplus \beta$ as a bundle
map from the trivial bundle $S^{k-1} \times \R^{2n}$
to $f^* TM$ where we trivialize $(TS^{k-1} \oplus \underline{\R}) \times \R^{2n-k}$
by viewing it as $T\R^{2n}|_{S^{k-1}}$ where $S^{k-1}$ is the unit sphere in $\R^k \subset \R^{2n}$.
An isotopy of HAT's is a smooth family
of HAT's $(f_t,\beta_t,\gamma_t)$.
If we have a framed isotropic sphere $a : U \rightarrow \partial M$,
then we have a HAT $(f,\beta,\gamma)$ as follows:
the map $f$ is given by $a|_{S^{k-1} \times \bracket{0}}$.
We have a canonical isomorphism from
$S^{k-1} \times \R^k$ to the normal bundle
of $S^{k-1}$ inside ${\mathcal J}^1 \left(S^{k-1} \times \R^{n-k}\right)$
because we identify the bundle $T^*S^{k-1} \oplus \underline{\R} \times \R^{n-k}$
with $\R^{2n}$ from the embedding of the unit sphere $S^{k-1}$
into $\R^k \subset \R^k \times \R^{n-k}$.
Hence we have a canonical framing of the normal bundle of $S^{k-1}$
inside ${\mathcal J}^1 \left(S^{k-1} \times \R^{n-k}\right)$.
By pushing forward this framing via the map $a$ we get our framing $\beta$.
A similar argument gives us our framing $\gamma$ if we embed our sphere in $\R^k \subset \C^{n}$.
By an h-principle (see \cite{EliashberMishachev:hprinciple},\cite{Eliashberg:steintopology}),
any HAT is isotopic through HAT's to the HAT associated
to a framed isotropic sphere (as long as $k<n$).
This same h-principle also ensures that if two
framed isotropic spheres are isotopic through HAT's
then they are isotopic through framed isotropic spheres
for $k<n$.
Any HAT gives us a unique Liouville domain $M'$
up to isotopy by attaching a $k$-handle
along some framed isotropic sphere which
is isotopic to this HAT.

From now on we will assume that $k<n$.
Let $(E,\pi)$ be a partial Lefschetz fibration.
The manifold $E$ is a manifold with corners. We can smooth these
corners so we have a Liouville subdomain $M\subset \widehat{E}$
such that $\widehat{M}$ is exact symplectomorphic to $\widehat{E}$.
We start with a HAT $(f,\beta,\gamma)$ modelling a $k$ handle on the contact boundary
of $M$. Let $M'$ be obtained from $M$ by
attaching a Weinstein handle along this HAT.

Let $F$ be a smooth fiber of $E$.
Near the horizontal boundary of $E$, we have that
the fibration looks like a product fibration
$\nhd(\partial F) \times \D \twoheadrightarrow \D$.
Let $(f',\beta',\gamma')$ be a HAT on $F$ modelling a $k$ handle.
We create a new partial Lefschetz fibration by attaching a handle to $F$
along this HAT creating $F'$ and then extending the fibration $E$ by gluing
$\overline{F' \setminus F} \times \D$ to $\text{nhd}(\partial F) \times \D$.


\begin{lemma} \label{lemma:lefschetzsubcriticalhandleattachment}
Given a HAT $(f,\beta,\gamma)$ on $M$ as above, we can find
a corresponding HAT $(f',\beta',\gamma')$ on $F$ such
that $M'$ is isotopic to $E'$ through Liouville domains (after smoothing
the corners of $E'$).
Hence $\widehat{M'}$ is exact symplectomorphic to $\widehat{E'}$.
\end{lemma}

This theorem will be proven in almost exactly the same way as
in \cite{Cieliebak:subcriticalsplit}.
The only difference is that in our case, we have
a partial Lefschetz fibration, whereas \cite{Cieliebak:subcriticalsplit}
has a product Lefschetz fibration.
We will now prove some preliminary Lemmas.

\begin{lemma} \label{lemma:homotopygroupchange}
Let $B$ be a manifold of dimension $m$ and let
$f : B \rightarrow \R$ be an exhausting Morse function
all of whose critical points have index less than $m-k-1$
and such that it only has finitely many critical points.
Let $B'$ be the union of all the stable manifolds of
all the critical points of $f$,
then the map $\pi_{k}(B \setminus B') \rightarrow \pi_{k}(B)$ is an isomorphism
and the map
$\pi_{k+1}(B \setminus B') \rightarrow \pi_{k+1}(B)$ is surjective.
\end{lemma}
\proof
First of all because there are only finitely many critical
points, $B$ is homotopic to $B_C := f^{-1}(-\infty,C]$
for some $C \gg 0$.
The function $-f$ has critical points of index greater
than $k+1$ hence $B_C$ is homotopic to
$f^{-1}(C)$ with cells of dimension greater
than $k+1$ attached.
Attaching a cell of dimension greater than $k+1$
does not change $\pi_k$.
We have that $B'$ is homotopic to $B_C$,
hence $\pi_{k}(B \setminus B') \rightarrow \pi_{k}(B)$ is an isomorphism.
Also attaching these cells can only add relations to $\pi_{k+1}$
which means that the map
$\pi_{k+1}(B \setminus B') \rightarrow \pi_{k+1}(B)$ is surjective.
\qed

Now we prove a preliminary Lemma which is very
similar to \cite[Lemma 2.1]{Cieliebak:handleattach}.
The boundary of $M$ is a smoothing of the boundary of $E$ inside $\widehat{E}$.
We can choose this smoothing so that $\partial M = \partial E$ outside
a small neighbourhood of the corners of $E$. In particular,
we can assume that $\D_{\frac{1}{2}} \times \partial F \subset \partial M \cap \partial_h E$
(Here $\D_{\frac{1}{2}}$ is the disc of radius $\frac{1}{2}$).
\begin{lemma} \label{lemma:isotopyhat}
There exists an embedding 
\[f_0 : S^{k-1} \rightarrow \{0\} \times \partial F
\subset \D_{\frac{1}{2}} \times \partial F \subset \partial M\]
which is homotopic within $\partial M$ to $f$.
\end{lemma}
\proof
We assume the dimension of $E$ is $2n > 2$.
The boundary of $E$ can written as the union of two
manifolds with boundary, one is the horizontal boundary
$\partial_h E = \D \times \partial F$ and the other is the
vertical boundary $\partial_v E = \pi^{-1}(\partial \D)$.
Here $k$ is less than $n$.
Even though $\partial E$ is a manifold with corners, it
is still homeomorphic to $\partial M$. 
All we need to show is that given a $k-1$ sphere in $\partial E$,
we can homotope it to a sphere in $\{0\} \times \partial F \subset \partial_h E$.
This is equivalent to showing that any $k-1$ sphere is homotopic to a $k-1$ sphere
in $\partial_h E$.
If we can show that the relative homotopy group
$\pi_{j}(\partial_v E,\partial_h E \cap \partial_v E)$ is zero for all $j \leq k-1$
then we are done.
The reason for this is that if we have a map $g : S^{k-1} \rightarrow \partial E$
then after generically perturbing $g$,
$g^{-1}(\partial_v E)$ is a codimension $0$ submanifold $Q$ with boundary
and so $(Q,\partial Q)$ is homotopic to a $(k-1,k-2)$ dimensional cell complex.
Hence $g : (Q,\partial Q) \rightarrow (\partial_v E, \partial_v E \cap \partial_h E)$
is homotopic through maps
$(Q,\partial Q) \rightarrow (\partial_v E, \partial_h E \cap \partial_v E)$ to a map
$(Q,\partial Q) \rightarrow (\partial_h E \cap \partial_v E, \partial_h E \cap \partial_v E)$.
This implies that $g$ is homotopic to a map into $\partial_h E$.
We can show
$\pi_{j}(\partial_v E,\partial_h E \cap \partial_v E) = 0$ for $j \leq k-1$ by proving that
$\pi_{j}(\partial_h E \cap \partial_v E) \rightarrow \pi_{j}(\partial_v E)$ is an isomorphism
for $j < n-2$ and a surjection for $j = n-2$ by a long exact sequence argument.
Because $\partial_v E$ is a fibration with fiber $F$ we have the following commutative diagram:
\[
\xy
(-3,0)*{}="A"; (25,0)*{}="B";(50,0)*{}="C"; (70,0)*{}="D";(100,0)*{}="E"; 
(-3,-10)*{}="G"; (25,-10)*{}="H";(50,-10)*{}="I"; (70,-10)*{}="J";(100,-10)*{}="K";
"A" *{\pi_{j}(\partial F)};
"B" *{\pi_{j}(\partial_v E \cap \partial_h E)};
"C" *{\pi_{j}(\partial S)};
"D" *{\pi_{j-1}(\partial F)};
"E" *{\pi_{j-1}(\partial_v E \cap \partial_h E)};
"G" *{\pi_{j}(F)};
"H" *{\pi_{j}(\partial_v E)};
"I" *{\pi_{j}(\partial S)};
"J" *{\pi_{j-1}(F)};
"K" *{\pi_{j-1}(\partial_v E).};
%
%
{\ar@{->} "A"+(7,0)*{};"B"-(14,0)*{}};
{\ar@{->} "B"+(14,0)*{};"C"-(7,0)*{}};
{\ar@{->} "C"+(7,0)*{};"D"-(10,0)*{}};
{\ar@{->} "D"+(9,0)*{};"E"-(16,0)*{}};
{\ar@{->} "G"+(6,0)*{};"H"-(8,0)*{}};
{\ar@{->} "H"+(8,0)*{};"I"-(7,0)*{}};
{\ar@{->} "I"+(7,0)*{};"J"-(9,0)*{}};
{\ar@{->} "J"+(8 ,0)*{};"K"-(10,0)*{}};
{\ar@{->} "A"+(0,-3)*{};"G"+(0,3)*{}};
{\ar@{->} "B"+(0,-3)*{};"H"+(0,3)*{}};
{\ar@{->} "C"+(0,-3)*{};"I"+(0,3)*{}};
{\ar@{->} "D"+(0,-3)*{};"J"+(0,3)*{}};
{\ar@{->} "E"+(0,-3)*{};"K"+(0,3)*{}};
"C"+(3,-5) *{\cong};
\endxy
\]
The horizontal arrows form long exact sequences coming from the fibration
and the vertical arrows are induced by the natural inclusion maps.
The morphism $\pi_j (\partial F) \rightarrow \pi_j(F)$
is an isomorphism for all $j < n-2$ and a surjection for $j = n-2$
by Lemma \ref{lemma:homotopygroupchange}.
Hence by a repeated application of the five lemma we have that
$\pi_{j}(\partial_h E \cap \partial v_E) \rightarrow \pi_{j}(\partial_v E)$ is an isomorphism for
$j < n-2$ and surjective for $j = n-2$.
This implies that $f$ is homotopic to some
\[f_0 : S^{k-1} \rightarrow \{0\} \times \partial F
\subset \D_{\frac{1}{2}} \times \partial F \subset \partial M.\]
\qed

\proof of Lemma \ref{lemma:lefschetzsubcriticalhandleattachment}.
By Lemma \ref{lemma:isotopyhat}, we can isotope $f$ to $\widehat{f}$ whose image
is contained in \[\{0\} \times \partial F
\subset \D_{\frac{1}{2}} \times \partial F \subset \partial M.\]
By using parallel transport techniques, we then have an isotopy of HAT's
joining $(f,\beta,\gamma)$ to $(\widehat{f},\beta',\gamma')$
for some $\beta',\gamma'$.
By looking at the proof of \cite[Lemma 2.3]{Cieliebak:subcriticalsplit},
we can isotope $\beta',\gamma'$ to $\widehat{\beta},\widehat{\gamma}$
with the following property:
there exists a HAT $(f_0,\beta_0,\gamma_0)$
on $\partial F$ so that
$f_0$ is the same as $\widehat{f}$ where we identify
$\partial F$ with $\{0\} \times \partial F$,
$\widehat{\beta} = \text{id}_\C \times \beta_0$
and $\widehat{\gamma} = \text{id}_\C \times \gamma_0$.
This is on the region $\D_{\frac{1}{2}} \times \partial F \subset \partial M$.

Let $F'$ be obtained from $F$ by attaching a handle along $(f_0,\beta_0,\gamma_0)$.
Then we can glue $\D \times \overline{F' \setminus F}$ to
$E$ along $\partial_h E = \D \times \partial F$ to obtain
a new partial Lefschetz fibration $E'$.
Here the gluing map is $\text{id} \times q$
where $q$ is the map gluing the handle along $(f_0,\beta_0,\gamma_0)$.
We have that $E'$ is a Liouville domain (after smoothing the corners slightly)
that is equal to $E$ with a handle attached to $(\widehat{f},\widehat{\beta},\widehat{\gamma})$.
Let $M'$ be obtained from $M$ by attaching a handle along $(f,\beta,\gamma)$.
Because $(f,\beta,\gamma)$
is isotopic to $(\widehat{f},\widehat{\beta},\widehat{\gamma})$ we get that
$E'$ is isotopic through Liouville domains to $M'$.
Hence $\widehat{E'}$ is exact symplectomorphic to $\widehat{M'}$ by \cite[Lemma 5]{SS:rama}.
This completes the Lemma.
\qed

\section{Appendix B : Algebraic and symplectic Lefschetz fibrations}
\label{section:algebraicsymplecticlefschetzfibration}

The aim of this section is to prove Theorem \ref{theorem:algebraiclefschetzfibration}.
Here is a statement of this theorem:
{\it For every algebraic Lefschetz fibration $p : A \rightarrow \C$,
where $q \in \C$ is a regular value,
there is a Lefschetz fibration $\pi : \widehat{E} \rightarrow \C$
such that $A$ is convex deformation equivalent to $\widehat{E}$
and $A \setminus p^{-1}(q)$ is convex deformation equivalent
to a Lefschetz fibration $\pi' : \widehat{E'} \rightarrow \C^*$
such that $\pi,\pi'$ have two identical Lefschetz cylindrical end components.
Also the fiber $\pi^{-1}(q)$ is convex deformation equivalent to $p^{-1}(q)$.
}

Before we prove this theorem, we need some preliminary
definitions and lemmas.

\begin{defn} \label{defn:partiallytrivializedfibration}
A {\it partially trivialized fibration} is a
manifold $M$, a $1$-form $\theta_M$ and a smooth map
$p_M : M \rightarrow S$ satisfying the following properties:
\begin{enumerate}
\item $d\theta_M$ is a symplectic form.
\item $p_M$ has finitely many singularities and away from
these singularities we have
that $d\theta_M$ restricted to the fibers is symplectic.
These singularities are also modelled on Lefschetz singularities
and are on distinct fibers.
\item All smooth fibers are exact symplectomorphic to some
finite type convex symplectic manifold $(F,\theta_F)$.
\item
We have well defined parallel transport maps between any
two fibers (i.e. the horizontal lift of any integrable vector field
on $S$ is still integrable away from the singularities).
\end{enumerate}
\end{defn}

A {\it deformation of partially trivialized fibrations} is defined to be
a smooth family $(M,\theta^t_M,p_M)$ of such fibrations
such that we have the following additional properties:
\begin{enumerate}
\item
Let $H_t$ be the plane distribution which is $d\theta^t_M$
orthogonal to the fibers of $p_M$ (away from the singularities).
If we take the fibration 
\[\widetilde{p} : M \times [0,1] \rightarrow S \times [0,1]\]
where $\widetilde{p}(x,t) = (p_M(x),t)$ and where the
connection is given by the horizontal plane distribution
$H_t \oplus \R.\frac{\partial}{\partial t}$, then
we require that this connection has well defined parallel transport maps as well.
\item There is an open subset $M_h$
and a smooth family of functions $\kappa_t$
such that for all $q \in S$,  $p_M^{-1}(q) \cap (M \setminus M_h)$
is compact and \[\theta^0_M|_{p_M^{-1}(q) \cap M_h} = \left(\theta^t_M - d\kappa_t\right)|_{p_M^{-1}(q) \cap M_h}.\]
Also $M \setminus M_h$ contains all the singularities of $p_M$.
\end{enumerate}

We say that $p_M$ is
{\it trivial at infinity} if there is an open subset
$M_h \subset M$ such that
\begin{enumerate}
\item the fibers of $p_M|_{M \setminus M_h}$
are compact
\item
There is a trivialization $S \times Q$
of $p_M|_{M_h}$ such that
$\theta_M|_{M_h} = \theta_F|_{Q} + \theta_S$
where $\theta_S$ is a $1$-form on $S$.
Here $\theta_F = \theta_M|_{p_M^{-1}(q)}$ and $Q$ is an open subset of
$p_M^{-1}(q)$ whose complement is compact.
\end{enumerate}
The fiber $p_M^{-1}(q)$ has an exhausting function
$f_F : p_M^{-1}(q) \rightarrow \R$ such that
$df_F(X_{\theta_F}) > 0$ for $f_F$ sufficiently large.
We view $f_F$ as a function on $M$ by pulling it back to the
trivialization $S \times Q$ and then extending it any way
we like inside $M$.
We say that $f_F$ is the {\it vertical cylindrical coordinate}.

\begin{lemma} \label{lemma:trivializedskeleton}
Suppose that our base $S$ is contractible.
Then $(p_M,\theta_M)$ is deformation equivalent to a new
partially trivialized fibration $(p'_M,\theta'_M)$
which is trivial at infinity.
\end{lemma}
\proof of Lemma \ref{lemma:trivializedskeleton}.
We assume that $M$ is connected.
Because $S$ is contractible we have a smooth family of maps
$\iota_t : S \hookrightarrow S$ ($t \in [0,1]$) satisfying
\begin{enumerate}
\item $\iota_0$ is the identity map.
\item $\iota_t$ is an embedding for all $t \in [0,1)$.
\item The image of $\iota_1$ is contained in a point $a \in S$.
\item $\iota_t(a) = a$ for all $t$.
\item For any open set containing $a$ there is a $T < 1$
such that $\iota_t(S)$ is contained in this subset
for $t > T$.
\end{enumerate}
Let $Q \subset p_M^{-1}(a)$ be an open subset such that its complement is compact
and so that if we parallel transport along the path
$\iota_{1-t}(x)$ starting at a point in $Q$ then it does not hit any singular point of $p_M$.
Let \[P^t_x : p_M^{-1}(\iota_1(x)) \cap Q \rightarrow  p_M^{-1}(\iota_{1-t}(x))\]
be the associated parallel transport maps along the path $\iota_{1-t}(x)$.
For $x \in S$ we define $Q_x$ to be
the image of $P^1_x$ and
$\widetilde{Q}$ to be the union $\cup_{x \in S} Q_x$.
We have a natural map $\pi_Q : \widetilde{Q} \rightarrow Q$
sending a point $q \in Q_x$ to $(P^1_x)^{-1}(q)$.
Hence we have a diffeomorphism
\[(\pi_Q, p_M) : \widetilde{Q} \rightarrow Q \times S.\]
The map $p_M$ is now the natural projection from
$Q \times S$ to $S$.
%
The surface $S$ has a natural volume form $V$ so that
if we look at the natural volume form on the horizontal plane distribution
then it is a positive multiple of $\pi^* V$.
We have a smooth family of maps $a_t : Q \times S \rightarrow Q \times S$
given by
\[a_t(x,y) = (x,\iota_{t}(y)).\]
We write $\nu_t$ to be the pullback of $\theta_M$ via this map.
The problem with $\nu_t$ is that $d\nu_t$ may not be a symplectic form for $t = 1$
because $\iota_t$ is not a diffeomorphism onto its image in this case.
Because $\nu_1$ splits up as a product in $Q \times S$, there is a $1$-form
$\theta'_S$ on the base such that the exterior derivative of
$\nu'_t := \nu_t + t p_M^* \theta'_S$ is a symplectic form for all $t$.
Note that $d\theta'_S$ has the same orientation as $V$.
Because $a_t$ is a fiberwise exact symplectomorphism, we have that
$\nu'_t = \theta_F + \beta_t + dR_t$
where $\theta_F = \theta_M|_{p^{-1}(q)}$, $\beta_t$ is a family
of $1$-forms that vanish in the fibers and $R_t$ is a smooth family of functions.
We also have
\[\theta_M = \nu'_0 = \theta_F + \beta_0 + dR_0\]
inside $Q \times S$.
Let $B : p_M^{-1}(q) \rightarrow \R$ be a bump function which is $0$
outside $Q$ and on a neighbourhood of $\partial Q$ and equal to $1$
outside some compact set. We view $B$ as a function on $M$ defining
it to be the pullback of $B$ on $Q \times S$ and zero outside $Q \times S$.
We define $\nu''_t$ to be equal to
\[ \theta_F + B \beta_t + (1- B)\beta_0 + d(BR_t + (1-B)R_0)\]
inside $Q \times S$ and $\theta_M$ outside $Q \times S$.
The problem is that in the region where $dB \neq 0$, we could have
that $d\nu''_t$ is not a symplectic form. But it is a symplectic form when
restricted to the fibers. This means we can find a $1$-form
$\theta''_S$ such that $d\theta''_S$ is a large volume form with the same orientation as $V$
and such that
$\theta^t_M := \nu''_t + t \theta''_S$ is a symplectic form for all $t$.

If we view $\theta^t_M$ as a $1$-form on $[0,1] \times M$ where $t$ parameterizes $[0,1]$
then we wish to show that the associated horizontal plane distribution ${\mathcal H}$
spanned by the horizontal plane distribution of $d\theta^t_M$ inside $\{t\} \times M$
and $\frac{\partial}{\partial t}$ has well defined parallel transport maps.
The base of the fibration is $[0,1] \times S$ and the fibration map is
$(\text{id},p_M)$.
We only need to show that this condition works inside $Q \times S$
(i.e. no paths in this region tangent to ${\mathcal H}$ can escape
to infinity if the projection of this path to $S$ is relatively compact).
Let $p(t) = (f(t),w(t)) \in [0,1] \times S$ be some path in the base $[0,1] \times S$ and $(f(0),q,w(0))$ a point in
\[[0,1] \times Q \times S.\]
We define $(\tilde{q}(t),\iota_{f(t)}(w(t)))$ to be the lift of the path $\iota_{f(t)}(w(t))$ to $Q \times S$ where we use the horizontal
plane distribution associated to $\theta_M$ (we will assume that this lift does not exit this region)
for all $s \in [0,1]$.
We define $\tilde{p}$ to be the path $(f(t),\tilde{q}(t),w(t))$.
This is tangent to ${\mathcal H}$ and a lift of the path $(f(t),w(t))$ in $[0,1] \times S$.
The reason why this is true is because 
\[\left(\frac{d}{dt}\tilde{q}(t),\frac{d}{dt}\left(\iota_{f(t)} w(t)\right)\right)\]
is a lift of the vector $f'(t) (\frac{d}{dt}(\iota_t)) + (\iota_{f(t)})_* \frac{d}{dt}(w(t))$.
But the horizontal lift of $f'(t)(\frac{d}{dt}(\iota_t))$ is $(0,f'(t)(\frac{d}{dt}(\iota_t))$
because of the way we trivialized
$Q \times S$ hence $(\frac{d}{dt}\tilde{q}(t),\frac{d}{dt}w(t))$ is a lift
with respect to $\theta^{f(t)}_M$ of $\frac{d}{dt}(w(t))$
which ensures that $\tilde{p}$ is tangent to ${\mathcal H}$.
Hence parallel transport maps for ${\mathcal H}$ are well defined.
This implies that $(M,\theta^t_M,p_M)$ is a deformation of partially trivialized fibrations.
Also $\theta^1_M$ is a product inside $Q \times S$ and so is trivial at infinity.
\qed

We say that that $M$ is an
{\it extremely convex fibration} if it is a
partially trivialized fibration and it satisfies the following
property:
there is a $1$-form $\theta_S$ such that
\begin{enumerate}
\item $(S,\theta_S)$ is a convex symplectic manifold and so is $(M,\theta_M)$.
\item There exists an exhausting function $f > 0$ such that
for every function $g$ with $g > 0$ and $g' \geq 0$
we have that $\theta_M + p_M^* (g \circ f)\theta_S$
has the structure of a convex symplectic manifold
which is convex deformation equivalent to $(M, \theta_M)$.
\item There is an almost complex structure on $M$
compatible with $d\theta_M$ and also one on $S$
compatible with $d\theta_S$ making $p_M$ holomorphic.
\end{enumerate}

\begin{lemma} \label{lemma:trivialfibrationisextremelyconvex}
Suppose that we have a partially trivialized fibration
\[p_M : M \rightarrow S\]
that is trivial at infinity.
Let $\theta_S$ be a $1$-form making $S$ into a finite type
convex symplectic manifold and let $f_S$
be an exhausting function such that
$df_S(X_{\theta_S}) > 0$ outside some compact set.
Then there is a $1$-form $\theta'$ on the base $S$
making $S$ into a convex symplectic manifold such that
$(M,\theta_M + p_M^* \theta',p_M)$ is an extremely convex fibration
and such that $d\theta'$ is a positive multiple of $d\theta_S$.

If $f_F$ is the vertical cylindrical coordinate coordinate of this fibration
that is trivial at infinity then for $c$ sufficiently large
we have that
${f_F}^{-1}(-\infty,c] \cap {f_S}^{-1}(-\infty,c]$
is a Lefschetz fibration whose completion is convex
deformation equivalent to
$(M,\theta_M + p_M^* \theta',p_M)$.
\end{lemma}
\proof of Lemma \ref{lemma:trivialfibrationisextremelyconvex}.
Let $C$ be a constant so that $df_S(X_{\theta_S}) > 0$
in the region $f_S^{-1}[C,\infty)$.
Note that all of the level sets of $f_S$ in $f_S^{-1}[C,\infty)$ are regular.
Let $\rho : \R \rightarrow \R$ be a function such that $\rho(x) = 1$ for $x < C+1$
and $\rho'(x) > 0$ for $x > C+1$.
Let $R$ be a vector field tangent to the level sets of $f_S$
in $f_S^{-1}[C,\infty)$  such that $\theta_S(R) > 0$.
Then $(d\theta_S)(X_{\theta_S},R) > 0$ in this region.
Hence if we have some function $g : \R \rightarrow (0,\infty)$
such that $g' \geq 0$ and $g'(x) = 0$
inside $f_S^{-1}(-\infty,C)$
then $(g \circ f_S)\theta_S$ is still a convex symplectic structure
on $S$.
This is because 
\[\left(d(g \circ f_S)\theta_S\right)(X_{\theta_S},R) = 
(d(g \circ f_S) \wedge \theta_S)(X_{\theta_S},R) +
(g\circ f_S) d\theta_S(X_{\theta_S},R) \]
\[ = ((g'\circ f_S) df_S(X_{\theta_S})) \theta_S(R) +
 (g\circ f_S)  d\theta_S(X_{\theta_S},R) > 0\]
which implies that $d((g \circ \rho \circ f_S)\theta_S)$ is still a symplectic form.

Let $q$ be a regular value of $p_M$.
We write $F$ as $p_M^{-1}(q)$ and $\theta_F$ as $\theta_M|_{F}$.
Because $p_M$ is trivial at infinity we have an open set
$Q \subset p_M^{-1}(q)$ whose complement is compact
and such that we have a region
$M_h \subset M$ where
\begin{enumerate}
\item it is diffeomorphic to $Q \times S$.
\item $p_M^{-1}(s) \cap (M \setminus M_h)$ is relatively compact for all $s \in S$.
\item $\theta_M = \psi_S + \theta_F$ in $Q \times S$
where $\psi_S$ is a $1$-form on the base.
\end{enumerate}
Because $(F,\theta_F)$ has the structure of a finite type convex symplectic
manifold, we have that
$df_F(X_{\theta_F}) > 0$
in the region $f_F^{-1}[D,\infty)$ for some $D \geq 0$.
%
We have (by abuse of notation) that $f_F$ is a function on $M$
such that it is the pullback of $f_F$ to $Q \times S$ via the natural projection to $Q$.
We can assume that it is less than $D$ outside $Q \times S$.

We define $M_c$ to be the intersection
\[(p_M^* f_S)^{-1}(-\infty,c] \cap f_F^{-1}(-\infty,c]\]
for $c \geq \text{max}(C,D)$.
%
Basically by \cite[Theorem 2.15]{McLean:symhomlef}
we have a function $G : \R \rightarrow (0,\infty)$ with $G' \geq 0$
 such that for all functions $g$ with $g' > 0$
and for all $c \geq \text{max}(C,D)$,
$(M_c,\theta_M + p_M^* ((G \circ \rho \circ f_S)\theta_S + (g \circ \rho \circ f_S) \theta_S))$
is a Liouville domain after smoothing the codimension $2$ corner.
Note that \cite[Theorem 2.15]{McLean:symhomlef}
really requires $G \circ \rho \circ f_S + g \circ \rho \circ f_S$ to be a large constant but we can adjust
it so that we only need $G', g' \geq 0$.
So if we choose $f = \rho \circ f_S$
and $\theta' := (G \circ \rho \circ f_S)\theta_S$ then
$(M,\theta_M + p_M^* \theta')$ has the structure of an extremely convex fibration.

We have that $M_c$ has the structure of a Lefschetz fibration for $c \geq \text{max}(C,D)$.
Also $d(f_F + p_M^* f_S) (X_{\theta_M + p_M^*(\theta')}) > 0$
on the boundary $\partial M_c$ and outside $M_c$.
Hence by \cite[Corollary 8.3]{McLean:affinegrowth} we get that
the completion of the Lefschetz fibration $M_c$ is convex deformation equivalent to $M$.
\qed

\begin{lemma} \label{lemma:surfaceexactvolumeforms}
Suppose we have some surface $S$ and two $1$-forms
$\theta_1,\theta_2$ on $S$ so that $d\theta_1$
and $d\theta_2$ are volume forms with the same orientation.
Suppose that both volume forms give $S$ infinite volume.
Let $F$ be any exhausting function on $S$.
Then there is a third $1$-form $\alpha$
and constants $c_1 < c_2 \cdots$ tending to infinity
so that
\begin{enumerate}
\item $d\alpha$ is a volume form with the same orientation as $d\theta_1$.
\item
$\alpha = \theta_1$ in the region
$F^{-1}[c_{4k},c_{4k+1}]$
\item
$\alpha = \theta_2$ in the region $F^{-1}[c_{4k+2},c_{4k+3}]$
for all $k \in \N$.
\end{enumerate}
\end{lemma}
\proof of Lemma \ref{lemma:surfaceexactvolumeforms}.
Let $\overline{\mathcal A}(c)$ be the maximum of the integrals of $d\theta^0_f$
and $d\theta^1_f$ over $F^{-1}(-\infty,c]$.
And let $\underline{\mathcal A}(c)$ be the minimum of these integrals.
Because these functions are increasing and tend to infinity
we have for each $c$ a new number which we denote by ${\mathcal B}(c)$
so that $\underline{\mathcal A}({\mathcal B}(c)) > \overline{\mathcal A}(c)$.
We define inductively $c_1 := 1$,
$c_{k+1} := \text{max}(c_k + 2, {\mathcal B}(c_k)+1))$.
We construct a volume form $\mu$ on $S$ as follows:
In the region $F^{-1}[c_{4k},c_{4k + 1}]$ we have $\mu = d\theta^0_f$
and $\mu = d\theta^1_f$ in the region $F^{-1}[c_{4k + 2},c_{4k+3}]$.
We construct $\mu$ in the other regions inductively.
Suppose (inductively) we have chosen $\mu$ in the region
$F^{-1}(-\infty,c_{4k+3}]$ so that the integrals of $\mu$ and $d\theta^0$
over $F^{-1}(-\infty,c_{4j + 1}]$ are equal for all $j \leq k$
and similarly the integrals of
$\mu$ and $d\theta^1$
over $F^{-1}(-\infty,c_{4j + 3}]$ are equal for all $j \leq k$.
We will now construct $\mu$ in the region $F^{-1}[c_{4k+3},c_{4(k+1)}]$.
The integral of $\mu$ over $F^{-1}(-\infty,c_{4k+3}]$
is strictly less than the integral of $d\theta^0_f$
over $F^{-1}(-\infty,c_{4(k+1)}]$.
This means we can choose $\mu$
in the region $F^{-1}[c_{4k+3},c_{4(k+1)}]$
so that it is a volume form in this region
and its integral over $F^{-1}(-\infty,c_{4(k+1)}]$
is the same as the integral of $d\theta^0_f$ over this region.
Similarly we can choose $\mu$ in the region
$F^{-1}[c_{4(k+1)+1},c_{4(k+1)+2}]$ to ensure
that the integrals of $\mu$ and $d\theta^1_f$
over $F^{-1}(-\infty,c_{4k + 3}]$ coincide.
Hence we have constructed $\mu$ on the whole of $S$.

Because $S$ is not a compact surface, we have that its second homology
group vanishes so $\mu = d\alpha$ for some $\alpha$.
Because the integrals of $\mu$ and $d\theta^0_f$
over $F^{-1}(-\infty,c_{4k + 1}]$ coincide we have that
$\alpha$ and $\theta^0_f$ differ by an exact $1$-form
in the region $F^{-1}[c_{4k},c_{4k + 1}]$, so we may as well assume
that $\alpha = \theta^0_f$ in this region for all $k$.
Similarly we can assume that $\alpha = \theta^1_f$ in the region
$F^{-1}[c_{4k+2},c_{4k+3}]$ for all $k$.
\qed

\begin{lemma} \label{lemma:trivialfibrationconvexdeformationequivalence}
Suppose that two extremely convex fibrations
$(M,\theta^0_M,p_M)$, 
$(M,\theta^1_M,p_M)$
are deformation equivalent as partially trivialized fibrations.
These extremely convex fibrations have associated $1$-forms
$\theta^i_S$ on the base $S$. We assume that they induce
the same orientation on this base.
Then they are convex deformation equivalent.
\end{lemma}
\proof of Lemma \ref{lemma:trivialfibrationconvexdeformationequivalence}.
Let $(M,\theta^t_M,p_M)$, $t \in [0,1]$ be this deformation.
By the definition of extremely convex fibration there is a function
$f_i$ such that for all functions $g>0$ with $g' \geq 0$
we have that $(M,\theta^i_M + (g \circ f_i) \theta^i_S)$
is a an extremely convex fibration convex deformation equivalent to
$(M,\theta^i_M)$.
So all we need to do is find a convex deformation from
$(M,\theta^0_M + (g \circ f_0)\theta^0_S,p_M)$ to
$(M,\theta^1_M + (g \circ f_1)\theta^1_S,p_M)$ for some $g$.
We will show that for some $g$ they are in fact
exact symplectomorphic and hence by Theorem
\cite[Lemma 8.3]{McLean:symhomlef} they are convex deformation equivalent.

We will construct this exact symplectomorphism in two stages.
In the first stage we will show that
$(M,\theta^0_M + (g \circ f_0)\theta^0_S)$ is exact symplectomorphic to
$(M,\theta^1_M + (g \circ f_0)\theta^0_S)$ for some sufficiently large $g$.
Then we will show that
$(M,\theta^1_M + (g \circ f_0)\theta^0_S)$ is exact symplectomorphic to
$(M,\theta^1_M + (g \circ f_1)\theta^1_S)$.

{\it Step 1}: In this step we will construct the first
exact symplectomorphism.
By definition of deformation we have that there is a region
$M_h$ and a smooth family of functions $\kappa_t$
such that $p_M^{-1}(q) \cap (M \setminus M_h)$ is compact
and 
\[a_t := \frac{d}{dt}\left( \theta^t_M  - d\kappa_t \right)\]
restricted to $p^{-1}(q) \cap M_h$ is zero for all regular values $q$ of $p_M$.
Hence we have that its $d(\theta^t_M + (g \circ f_0)\theta^0_S)$-dual
$X^g_{a_t}$ is tangent to the horizontal plane distribution
$H_t$ associated to $\theta^t_M$ in the region $M_h$.
In order to construct our exact symplectomorphism we need
to show that the vector field $X^g_{a_t}$
is integrable for some $g$.

We have that $X^g_{a_t}$ is tangent to the horizontal plane distribution
in the region $M_h$ and the horizontal plane distribution gives
us an integrable connection. This means that
points cannot flow to infinity in the vertical direction
(i.e. inside $p_M^{-1}(K)$ where $K$ is some compact set).
In order to show that they don't flow to infinity
in the horizontal direction we will put certain bounds on
the function $df_0(X^g_{a_t})$ for $g$
sufficiently large.
Let $\gamma$ be an exhausting function on $M$.
We will define:
\[M^\gamma_c := \gamma^{-1}(-\infty,c],\]
\[M^{f_0}_c := p_M^{-1}f_0^{-1}(-\infty,c],\]
\[S^{f_0}_c := f_0^{-1}(-\infty,c].\]
Suppose inductively we have found a function $g_i : \R \rightarrow \R$ such that
the time $t$ flow with respect to $X^{g_i}_t$ of the compact set $M^\gamma_i$
for $0 \leq t \leq 1$
is contained in the region $M^{f_0}_{c_i}$
for some $c_i$.
We assume that $c_i$ is large enough so that $M^{f_S}_{c_i-1}$ contains all the singular values of $p_M$.
Changing the function $g_i$ outside $(-\infty,c_i]$ does not change the above
property. Hence we will define $g_{i+1}$ to be equal to $g_i$
in the region $(-\infty,c_i]$ and something else outside this region.
There is a family of functions $G_t : M \rightarrow (0,\infty)$
(defined away from the singularities of $p_M$)
such that $d\theta^t_M|_{H_t} = G_t.(p_M)^*d\theta_S|_{H_t}$.
The reason why $G_t$ is positive is because
$p_M$ is holomorphic with respect to some compatible
almost complex structure.
Choose $c_{i+1} > c_i + 1$ so that
$M^{f_0}_{c_{i+1}-2}$ contains
$M^\gamma_{i+1}$.
We define $X_{a_t}$ to be the $d\theta^t_M$-dual
of $a_t$.
We define $\widehat{g}: (-\infty,c_{i+1}) \rightarrow \R$
to be equal to $g_i$ inside $(-\infty,c_i]$ and so that $\widehat{g}$
tends to infinity as we reach $c_{i+1}$.
We can choose $\widehat{g}$ so that $\frac{1}{\widehat{g}}$
extends to a smooth function on the whole of $\R$.
This implies that $\frac{1}{\widehat{g} \circ f_0}$ extends to a smooth function on the
whole of $S$.
We have
\[X^{\widehat{g}}_{a_t} = \frac{\frac{G_t}{\widehat{g} \circ f_0}}{\frac{G_t}{\widehat{g} \circ f_0} + 1} X_{a_t}\]
also extends to a smooth family of smooth vector fields
$\Xi_t$ on the whole of $M$.
Because $\Xi_t$ is zero on the boundary of $M^{f_0}_{c_{i+1}}$ we have that
no flowline can pass through this boundary hence $\Xi_t$
inside $M^{f_0}_{c_{i+1}}$ must be integrable as it is tangent to the horizontal plane distribution $H_t$.
Let $\overline{K}$ be the union over all $t \in [0,1]$ of the time $t$ flow of $M^\gamma_{i+1}$ under $\Xi_t$.
We define $g_{i+1} : S \rightarrow \R$ to be equal to $\widehat{g}$ inside $(-\infty,c_i] \cup f_0(p_M(\overline{K}))$
and anything else we like outside this region.
Hence the time $t$ flow of $M^\gamma_{i+1}$ under $X^{g_{i+1}}_{a_t}$ is contained in $M^{f_0}_{c_{i+1}}$
for all $t \in [0,1]$.
We now define $\widetilde{g}$ to be equal to
$g_i$ in the region $(c_{i-1},c_i]$
for all $i$.
Because $M = \cup_i M^\gamma_i$ and the time $1$ flow
along $X^{\widetilde{g}}_{a_t}$ of $M^\gamma_j$ exists for all $j$,
we get that the vector field
$X^{\widetilde{g}}_{a_t}$ is integrable.
This gives us an exact symplectomorphism from
$(M,\theta^0_M + p_M^* (\widetilde{g} \circ f_0)\theta^0_S)$
to
$(M,\theta^1_M + p_M^* (\widetilde{g} \circ f_0)\theta^0_S)$.
We can perform the above construction so that $(\widetilde{g} \circ f_0)d\theta^0_S$
and $(\widetilde{g} \circ f_1)d\theta^1_S$ gives $S$ infinite volume

{\it Step 2}:
We will define $\theta''_i$ to be equal to $(\widetilde{g} \circ f_i) \theta^i_S$.
In this step we will show that
$(M,\theta^1_M + \theta''_0)$ is exact symplectomorphic to
$(M,\theta^1_M + \theta''_1)$ as long as the integral of
$\theta''_i$ over $S$
is infinite.
Let $F$ be an exhausting function on $S$.
By Lemma \ref{lemma:surfaceexactvolumeforms}, there is a third $1$-form $\alpha$
and constants $c_1 < c_2 \cdots$ tending to infinity
so that
\begin{enumerate}
\item $d\alpha$ is a volume form with the same orientation as $d\theta_1$.
\item
$\alpha = \theta''_0$ in the region
$F^{-1}[c_{4k},c_{4k+1}]$
\item
$\alpha = \theta''_1$ in the region $F^{-1}[c_{4k+2},c_{4k+3}]$
for all $k \in \N$.
\end{enumerate}

Let $\psi_t := (1-t)\theta''_0 + t\alpha$.
Let $\beta_t := \frac{d}{dt}\left( p_M^* \psi_t \right)$,
and let $X_{\beta_t}$ be the $d(\theta^1_M + p_M^* \psi_t)$-dual
of $\beta_t$.
We have that $\beta_t = 0$ in the region
$F^{-1}[c_{4k},c_{4k + 1}]$ and $X_{\beta_t}$
is tangent to the horizontal plane distribution.
This means that $X_{\beta_t}$ must be integrable
as any flowline cannot pass through the region
$F^{-1}[c_{4k},c_{4k + 1}]$ and it cannot travel to infinity
in the fiber direction because the horizontal plane distribution
has well defined parallel transport maps.
This means that $(M,\theta^1_M + p_M^* \theta''_0)$
and $(M,\theta^1_M + p_M^* \alpha)$ are exact symplectomorphic.
Similar reasoning ensures that $(M,\theta^1_M + p_M^* \theta''_1)$
and $(M,\theta^1_M + p_M^* \alpha)$ are exact symplectomorphic.
Hence $(M,\theta^1_M + p_M^* \theta''_0)$ and
$(M,\theta^1_M + p_M^* \theta''_1)$ are exact symplectomorphic.

In Step 1 we already established that
$(M,\theta^0_M + p_M^* (\widetilde{g} \circ f_0)\theta^0_S)$ and
$(M,\theta^1_M + p_M^* (\widetilde{g} \circ f_0)\theta^0_S)$ are exact symplectomorphic
for some $\widetilde{g}$ so that $(S,(\widetilde{g} \circ f_i)d\theta^i_S)$ has infinite volume for $i = 0,1$.
This means that
$(M,\theta^0_M + p_M^* (\widetilde{g} \circ f_0)\theta^0_S)$ and
$(M,\theta^1_M + p_M^* (\widetilde{g} \circ f_1)\theta^1_S)$ are exact symplectomorphic.
Hence
$(M,\theta^0_M)$ is convex deformation equivalent
to $(M,\theta^1_M)$.
\qed

We will put a new symplectic structure on our affine variety $A$
making it into a finite type convex symplectic manifold
which is convex deformation equivalent to the standard one coming from its embedding in $\C^N$.
Let $X$ be a compactification of $A$ by a smooth normal
crossing divisor $D$.
Let $s$ be a section of some ample line bundle $L$
such that $s^{-1}(0) = D$.
Choose some metric $\|.\|$ on $L$ such that
its curvature form $F$ has the property that $iF$ is a positive $(1,1)$ form.
The $1$-form $\theta_A := -d^c \log{\|s\|}$ makes $A$ into a convex
symplectic manifold.
Also by \cite[Section 4b]{Seidel:biasedview} we have that
if we make other choices of ample line bundle $L$, section $s$
and curvature form then we get another convex symplectic manifold
convex deformation equivalent to $(A,\theta_A)$.
By looking at the proof of
\cite[Lemma 2.1]{McLean:affinegrowth} we have that $(A,\theta_A)$
is convex deformation equivalent to $(A,\theta'_A)$ where $\theta'_A$
is equal to $\sum_i r_i^2 d\vartheta_i$ restricted to $A \subset \C^N$.
Here $(r_i,\vartheta_i)$ are polar coordinates for the $i$'th $\C$ factor
in $\C^N$.
From now on we will use the convex symplectic structure $\theta_A$
unless stated otherwise.

\proof of Theorem \ref{theorem:algebraiclefschetzfibration}.
We will first show that $p : A \rightarrow \C$
and hence $p|_{A \setminus p^{-1}(q)}$
has the structure of a partially trivialized fibration.
These fibrations have finitely many singularities
and the fibers are holomorphic away from these singularities.
Because the complex structure on $A$ is compatible with the symplectic
form this implies that $d\theta_A$ restricted to these
fibers is a symplectic form (away from the singularities of $p$).
By the results in
\cite[Section 2]{FukayaSeidelSmith:exactlagragiancotangent}
we have that parallel transport maps are well defined.
Because the base $\C$ is connected this implies that
all the smooth fibers are exact symplectomorphic
to the fiber $p^{-1}(q)$ which has the structure of a
finite type convex symplectic manifold
because $p^{-1}(q)$ has the structure of a smooth affine variety.
Hence both  $p : A \rightarrow \C$
and hence $p|_{A \setminus p^{-1}(q)}$
have the structure of a partially trivialized fibration.

We also need to show that they have the structure of an
extremely convex fibration.
Let $\phi_A$ be the exhausting plurisubharmonic function on $A$
such that $-d_c \phi_A = \theta_A$.
Choose some exhausting plurisubharmonic function on the base
$\C$ then pull it back to $A$ and call it $\psi$.
For any function $g$ with $g'>0$ and $g'' \geq 0$ we have that
$g(\psi)$ is also an exhausting plurisubharmonic function.
Hence $\phi_A + g(\psi)$ is an exhausting plurisubharmonic function and so
$(A, -d_c (\phi_A + g\psi))$ has the structure of a convex
symplectic manifold.
It is convex deformation equivalent to $A$ by work from
\cite{EliahbergGromov:convexsymplecticmanifolds}
(see the proof of Theorem 1.4.A).
This implies that $(A,\theta_A,p)$ has the structure
of an extremely convex fibration.
Similar reasoning ensures that
$(A \setminus p^{-1}(q))$ has the structure of an
extremely convex fibration as we only use the fact that
$p$ is holomorphic and that $A \setminus p^{-1}(q)$
is a Stein manifold.

Because the base is contractible
we have that $p$ is deformation equivalent to
a partially trivialized fibration which is trivial at infinity 
by Lemma \ref{lemma:trivializedskeleton}.
By restricting this deformation to
$A \setminus p^{-1}(q)$ we then get that $A \setminus p^{-1}(q)$
is also deformation equivalent to a partially trivialized
fibration at infinity.
Let $(\widetilde{W},\widetilde{P})$ and
$(\widetilde{W'},\widetilde{P}')$ be these fibrations.
We have that $\widetilde{W}' = \widetilde{W} \setminus \widetilde{P}^{-1}(q)$
and $\widetilde{P}' = \widetilde{P}|_{\widetilde{W}'}$.
By Lemma \ref{lemma:trivialfibrationisextremelyconvex},
we can ensure that both of these
fibrations that are trivial at infinity are deformation
equivalent to a partially trivialized fibration
that is extremely convex and trivial at infinity.
Let $(W,P)$ and $(W',P')$ be these fibrations respectively.
These are the same as the fibrations 
$(\widetilde{W},\widetilde{P})$ and
$(\widetilde{W'},\widetilde{P}')$ respectively
except that their $1$-forms are obtained by adding a pullback
of some $1$-form on the base.
This means that their parallel transport maps are identical.
In particular the parallel transport maps of
$(W',P')$ are the same as those from $(W,P)$ restricted to $W'$.
By Lemma \ref{lemma:trivialfibrationconvexdeformationequivalence}
we get that $A$ is convex deformation equivalent to $W$
and $A \setminus p^{-1}(q)$ is convex deformation equivalent to
$W'$.

Let $(r,\vartheta)$ be polar coordinates on $\C$.
We have that
$(\C,\theta_{\C} := r^2 d\vartheta)$
gives $\C$ a finite type convex symplectic structure
so that $d(r^2)(X_{r^2 d\vartheta}) > 0$ for $r^2 > 0$.
We also have a finite type convex symplectic structure
on $\C^*$ with $1$-form $\theta_{\C^*} := (r^2-1)d\vartheta$ and function
$a(r) := r^2 + \frac{1}{r^2}$ such that $da(X_{\theta_{\C^*}}) > 0$
for $r \neq 1$.
Because $(W',P')$ is a subfibration of $(W,P)$
such that their trivializations at infinity coincide, we can assume that their respective
vertical cylindrical coordinates $f_F$ are identical
(i.e. the vertical cylindrical coordinate on $W'$ is the restriction
of the vertical cylindrical coordinate of $W$ to $W'$).
This means by Lemma \ref{lemma:trivialfibrationisextremelyconvex}
that (for $c > 1$ large enough)
\[\{f_F \leq c\} \cap \{P^* r^2 \leq c\}\]
is a Lefschetz fibration with $1$-form $\theta_W$
whose completion is convex deformation equivalent to $W$
and
\[\{f_F \leq c\} \cap \{ a(r) \leq c + \frac{1}{c}\}
= \{f_F \leq c\} \cap \{{P'}^* r^2 \leq c\} \cap \{{P'}^* \frac{1}{r^2} \leq c\}\]
is a Lefschetz fibration with $1$-form $\theta_{W'}$
for $c$ sufficiently large.
The completions $\widehat{W_c}$ and $\widehat{W'_c}$
have identical Lefschetz cylindrical end components
because $\theta_W$ and $\theta_{W'}$ differ by a $1$-form
that is the pullback of a $1$-form from the base $\C^*$.
\qed

\section{Appendix C: A maximum principle}

Let $\phi : \widehat{F} \rightarrow \widehat{F}$
be a compactly supported exact symplectomorphism,
and let $\pi_\phi : M_\phi \twoheadrightarrow S^1$ be its mapping torus.
Let $\alpha_\phi$ be its contact form.
Let $(Q,d\theta_Q)$ be a symplectic manifold whose boundary is $M_\phi$ such that:
\begin{enumerate}
\item there exists a closed $1$-form $\beta$ on $Q$ such that
$\beta|_{M_\phi} = \pi_\phi^* d\vartheta$ where $\vartheta$ is the angle
coordinate for $S^1$.
\item there is a neighbourhood $[1,1+\epsilon) \times M_\phi$
of $M_\phi$ with $d\theta_Q = dr \wedge d\vartheta + d\alpha_\phi$ where $r$
parameterizes the interval $[1,1+\epsilon)$
and we write $d\vartheta$ by abuse of notation as the pullback of
$d\vartheta$ via the natural map to $S^1$.
\item $\theta_Q - \alpha_\phi$ restricted to each fiber of $M_\phi$ is exact.
\end{enumerate}
Let $J$ be an almost complex structure on $Q$ compatible with the symplectic
form $d\theta_Q$ such that near $M_\phi$,
the natural map
\[(\text{id},\pi_\phi) : [1,1+\epsilon) \times M_\phi \twoheadrightarrow
[1,1+\epsilon) \times S^1\]
is $(J,j)$ holomorphic where $j$ is the standard complex structure
on the annulus $[1,1+\epsilon) \times S^1 \subset \C / 2\pi \Z$.
The region inside the mapping torus near infinity is diffeomorphic to
\[ [1,\infty) \times \partial F \times S^1\]
with $\alpha_\phi = d\vartheta + r_F d\alpha_F$ where $r_F$ parameterizes $[1,\infty)$ and
$\alpha_F$ is the contact form on the boundary of $F$.
Hence near this region inside $Q$ we have that $d\theta_Q = dr \wedge d\vartheta + d(r_F \alpha_F)$.
Let $g : \R \rightarrow \R$ be a function which is $0$ near $1$. We define $g(r_F)$ to be a function defined near $M_\phi$
which is equal to $g(r_F)$ when $r_F$ is well defined and $0$ elsewhere.
Let $H$ be a Hamiltonian such that $H = r + g(r_F)$ on a neighbourhood
of $M_\phi$ and $H > 0$.
\begin{lemma} \label{lemma:lefschetzmaximumprinciple2}
Let $S$ be an oriented surface with boundary and
$\gamma$ a $1$-form with $d\gamma \geq 0$.
Suppose we have map from $S$ into $Q$
such that $\partial S$ maps to $M_\phi$.
If $u$ satisfies:
\[(du - X_H \otimes \gamma)^{0,1} = 0\]
near $M_\phi$
then $S$ is contained in $M_\phi$.
\end{lemma}

Actually we can prove a more general version of this Lemma
where the Hamiltonian $H$ can vary with respect to $S$
(i.e. be a function $H : S \times Q \rightarrow \R$)
as long as it is equal to $r + g(r_F)$ near $M_\phi$.
This is useful if we are considering various continuation map equations.

\proof of Lemma \ref{lemma:lefschetzmaximumprinciple2}.
The proof is very similar to the one in 
\cite[Lemma 7.2]{SeidelAbouzaid:viterbo}.
We suppose for a contradiction that $S$ is not contained in $M_\phi$.
We shift $M_\phi$ slightly maybe to a level set $\{r = 1+ \delta\}$ so that it intersects $S$ transversally and so that
$M_\phi \cap S$ is still non-zero.
This gives us a new symplectic manifold with boundary and function $r$,
which we will now write as $Q$ and $r$ by abuse of notation.
This change also means we need to add some constant to $H$ so that $H = r + g(r_F)$ on a neighbourhood of $M_\phi$ again.
So from now on we assume that $S$ intersects $\partial Q = M_\phi$ transversally.
By adding some multiple of $\beta$
to $\theta_Q$ and some exact $1$-form
we can assume that
$\theta_Q = rd\vartheta + \alpha_\phi$ near $M_\phi$.

Because $d\gamma \geq 0$ we know that
\[\int_{\partial S} \gamma \geq 0.\]
Stokes' theorem ensures that $\int_{\partial S} u^* d\vartheta = \int_{\partial S} \beta = 0$.
Near $M_\phi$, $X_H$ is equal to the horizontal lift $-\tilde{\frac{\partial}{\partial \vartheta}}$
of $-\frac{\partial}{\partial \vartheta}$ plus some vector field $X$ tangent to the fibers of our mapping torus.
This means we have
\[\int_{\partial S} \gamma =
\int_{\partial S} u^*d\vartheta - ( -d\vartheta\left(\tilde{\frac{\partial}{\partial \vartheta}} + X\right) ) \gamma\]
\[=\int_{\partial S} d\vartheta \circ (du - X_H \otimes \gamma)\]
\[=\int_{\partial S} d\vartheta \circ J \circ (du - X_H \otimes \gamma) \circ (-j).\]
Now $d\vartheta \circ J = dr$ 
which means that $d\vartheta(JX_H) = 0$.
Hence
\[\int_{\partial S} \gamma = -\int_{\partial S}(dr \circ du \circ j).\]
If $\xi$ is a vector tangent to the boundary of $S$ which is positively
oriented then $j \xi$ points inwards along $S$.
Because $S$ is transverse to $\partial Q$ we have
$-(dr \circ du \circ j)( \xi) < 0$.
But this means that
\[\int_{\partial S} \gamma < 0\] which is a contradiction.
Hence $S$ must be contained in $M_\phi$.
\qed

\bibliography{references}

\end{document}

%% file: bumpfunctions.pstex_t
\begin{picture}(0,0)%
\includegraphics{bumpfunctions}%
\end{picture}%
\setlength{\unitlength}{2960sp}%
\begingroup\makeatletter\ifx\SetFigFont\undefined%
\gdef\SetFigFont#1#2#3#4#5{%
  \reset@font\fontsize{#1}{#2pt}%
  \fontfamily{#3}\fontseries{#4}\fontshape{#5}%
  \selectfont}%
\fi\endgroup%
\begin{picture}(5577,1753)(1486,-5852)
\put(3901,-5761){\makebox(0,0)[lb]{\smash{{\SetFigFont{12}{14.4}{\rmdefault}{\mddefault}{\updefault}{\color[rgb]{0,0,0}$1-\epsilon_S$}%
}}}}
\put(5776,-4636){\makebox(0,0)[lb]{\smash{{\SetFigFont{10}{12.0}{\rmdefault}{\mddefault}{\updefault}{\color[rgb]{0,0,0}$\nu$}%
}}}}
\put(4726,-4486){\makebox(0,0)[lb]{\smash{{\SetFigFont{10}{12.0}{\rmdefault}{\mddefault}{\updefault}{\color[rgb]{0,0,0}$\mu$}%
}}}}
\put(6076,-5686){\makebox(0,0)[lb]{\smash{{\SetFigFont{12}{14.4}{\rmdefault}{\mddefault}{\updefault}{\color[rgb]{0,0,0}$1$}%
}}}}
\put(1501,-4336){\makebox(0,0)[lb]{\smash{{\SetFigFont{12}{14.4}{\rmdefault}{\mddefault}{\updefault}{\color[rgb]{0,0,0}$1$}%
}}}}
\put(5251,-5686){\makebox(0,0)[lb]{\smash{{\SetFigFont{12}{14.4}{\rmdefault}{\mddefault}{\updefault}{\color[rgb]{0,0,0}$1-\epsilon$}%
}}}}
\end{picture}%

%% file: Lefschetzfunction.pstex_t
\begin{picture}(0,0)%
\includegraphics{Lefschetzfunction}%
\end{picture}%
\setlength{\unitlength}{1184sp}%
\begingroup\makeatletter\ifx\SetFigFont\undefined%
\gdef\SetFigFont#1#2#3#4#5{%
  \reset@font\fontsize{#1}{#2pt}%
  \fontfamily{#3}\fontseries{#4}\fontshape{#5}%
  \selectfont}%
\fi\endgroup%
\begin{picture}(7977,7059)(3586,-7666)
\put(6151,-7561){\makebox(0,0)[lb]{\smash{{\SetFigFont{7}{8.4}{\rmdefault}{\mddefault}{\updefault}{\color[rgb]{0,0,0}$2$}%
}}}}
\put(9676,-7486){\makebox(0,0)[lb]{\smash{{\SetFigFont{6}{7.2}{\rmdefault}{\mddefault}{\updefault}{\color[rgb]{0,0,0}$m$}%
}}}}
\put(10351,-7486){\makebox(0,0)[lb]{\smash{{\SetFigFont{6}{7.2}{\rmdefault}{\mddefault}{\updefault}{\color[rgb]{0,0,0}$m+1$}%
}}}}
\put(6001,-3136){\makebox(0,0)[lb]{\smash{{\SetFigFont{6}{7.2}{\rmdefault}{\mddefault}{\updefault}{\color[rgb]{0,0,0}Slope $\epsilon$.}%
}}}}
\put(8776,-886){\makebox(0,0)[lb]{\smash{{\SetFigFont{6}{7.2}{\rmdefault}{\mddefault}{\updefault}{\color[rgb]{0,0,0}Slope $\lambda$}%
}}}}
\put(3601,-5536){\makebox(0,0)[lb]{\smash{{\SetFigFont{6}{7.2}{\rmdefault}{\mddefault}{\updefault}{\color[rgb]{0,0,0}$1$}%
}}}}
\end{picture}%